\documentclass[smallextended]{svjour3}       

\smartqed  
\usepackage{graphicx}

\usepackage{amssymb}
\usepackage{latexsym}

\usepackage{amsmath,amstext,amsthm,amssymb,bbm,amsbsy}
\usepackage{mathtools}
\usepackage{tabularx}
\usepackage{subcaption}
\usepackage{dsfont}
\usepackage{xcolor}

\usepackage[ruled,linesnumbered]{algorithm2e}

\usepackage[backref=page]{hyperref}
\renewcommand*{\backref}[1]{}
\renewcommand*{\backrefalt}[4]{%
    \ifcase #1 (Not cited.)%
    \or        (Cited on page~#2.)%
    \else      (Cited on pages~#2.)%
    \fi}

\newtheorem{prop}{Proposition}[section]
\newtheorem{lem}{Lemma}[section]
\newtheorem{rem}{Remark}[section]

\DeclareMathOperator*{\argmin}{arg\,min}

\newcommand{\bx} {\boldsymbol{x}}
\newcommand{\by} {\boldsymbol{y}}
\newcommand{\bz} {\boldsymbol{z}}
\newcommand{\bv} {\boldsymbol{v}}
\newcommand{\bu} {\boldsymbol{u}}
\newcommand{\bd} {\boldsymbol{d}}
\newcommand{\bp} {\boldsymbol{p}}
\newcommand{\bq} {\boldsymbol{q}}
\newcommand{\ba} {\boldsymbol{a}}
\newcommand{\bb} {\boldsymbol{b}}

\newcommand{\be} {\boldsymbol{e}}
\newcommand{\balp} {\boldsymbol{\alpha}}

\newcommand{\R} {\mathbb{R}}
\newcommand{\N} {\mathbb{N}}

\newcommand{\prox} {\mathrm{prox}}

\newcommand{\co} {\mathrm{co~}}

\newcommand{\unitsim} {\Delta}
\newcommand{\Rn} {\R^n}

\newcommand{\ind} {I}

\newcommand{\opttraj} {\gamma}
\newcommand{\opttrajhd} {\boldsymbol\opttraj}
\newcommand{\initpos} {p}
\newcommand{\initposhd} {\boldsymbol{p}}

\newcommand{\initcond}{\Phi}
\newcommand{\potentfn}{U}
\newcommand{\valuefn}{V}
\newcommand{\matP}{P}

\begin{document}

\title{Hopf-type representation formulas and efficient algorithms for certain high-dimensional optimal control problems\thanks{This research is supported by  DOE-MMICS SEA-CROGS DE-SC0023191, NSF 1820821, and AFOSR MURI FA9550-20-1-0358. P.C. is supported by the SMART Scholarship, which is funded by USD/R\&E (The Under Secretary of Defense-Research and Engineering), National Defense Education Program (NDEP) / BA-1, Basic Research.
Authors' names are given in last/family name alphabetical order.}}

\author{Paula Chen \and J\'er\^ome Darbon  
        \and Tingwei Meng
}

\institute{Paula Chen \at
        Division of Applied Mathematics, Brown University, Providence, RI, USA  \\
              \email{paula\_chen@brown.edu}           
              \and
        J\'er\^ome Darbon \at
              Division of Applied Mathematics, Brown University, Providence, RI, USA \\
              \email{jerome\_darbon@brown.edu}           
           \and
        Tingwei Meng
        \at
      Department of Mathematics, UCLA, Los Angeles, CA, USA \\
      \email{tingwei@math.ucla.edu}
}

\date{Received: date / Accepted: date}

\maketitle

\begin{abstract}
    Two key challenges in optimal control include efficiently solving high-dimensional problems and handling optimal control problems with state-dependent running costs. In this paper, we consider a class of optimal control problems whose running costs consist of a quadratic on the control variable and a convex, non-negative, piecewise affine function on the state variable. We provide the analytical solution for this class of optimal control problems as well as a Hopf-type representation formula for the corresponding Hamilton-Jacobi partial differential equations. Finally, we propose efficient numerical algorithms based on our Hopf-type representation formula, convex optimization algorithms, and min-plus techniques. We present several high-dimensional numerical examples, which demonstrate that our algorithms overcome the curse of dimensionality. We also describe a field-programmable gate array (FPGA) implementation of our numerical solver whose latency scales linearly in the spatial dimension and that achieves approximately a 40 times speedup compared to a parallelized central processing unit (CPU) implementation. Thus, our numerical results demonstrate the promising performance boosts that FPGAs are able to achieve over CPUs. As such, our proposed methods have the potential to serve as a building block for solving more complicated high-dimensional optimal control problems in real-time.
\end{abstract}

\section{Introduction}
Optimal control problems find applications in many practical problems, including trajectory planning~\cite{Coupechoux2019Optimal,Rucco2018Optimal,Hofer2016Application,Delahaye2014Mathematical,Parzani2017HJB,Lee2021Hopf}, robot manipulator control~\cite{lewis2004robot,Jin2018Robot,Kim2000intelligent,Lin1998optimal,Chen2017Reachability}, and humanoid robot control~\cite{Khoury2013Optimal,Feng2014Optimization,kuindersma2016optimization,Fujiwara2007optimal,fallon2015architecture,denk2001synthesis}. 
We formulate a general continuous finite time horizon optimal control problem mathematically as follows:
\begin{equation}\label{eqt:intro_optctrl}
    \valuefn(\bx,t) = \min\left\{\int_0^t
    \ell(\bx(s), s, \balp(s)) ds
    + \initcond(\bx(0)) \right\},
\end{equation}
where $\bx(\cdot)\colon [0,t]\to \Rn$ is a trajectory satisfying the following backward ordinary differential equation (ODE):
\begin{equation}\label{eqt:intro_optctrl_ode}
    \begin{dcases}
    \dot{\bx}(s) = f(\bx(s),s,\balp(s)) & s\in (0,t),\\
    \bx(t) = \bx.
    \end{dcases}
\end{equation}
In the optimal control problem~\eqref{eqt:intro_optctrl}, the variables $\bx\in\Rn$ and $t\in (0,+\infty)$ denote the terminal position and the time horizon, respectively. 
Let $A$ be the control space, which is a subset of a Euclidean space.
Then, the function $\ell\colon \Rn\times [0,t]\times A\to\R$ is called the running cost, the function $\initcond\colon \Rn\to\R$ is called the initial cost, and the objective function in~\eqref{eqt:intro_optctrl} is called the cost of a control $\balp\colon [0,t]\to A$ and the corresponding trajectory $\bx(\cdot)$.
Under some assumptions, the value function $\valuefn$, as defined in~\eqref{eqt:intro_optctrl}, solves the following Hamilton-Jacobi partial differential equation (HJ PDE):
\begin{equation} \label{eqt: HJ}
\begin{dcases} 
\frac{\partial \valuefn}{\partial t}(\bx,t)+H(\bx,t, \nabla_{\bx}\valuefn(\bx,t)) = 0 & \bx\in\mathbb{R}^{n}, t\in(0,+\infty),\\
\valuefn(\bx,0)=\initcond(\bx) & \bx\in\mathbb{R}^{n},
\end{dcases}
\end{equation}
where the Hamiltonian $H\colon \Rn\times [0,T]\times \Rn\to\R$ is defined using the functions $f$ and $\ell$ in the optimal control problem~\eqref{eqt:intro_optctrl} and the initial condition is given by the initial cost $\initcond$. Solving the optimal control problem \eqref{eqt:intro_optctrl} and solving the corresponding HJ PDE \eqref{eqt: HJ} are intrinsically linked. For example, it is well-known that the optimal control in~\eqref{eqt:intro_optctrl} can be recovered from the spatial gradient $\nabla_{\bx}\valuefn(\bx,t)$ of the viscosity solution $\valuefn$ to the HJ PDE~\eqref{eqt: HJ} (see~\cite{Bardi1997Optimal}, for instance).

An active area of research in optimal control and the study of HJ PDEs is the development of numerical methods for high-dimensional problems.
Many practical engineering applications are formulated in high dimensions. For example, multi-agent path planning problems involve several agents, and each agent has several degrees of freedom, such as positions, velocities, and angles. As a result, the corresponding state spaces for these problems have high dimension (usually greater than five). 
However, the computational complexity of standard grid-based numerical algorithms for solving HJ PDEs, such as ENO \cite{Osher1991High}, WENO \cite{Jiang2000Weighted}, and DG \cite{Hu1999Discontinuous}, 
scales exponentially with respect to the dimension. This exponential scaling in dimension is often referred to as the ``curse of dimensionality"~\cite{bellman1961adaptive}. Due to the curse of dimensionality, these grid-based methods are infeasible for solving high-dimensional problems, e.g., for dimensions greater than five.
Several grid-free methods have been proposed to overcome or mitigate the curse of dimensionality, which include, but are not limited to,
optimization methods~\cite{darbon2015convex,Darbon2016Algorithms,darbon2019decomposition,darbon2021hamilton,yegorov2017perspectives,Lee2021Computationally}, 
max-plus methods~\cite{akian2006max,akian2008max,dower2015maxconference,Fleming2000Max,gaubert2011curse,McEneaney2006maxplus,McEneaney2007COD,mceneaney2008curse,mceneaney2009convergence}, 
tensor decomposition techniques \cite{dolgov2019tensor,horowitz2014linear,todorov2009efficient}, sparse grids \cite{bokanowski2013adaptive,garcke2017suboptimal,kang2017mitigating}, 
polynomial approximation \cite{kalise2019robust,kalise2018polynomial}, 
model order reduction \cite{alla2017error,kunisch2004hjb},
dynamic programming and reinforcement learning \cite{alla2019efficient,bertsekas2019reinforcement,zhou2021actor}, and neural networks \cite{bachouch2018deep,bansal2020deepreach,Djeridane2006Neural,jiang2016using,Han2018Solving,hure2018deep,hure2019some,lambrianides2019new,Niarchos2006Neural,reisinger2019rectified,royo2016recursive,Sirignano2018DGM,Li2020generating,darbon2020overcoming,Darbon2021Neural,darbon2021neuralcontrol,nakamurazimmerer2021adaptive,NakamuraZimmerer2021QRnet,jin2020learning,JIN2020Sympnets,onken2021neural}. 

However, many grid-free methods still rely on approximations. Instead of approximating the solution space by a finite-dimensional space (as grid-based methods do), grid-free methods often approximate the original optimal control problem by some simpler, more easily computable optimal control problems. In doing so, the solution to the original problem is approximated using the solutions to the simpler ones.
Thus, an important research direction is to enlarge the class of optimal control problems with easily computable solutions; such problems and their corresponding exact solvers can then serve as building blocks for solving more complicated optimal control problems. Some well-known techniques for solving optimal control problems that often serve as these building blocks include: the linear-quadratic regulator (LQR)~\cite{Li2004iterative,Sideris2005efficient,McEneaney2006maxplus,Coupechoux2019Optimal}, which corresponds to optimal control problems with certain quadratic running costs and initial costs; the Hopf and Lax-Oleinik representation formulas~\cite{darbon2015convex,darbon2019decomposition,Darbon2016Algorithms,yegorov2017perspectives}, which correspond to optimal control problems whose running costs do not depend on the state variable; and the max-plus (or min-plus) technique~\cite{akian2006max,akian2008max,dower2015maxconference,Fleming2000Max,gaubert2011curse,McEneaney2006maxplus,McEneaney2007COD,mceneaney2008curse,mceneaney2009convergence}, which corresponds to optimal control problems whose running costs or initial costs are the maximum (or minimum) of several simpler functions. However, there are still many more classes of optimal control problems that cannot be solved (exactly) using these techniques. For example, optimal control problems with state-dependent running costs are, in general, difficult to solve without approximations.
To this end,~\cite{Chen2021Lax} recently provided a Lax-Oleinik-type formula and corresponding exact numerical solver for certain optimal control problems with running costs quadratic in the state variable and certain constraints on the control variable. However, to our knowledge, there is no numerically-computable representation formula in the literature for optimal control problems with non-quadratic, state-dependent running costs.

In this paper, we consider a class of optimal control problems whose running costs consist of a quadratic on the control variable and a convex, non-negative, piecewise affine function on the state variable. 
We provide the analytical solution to this class of optimal control problems as well as a Hopf-type representation formula for the corresponding HJ PDEs. Moreover, we show that the analytical solutions to these problems with convex initial costs and certain non-convex initial costs are easily and efficiently computable in high dimensions using convex optimization algorithms and min-plus techniques. As such, the results of this paper enlarge the class of easily computable optimal control problems, and thus, our proposed methods have the potential to serve as a building block for solving more complicated optimal control problems. More specifically, since the running cost is non-smooth with respect to the state variable, our proposed methods could be helpful in solving some non-smooth optimal control problems.

The organization of this paper is as follows. In Section~\ref{sec: theory}, we present the class of optimal control problems and HJ PDEs considered in this paper as well as the analytical solutions of these problems. 
More specifically, we analyze the one-dimensional problems in Section~\ref{sec:HJ1d}, we consider a class of separable high-dimensional problems in Section~\ref{sec:HJhd_separable}, and we provide the Hopf-type representation formula for the general high-dimensional case in Section~\ref{sec:HJhd_general}. In Section~\ref{sec:HJhd_minplus}, we use min-plus techniques to extend the Hopf-type formula from Section~\ref{sec:HJhd_general} to solve the general high-dimensional problem with a certain class of non-convex initial costs.
In Section~\ref{sec:admm}, we propose efficient numerical solvers for these problems and present some high-dimensional numerical results. Quadratic initial costs are considered in Section~\ref{subsec:numerical_quad}, and the corresponding numerical solver serves as a building block for the algorithm in Section~\ref{sec:numerical_convex}, which handles more general convex initial costs. Then, in Section~\ref{sec: ADMM_nonconvex}, we generalize our proposed algorithms from the previous sections to handle the class of non-convex initial costs discussed in Section~\ref{sec:HJhd_minplus}. Several high-dimensional numerical results as well as the computational runtime for each example are provided in each of these subsections to demonstrate the performance of our proposed algorithms. In Section~\ref{sec:fpga}, we present an implementation of the numerical solver from Section~\ref{subsec:numerical_quad} on a field-programmable gate array (FPGA) and some corresponding numerical results, which demonstrate the promising performance boosts FPGAs are able to achieve in comparison to CPUs.
Finally, in Section~\ref{sec: conclusion}, we make some concluding remarks and list some possible future directions.
Some technical lemmas and computations for the proofs and the numerical algorithms are provided in the Appendix.
\def\pzero{p}

\section{Analytical solutions}\label{sec: theory}

In this section, we provide the analytical solution to the following 
optimal control problem:
\begin{equation}\label{eqt: result_optctrl1_hd_general}
    \valuefn(\bx,t) = \inf \left\{\int_0^t \left(\frac{1}{2}\|\dot{\bx}(s)\|_{M^{-1}}^2 - \potentfn(\bx(s))\right) ds + \initcond(\bx(0)) \colon \bx(t) = \bx\right\},
\end{equation}
where $t>0$ is the time horizon, $\bx\in\Rn$ is the terminal position, $\bx(\cdot)\colon [0,t]\to\Rn$ is a locally Lipschitz function, and $\dot{\bx}(s)$ denotes its the derivative at time $s$, which exists at $s\in (0,t)$ almost everywhere.
Here, $M$ is a positive definite matrix with $n$ rows and $n$ columns. The matrix $M$ and its inverse respectively define the norms $\|\cdot\|_M\colon \Rn \to\R$ and $\|\cdot\|_{M^{-1}}\colon \Rn \to\R$ by 
\begin{equation*}
    \|\bv\|_M:= \sqrt{\langle \bv, M\bv\rangle}, \quad 
    \|\bv\|_{M^{-1}}:= \sqrt{\langle \bv, M^{-1}\bv\rangle}
    \quad\forall\, \bv\in\Rn,
\end{equation*}
where $\langle\cdot,\cdot\rangle$ denotes the standard Euclidean inner product in $\Rn$. 
The potential energy is given by $\potentfn\colon\Rn\to(-\infty,0]$, which is a piecewise affine concave function satisfying some assumptions. The initial cost is given by the continuous function $\initcond\colon\Rn\to\R$.
In the remainder of this paper, if not mentioned specifically, we use bold characters to denote high-dimensional vectors in $\Rn$, and we use $x_i$ to denote the $i$-th component of a high-dimensional vector $\bx\in \Rn$. 
To avoid the ambiguity of a trajectory and a vector, we use $\bx(\cdot)$ to denote the trajectory, which is a function of the time variable, and we use $\bx$ to denote the vector.

The corresponding HJ PDE reads:
\begin{equation} \label{eqt: HJhd_1_general}
\begin{dcases} 
\frac{\partial \valuefn}{\partial t}(\bx,t)+\frac{1}{2}\|\nabla_{\bx}\valuefn(\bx,t)\|_M^2 + \potentfn(\bx) = 0 & \bx\in\mathbb{R}^n, t\in(0,+\infty),\\
\valuefn(\bx,0)=\initcond(\bx) & \bx\in\mathbb{R}^n,
\end{dcases}
\end{equation}
where the potential energy $\potentfn$, the matrix $M$, and the initial data $\initcond$ are the corresponding quantities in~\eqref{eqt: result_optctrl1_hd_general}.

In what follows, we provide analytical solutions to the optimal control problem~\eqref{eqt: result_optctrl1_hd_general} and the corresponding HJ PDE~\eqref{eqt: HJhd_1_general} under various assumptions. In Section~\ref{sec:HJ1d}, we solve the one-dimensional problems with convex initial cost $\initcond$. Then, in Section~\ref{sec:HJhd_separable}, we use the functions defined in Section~\ref{sec:HJ1d} to solve the high-dimensional problems, where $M$ is the identity matrix, $\initcond$ is convex, and $\potentfn$ has a specific form. In Section~\ref{sec:HJhd_general}, we solve the high-dimensional problems, where $\initcond$ is convex and $M$ and $\potentfn$ satisfy more general assumptions. Finally, in Section~\ref{sec:HJhd_minplus}, we consider a certain class of non-convex initial costs $\initcond$, and we generalize the representation formulas provided in Sections~\ref{sec:HJ1d},~\ref{sec:HJhd_separable}, and~\ref{sec:HJhd_general} to handle this case using min-plus techniques.

\subsection{One-dimensional case} \label{sec:HJ1d}
In this section, we solve the one-dimensional versions of the problems~\eqref{eqt: result_optctrl1_hd_general} and~\eqref{eqt: HJhd_1_general}.
More specifically, we consider the following one-dimensional optimal control problem:
\begin{equation}\label{eqt: result_optctrl1_1d}
    \valuefn(x,t) = \inf \left\{\int_0^t \left(\frac{(\dot{x}(s))^2}{2} - \potentfn(x(s))\right) ds + \initcond(x(0)) \colon x(t) = x\right\},
\end{equation}
where $\potentfn\colon \R\to(-\infty,0]$ is the 1-homogeneous concave function defined by
\begin{equation}\label{eqt: result_defV}
    \potentfn(x) = \begin{dcases}
    -ax & x\geq 0,\\
    bx & x< 0,
    \end{dcases}
\end{equation}
for some positive constants $a$ and $b$.
The corresponding HJ PDE reads:
\begin{equation} \label{eqt: result_HJ1_1d}
\begin{dcases} 
\frac{\partial \valuefn}{\partial t}(x,t)+\frac{1}{2}(\nabla_{x}\valuefn(x,t))^2 + \potentfn(x) = 0 & x\in\mathbb{R}, t\in(0,+\infty),\\
\valuefn(x,0)=\initcond(x) & x\in\mathbb{R}.
\end{dcases}
\end{equation}

In this section, we provide the analytical solution to the one-dimensional optimal control problem~\eqref{eqt: result_optctrl1_1d}. We also present a Hopf-type representation formula for the viscosity solution to the one-dimensional HJ PDE~\eqref{eqt: result_HJ1_1d}.

\begin{figure}[htbp]
    \centering
    \includegraphics[width=0.7\textwidth]{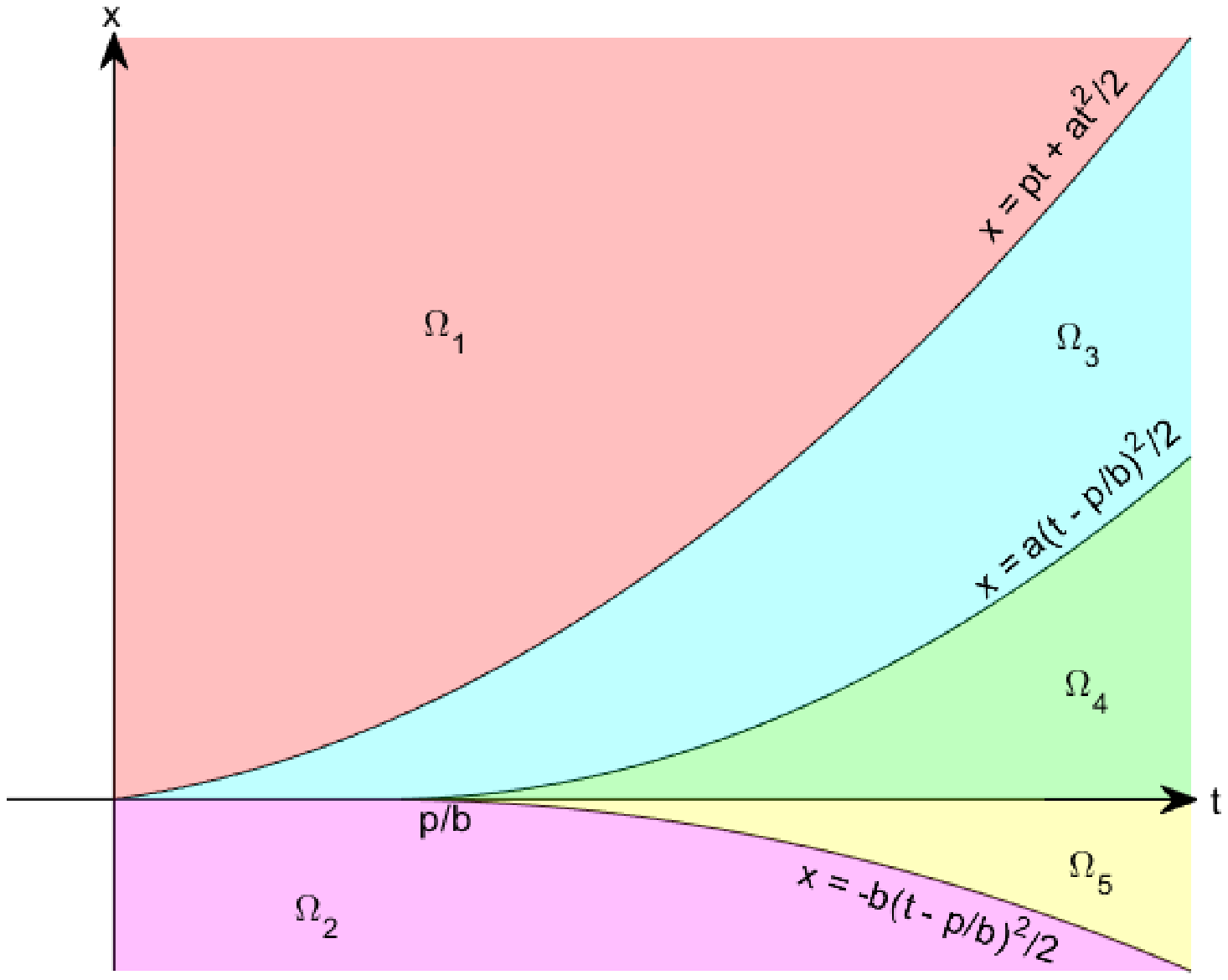}
    \caption{An illustration of a two-dimensional slice of the sets $\Omega_1,\Omega_2,\Omega_3, \Omega_4, \Omega_5$ on the $tx$-plane.}
    \label{fig:plot_of_domain}
\end{figure}

First, we consider the case when the initial cost $\initcond$ is linear, i.e., when $\initcond(x) = \pzero x$ holds for some $\pzero\in\R$. In this case, we denote the solution to the HJ PDE~\eqref{eqt: result_HJ1_1d} by $\R\times [0,+\infty)\ni (x,t)\mapsto \valuefn(x,t; \pzero, a,b)$, and we denote the optimal trajectory in~\eqref{eqt: result_optctrl1_1d} by $[0,t]\ni s\mapsto \opttraj(s;x,t,\pzero,a,b)\in\R$. The function $\valuefn$ is a function of $x$ and $t$, which are, respectively, the terminal position and the time horizon in~\eqref{eqt: result_optctrl1_1d}. The function $\valuefn$ has three parameters: $p\in\R$, which is the slope of the initial cost $\initcond$, and $a,b>0$, which are the two positive parameters in the potential function $\potentfn$. 
The function $\opttraj$ is a function of the current (running) time $s$ with five parameters: $x\in\R$, $t>0$, $\pzero\in\R$, and $a,b>0$, which have the same meaning as the corresponding variables and parameters in the function $\valuefn$.

If $\pzero\geq 0$, we define the function $\R\times [0,+\infty)\ni (x,t)\mapsto \valuefn(x,t; \pzero, a,b)\in\R$ as follows:
\begin{equation} \label{eqt: result_S1_1d}
    \valuefn(x,t; \pzero,a,b) := \begin{dcases}
    f_1(x,t; \pzero,a,b) & (x,t,\pzero)\in\Omega_1,\\
    f_2(x,t; \pzero,a,b) & (x,t,\pzero)\in\Omega_2,\\
    f_3(x,t; \pzero,a,b) & (x,t,\pzero)\in\Omega_3,\\
    f_4(x,t; \pzero,a,b) & (x,t,\pzero)\in\Omega_4,\\
    f_5(x,t; \pzero,a,b) & (x,t,\pzero)\in\Omega_5,
    \end{dcases}
\end{equation}
where the five regions $\{\Omega_i\}_{i=1}^5\subset \R\times [0,+\infty)\times [0,+\infty)$ are defined by
\begin{equation}\label{eqt: domain_subregions}
    \begin{split}
    \Omega_1 &:= \left\{(x,t,\pzero)\in\R\times [0,+\infty)\times [0,+\infty) \colon x\geq \pzero t + \frac{a}{2}t^2\right\},\\
    \Omega_2 &:= \left\{(x,t,\pzero)\colon x< 0, \,t< \frac{\pzero}{b}\right\} \bigcup \left\{(x,t,\pzero)\colon x< -\frac{b}{2}\left(t-\frac{\pzero}{b}\right)^2, \,t\geq \frac{\pzero}{b}\right\},\\
    \Omega_3 &:= \left\{(x,t,\pzero)\colon 0\leq x< \pzero t +\frac{a}{2}t^2, \,t< \frac{\pzero}{b}\right\}\\
    &\quad\quad\quad\quad \bigcup \left\{(x,t,\pzero)\colon \frac{a}{2}\left(t-\frac{\pzero}{b}\right)^2\leq x< \pzero t + \frac{a}{2}t^2,\, t\geq \frac{\pzero}{b}\right\},\\
    \Omega_4 &:= \left\{(x,t,\pzero)\in\R\times [0,+\infty)\times [0,+\infty)\colon 0\leq x< \frac{a}{2}\left(t-\frac{\pzero}{b}\right)^2, \,t\geq \frac{\pzero}{b}\right\},\\
    \Omega_5 &:= \left\{(x,t,\pzero)\in\R\times [0,+\infty)\times [0,+\infty)\colon -\frac{b}{2}\left(t-\frac{\pzero}{b}\right)^2\leq x< 0,\, t\geq \frac{\pzero}{b}\right\},
    \end{split}
\end{equation}
and the five functions $f_1, f_2, f_3, f_4, f_5$ are defined by 
\begin{equation}\label{eqt: case1_deff}
    \begin{split}
        f_1(x,t; \pzero,a,b) &:= -\frac{a^2}{6}t^3 - \frac{a}{2} \pzero t^2 + atx - \frac{1}{2}\pzero^2 t + \pzero x,\\
        f_2(x,t; \pzero,a,b) &:= -\frac{b^2}{6}t^3 + \frac{b}{2}\pzero t^2 - \frac{1}{2}\pzero^2 t + \pzero x - btx,\\
        f_3(x,t; \pzero,a,b) &:= \frac{a+b}{3(a+2b)^2}\left((bt-\pzero)^3 + ((bt-\pzero)^2 + 2x(a+2b))^{3/2}\right) \\
        &\quad\quad 
        - \frac{b}{a+2b}(bt-\pzero)x - \frac{b^2}{6}t^3 + \frac{b}{2}\pzero t^2 - \frac{1}{2}\pzero^2t,\\
        f_4(x,t; \pzero,a,b) &:= \frac{a^2}{3}\left(\frac{2x}{a}\right)^{3/2} - \frac{\pzero^3}{6b},\\
        f_5(x,t; \pzero,a,b) &:= \frac{b^2}{3}\left(\frac{-2x}{b}\right)^{3/2} -\frac{\pzero^3}{6b}.
    \end{split}
\end{equation}
An illustration of a two-dimensional slice of the sets $\Omega_i, i = 1, \dots, 5$ for a fixed $p\geq0$ is shown in Figure~\ref{fig:plot_of_domain}.

Now, we define the function $[0,t]\ni s\mapsto \opttraj(s;x,t,\pzero,a,b)\in\R$ for $\pzero\geq 0$. We define the function $\opttraj$ in five cases, which correspond to the five lines in \eqref{eqt: result_S1_1d}, as follows:
\begin{enumerate}
    \item When $(x,t,\pzero)\in\Omega_1$ holds, which corresponds to the first line in~\eqref{eqt: result_S1_1d}, we define $\opttraj(s;x,t,\pzero,a,b)$ by
    \begin{equation}\label{eqt: optctrl_defx_1}
        \opttraj(s;x,t,\pzero,a,b):= x - \pzero(t-s) - \frac{a}{2}(t^2-s^2) \quad\quad \forall s\in[0,t].
    \end{equation}
    In this case, we have 
    $\opttraj(s;x,t,\pzero,a,b) \geq 0$ for all $s\in[0,t]$. 
    
    \item When $(x,t,\pzero)\in\Omega_2$ holds, which corresponds to the second line in~\eqref{eqt: result_S1_1d}, we define $\opttraj(s;x,t,\pzero,a,b)$ by
    \begin{equation}\label{eqt: optctrl_defx_2}
        \opttraj(s;x,t,\pzero,a,b):= x - \pzero(t-s) + \frac{b}{2}(t^2-s^2) \quad\quad \forall s\in[0,t].
    \end{equation}
    In this case, we have 
    $\opttraj(s;x,t,\pzero,a,b) \leq 0$ for all $s\in[0,t]$. 
    
    \item 
    When $(x,t,\pzero)\in\Omega_3$ holds, which corresponds to the third line in~\eqref{eqt: result_S1_1d}, we define $\opttraj(s;x,t,\pzero,a,b)$ by
    \begin{equation}\label{eqt: optctrl_defx_3}
    \opttraj(s;x,t,\pzero,a,b) :=
    \begin{dcases}
    -\pzero(\tau - s) + \frac{b}{2}(\tau^2 - s^2) & s\in[0,\tau),\\
    (\pzero - b\tau)(s-\tau) + \frac{a}{2}(s-\tau)^2 & s\in[\tau, t],
    \end{dcases}    
    \end{equation}
    where $\tau \in\R$ is defined by
    \begin{equation}\label{eqt:def_opttraj_case3_tau}
        \tau := \frac{(a+b)t+\pzero - \sqrt{(bt-\pzero)^2 + 2(2b+a)x}}{2b+a}.
    \end{equation}
    By straightforward calculation using $0\leq x\leq \pzero t +\frac{a}{2}t^2$, we have $\tau \in [0,t]$, and the function $s\mapsto\opttraj(s;x,t,\pzero,a,b)$ is continuous in this case.
    Here, the trajectory $\opttraj$ is divided into two parts: $\opttraj(s; x,t,\pzero,a,b)\leq 0$ for $s\in [0,\tau)$ and $\opttraj(s; x,t,\pzero,a,b)\geq 0$ for $s\in [\tau,t]$.
    
    \item 
    When $(x,t,\pzero)\in\Omega_4$ holds, which corresponds to the fourth line in~\eqref{eqt: result_S1_1d}, we define $\opttraj(s;x,t,\pzero,a,b)$ by
    \begin{equation}\label{eqt: optctrl_defx_4}
    \opttraj(s;x,t,\pzero,a,b):=
        \begin{dcases}
        -\frac{1}{2b}(\pzero - bs)^2 & s\in\left[0, \frac{\pzero}{b}\right),\\
        0 & s\in\left[\frac{\pzero}{b}, t-\sqrt{\frac{2x}{a}}\right),\\
        \frac{a}{2}\left(s-t+\sqrt{\frac{2x}{a}}\right)^2 & s\in\left[t-\sqrt{\frac{2x}{a}},t\right].
        \end{dcases}
    \end{equation}
    By straightforward calculation using $(x,t,\pzero)\in\Omega_4$, we have $0\leq \frac{p}{b}\leq t-\sqrt{\frac{2x}{a}}\leq t$ and the function $s\mapsto\opttraj(s;x,t,\pzero,a,b)$ is continuous.
    Therefore, in this case, the trajectory $\opttraj$ is divided into three parts: $\opttraj(s; x,t,\pzero,a,b)\leq 0$ in the first time period $s\in[0,\frac{\pzero}{b})$, it remains zero in the second time period $\left[\frac{\pzero}{b}, t-\sqrt{\frac{2x}{a}}\right)$, and it becomes non-negative in the third time period $\left[t-\sqrt{\frac{2x}{a}},t\right]$. 
    
    \item When $(x,t,\pzero)\in\Omega_5$ holds, which corresponds to the fifth line in~\eqref{eqt: result_S1_1d}, we define $\opttraj(s;x,t,\pzero,a,b)$ by
    \begin{equation}\label{eqt: optctrl_defx_5}
    \opttraj(s;x,t,\pzero,a,b):=
        \begin{dcases}
        -\frac{1}{2b}(\pzero - bs)^2 & s\in\left[0, \frac{\pzero}{b}\right),\\
        0 & s\in\left[\frac{\pzero}{b}, t-\sqrt{\frac{2|x|}{b}}\right),\\
        -\frac{b}{2}\left(s-t+\sqrt{\frac{2|x|}{b}}\right)^2 & s\in\left[t-\sqrt{\frac{2|x|}{b}},t\right].
        \end{dcases}
    \end{equation}
    By straightforward calculation using $(x,t,\pzero)\in\Omega_5$, we have $0\leq \frac{p}{b}\leq t-\sqrt{\frac{2|x|}{b}}\leq t$ and the function $s\mapsto\opttraj(s;x,t,\pzero,a,b)$ is continuous.
    Therefore, in this case, the trajectory is divided into three parts: the trajectory $\opttraj$ is negative in the first time period $\left[0, \frac{\pzero}{b}\right)$, it remains zero in the second time period $\left[\frac{\pzero}{b}, t-\sqrt{\frac{2|x|}{b}}\right)$, and it becomes non-positive again in the third time period $\left[t-\sqrt{\frac{2|x|}{b}},t\right]$.
\end{enumerate}

So far, we have provided the analytical solution of the one-dimensional optimal control problem~\eqref{eqt: result_optctrl1_1d} and the corresponding HJ PDE~\eqref{eqt: result_HJ1_1d}, with initial cost $\initcond(x)=\pzero x$ for some $\pzero\geq 0$. 
When the initial cost is $\initcond(x)=\pzero x$ for some $\pzero<0$, we define the function $\R\times [0,+\infty)\ni (x,t)\mapsto \valuefn(x,t; \pzero, a,b)\in\R$ by
\begin{equation}\label{eqt: result_S1_1d_negative}
    \valuefn(x,t; \pzero, a,b) := \valuefn(-x,t; -\pzero, b,a) \quad \forall x\in\R, t\geq 0,
\end{equation}
where $\valuefn(-x,t; -\pzero, b,a)$ on the right-hand side is defined by~\eqref{eqt: result_S1_1d} since $-\pzero$ is positive. 
Similarly, the optimal trajectory $[0,t]\ni s\mapsto \opttraj(s;x,t,\pzero,a,b)\in\R$ for $p<0$ is defined by
\begin{equation}\label{eqt: optctrl_defx_neg}
    \opttraj(s; x,t,\pzero,a,b) := -\opttraj(s; -x,t,-\pzero,b,a) \quad \forall\, s\in[0,t],
\end{equation}
where the right-hand side is well-defined using~\eqref{eqt: optctrl_defx_1},~\eqref{eqt: optctrl_defx_2},~\eqref{eqt: optctrl_defx_3},~\eqref{eqt: optctrl_defx_4}, and~\eqref{eqt: optctrl_defx_5} for different cases.

Now, we consider a general convex initial cost $\initcond\colon\R\to\R$. 
The corresponding HJ PDE is solved using the following Hopf-type formula:
\begin{equation}\label{eqt: result_Hopf1_1d}
    \valuefn(x,t) = \sup_{p\in \R} \{\valuefn(x,t; p, a,b) - \initcond^*(p)\} \quad\forall\, x\in\R,\, t\geq 0,
\end{equation}
where the function $\valuefn(x,t; p, a,b)$ on the right-hand side is defined by~\eqref{eqt: result_S1_1d} and~\eqref{eqt: result_S1_1d_negative}, and the function $\initcond^*$ on the right-hand side is the Legendre-Fenchel transform of the initial cost $\initcond$.
By a one-dimensional corollary of Lemma~\ref{lem: prop_barp_hopf_hd}, the function value $\valuefn(x,t)$ defined in~\eqref{eqt: result_Hopf1_1d} is finite and the maximizer in~\eqref{eqt: result_Hopf1_1d} exists. Moreover, for any positive time horizon $t>0$, the maximizer is unique and we denote the unique maximizer by $p^*(x,t)\in\R$.
Then, the optimal trajectory $[0,t]\ni s\mapsto \opttraj(s;x,t)\in\R$ is defined by 
\begin{equation}\label{eqt: optctrl_traj_J}
    \opttraj(s; x,t) := \opttraj(s; x,t,p^*(x,t), a,b) \quad\forall s\in[0,t],
\end{equation}
where the function $\opttraj(s; x,t,p, a,b)$ on the right-hand side is defined by~\eqref{eqt: optctrl_defx_1},~\eqref{eqt: optctrl_defx_2}, \eqref{eqt: optctrl_defx_3},~\eqref{eqt: optctrl_defx_4},~\eqref{eqt: optctrl_defx_5}, and~\eqref{eqt: optctrl_defx_neg} for different cases of $x,t$ and $p$.

Next, we provide some theoretical properties for the functions $\valuefn$ and $\opttraj$ defined above.
In Proposition~\ref{prop: SJ_solves_HJ1}, we prove that,  under some assumptions, the function $\valuefn$ defined above is indeed the unique viscosity solution to the HJ PDE~\eqref{eqt: result_HJ1_1d}. In Proposition~\ref{prop:optctrl_1d}, we show that the function $\opttraj$ is the unique optimal trajectory in~\eqref{eqt: result_optctrl1_1d}, whose optimal value equals $\valuefn(x,t)$.

\begin{prop} \label{prop: SJ_solves_HJ1}
Let $\initcond\colon \R\to\R$ be a convex function. Let $a,b>0$ be positive constants and $\potentfn\colon \R\to\R$ be defined by~\eqref{eqt: result_defV} with parameters $a$ and $b$. 
Let $\valuefn\colon\Rn\times [0,+\infty)\to\R$ be the function defined in~\eqref{eqt: result_Hopf1_1d}. Then, the following statements hold:
\begin{itemize}
    \item[(a)] The function $\valuefn$ is a continuously differentiable solution to the HJ PDE~\eqref{eqt: result_HJ1_1d}. 
    \item[(b)] If $\initcond$ satisfies
\begin{equation} \label{eqt: HJ1_1d_condJ}
    |\initcond(x) - \initcond(y)|\leq C|x-y|(1 + |x|^\delta + |y|^\delta) \quad\quad \forall x,y\in\R,
\end{equation}
for some constants $C\geq 0$ and $\delta \geq 0$, then the function $\valuefn$ is the unique viscosity solution to the HJ PDE~\eqref{eqt: result_HJ1_1d} in the solution set $\mathcal{G}$ defined by
\begin{equation*}
\begin{split}
    \mathcal{G} := \{&W \in C(\R\times [0, +\infty))\colon
    \|W\|_R < \infty \ \forall R > 0, \exists \alpha\in\R  s.t.\ \forall T>0, \\
    & \, \exists C_T\in\R \ s.t.\ W(x,t)- \alpha x \geq C_T \, \forall x\in\R,\forall t\in[0,T]\},
\end{split}
\end{equation*}
where $\|W\|_R$ is defined by $\|W\|_R := \sup\{|W(x, t)| + |q|\colon (x, t)\in B_R(\R)\times [0,R], q \in D_x^-W(x, t)\}$, the set $B_R(\R)$ denotes the closed ball in $\R$ with center $0$ and radius $R$,
and $D_x^-W(x, t)$ denotes the subdifferential of $W$ with respect to $x$ at $(x,t)$.
\end{itemize}
\end{prop}
\begin{proof}
This is a corollary of Proposition~\ref{prop: HJ1_hd_1}.
\end{proof}

\begin{prop} \label{prop:optctrl_1d}
Let $\initcond\colon \R\to\R$ be a convex function satisfying~\eqref{eqt: HJ1_1d_condJ} and $\potentfn$ be the function defined in~\eqref{eqt: result_defV} with parameters $a,b>0$. Let $x\in\R$ and $t>0$. Then, the unique optimal trajectory for the optimal control problem~\eqref{eqt: result_optctrl1_1d} is given by the function $[0,t]\ni s\mapsto \opttraj(s;x,t)\in\R$ defined in~\eqref{eqt: optctrl_traj_J}.
Moreover, the optimal value of the optimal control problem~\eqref{eqt: result_optctrl1_1d} equals $\valuefn(x,t)$, as defined in~\eqref{eqt: result_Hopf1_1d}.
\end{prop}
\begin{proof}
This is a corollary of Proposition~\ref{prop: optctrl1_hd}.
\end{proof}

\begin{rem}
If $\initcond$ is a linear function, i.e., there exists a scalar $p\in\R$ such that $\initcond(x)=px$ holds for all $x\in\R$, then $\initcond$ satisfies the assumption~\eqref{eqt: HJ1_1d_condJ}. In this case, Proposition~\ref{prop: SJ_solves_HJ1} shows that the function $(x,t)\mapsto \valuefn(x,t;p,a,b)$ defined in~\eqref{eqt: result_S1_1d} and~\eqref{eqt: result_S1_1d_negative} (where $p$ is the slope of $\initcond$) is the unique continuously differentiable solution in the solution set $\mathcal{G}$ to the HJ PDE~\eqref{eqt: result_HJ1_1d} with this linear initial data $\initcond$. Moreover, Proposition~\ref{prop:optctrl_1d} shows that the trajectory $s\mapsto\opttraj(s;x,t,p,a,b)$ defined by~\eqref{eqt: optctrl_defx_1},~\eqref{eqt: optctrl_defx_2},~\eqref{eqt: optctrl_defx_3},~\eqref{eqt: optctrl_defx_4},~\eqref{eqt: optctrl_defx_5}, and~\eqref{eqt: optctrl_defx_neg} for different cases is the unique optimal trajectory of the optimal control problem~\eqref{eqt: result_optctrl1_1d}, whose optimal value equals $\valuefn(x,t;p,a,b)$.
\end{rem}

\subsection{Separable high-dimensional case}\label{sec:HJhd_separable}
In this section, we consider a special case of the high-dimensional problems~\eqref{eqt: result_optctrl1_hd_general} and~\eqref{eqt: HJhd_1_general}.
To be specific, we consider the following high-dimensional optimal control problem:
\begin{equation}\label{eqt: result_optctrl1_hd}
    \valuefn(\bx,t) = \inf \left\{\int_0^t \left(\frac{1}{2}\|\dot{\bx}(s)\|^2 - \potentfn(\bx(s))\right) ds + \initcond(\bx(0)) \colon \bx(t) = \bx\right\}.
\end{equation}
In the optimal control problem~\eqref{eqt: result_optctrl1_hd}, the initial cost $\initcond\colon\Rn\to\R$ is a convex function, and the function $\potentfn\colon\R^n\to(-\infty, 0]$ satisfies
\begin{equation*}
    \potentfn(\bx) = \sum_{i=1}^n \potentfn_i(x_i) \quad\forall\, \bx = (x_1,\dots, x_n)\in\R^n,
\end{equation*}
where each function $\potentfn_i\colon\R\to (-\infty,0]$ is a 1-homogeneous concave function given by~\eqref{eqt: result_defV} with positive constants $a_i$ and $b_i$. 
The corresponding HJ PDE reads:
\begin{equation} \label{eqt: result_HJ1_hd}
\begin{dcases} 
\frac{\partial \valuefn}{\partial t}(\bx,t)+\frac{1}{2}\|\nabla_{\bx}\valuefn(\bx,t)\|^2 + \potentfn(\bx) = 0 & \bx\in\mathbb{R}^n, t\in(0,+\infty),\\
\valuefn(\bx,0)=\initcond(\bx) & \bx\in\mathbb{R}^n,
\end{dcases}
\end{equation}
where the initial condition $\initcond$ and the potential energy $\potentfn$ are the corresponding functions in~\eqref{eqt: result_optctrl1_hd}.

We will see later that each component of the optimal trajectory is independent from each other as long as the initial momentum $\bp^*$ is determined. In other words, the computation of the optimal trajectory can be done in parallel, and hence, we call this problem separable.

The solution $\valuefn\colon \R^n\times\R\to\R$ is defined by the following Hopf-type formula:
\begin{equation}\label{eqt: result_Hopf1_hd}
    \valuefn(\bx,t) := \sup_{\bp\in\R^n} \left\{\sum_{i=1}^n \valuefn(x_i, t; p_i, a_i,b_i) - \initcond^*(\bp)\right\} \quad \forall\, \bx\in\Rn, \, t\geq 0,
\end{equation}
where the function $(x_i, t)\mapsto \valuefn(x_i, t; p_i,a_i,b_i)$ on the right-hand side is defined in~\eqref{eqt: result_S1_1d} and~\eqref{eqt: result_S1_1d_negative}. 
By Lemma~\ref{lem: prop_barp_hopf_hd}, the function value $\valuefn(\bx,t)$ defined in~\eqref{eqt: result_Hopf1_hd} is finite and the maximizer in~\eqref{eqt: result_Hopf1_hd} exists. Moreover, for a positive time horizon $t>0$, the maximizer is unique and we denote the unique maximizer by $\bp^* = (p_1^*, \dots, p_n^*)\in\R^n$.
Define the trajectory $[0,t]\ni s\mapsto \opttrajhd(s; \bx,t)\in\R^n$ by 
\begin{equation}\label{eqt: optctrlhd_traj_J}
    \opttrajhd(s; \bx,t) := \left(\opttraj(s; x_1,t,p_1^*, a_1,b_1), \dots, \opttraj(s; x_n,t,p_n^*, a_n,b_n)\right) \quad \forall\, s\in[0,t],
\end{equation}
where the $i$-th element $\opttraj(s; x_i,t,p_i^*, a_i,b_i)$ on the right-hand side is the one-dimensional trajectory defined in~\eqref{eqt: optctrl_defx_1},~\eqref{eqt: optctrl_defx_2},~\eqref{eqt: optctrl_defx_3},~\eqref{eqt: optctrl_defx_4},~\eqref{eqt: optctrl_defx_5}, and~\eqref{eqt: optctrl_defx_neg} for different cases of $x_i,t$ and $p_i^*$. Note that the components of $\opttrajhd(s; \bx,t)$ are independent from each other, and hence they can be computed in parallel as long as $\bp^*$ is known. Thus, we call this problem separable.

Now, we provide some theoretical guarantees for the functions $\valuefn$ and $\opttrajhd$ defined above. Proposition~\ref{prop: HJ1_hd_1} shows that the function $\valuefn$ is indeed the unique viscosity solution to the HJ PDE~\eqref{eqt: result_HJ1_hd} under some assumptions. Moreover, Proposition~\ref{prop: optctrl1_hd} proves that the function $\opttrajhd$ is the unique optimal trajectory in~\eqref{eqt: result_optctrl1_hd} under some assumptions and that the corresponding optimal value equals $\valuefn(\bx,t)$.

\begin{prop} \label{prop: HJ1_hd_1}
Let $\initcond\colon \Rn\to\R$ be a convex function. Let $\{a_i,b_i\}_{i=1}^n$ be positive constants and $\potentfn_i\colon \R\to\R$ be defined by~\eqref{eqt: result_defV} with constants $a_i$ and $b_i$ for each $i\in\{1,\dots, n\}$. Let $\valuefn\colon \R^n\times [0,+\infty)\to\R$ be the function defined in~\eqref{eqt: result_Hopf1_hd}. Then, the following results hold:
\begin{itemize}
    \item[(a)] The function $\valuefn$ is a continuously differentiable solution to the HJ PDE~\eqref{eqt: result_HJ1_hd}.
    \item[(b)] Furthermore, assume $\initcond$ satisfies \begin{equation} \label{eqt: HJ1_hd_condJ}
    |\initcond(\bx) - \initcond(\by)|\leq C\|\bx-\by\|(1 + \|\bx\|^\delta + \|\by\|^\delta) \quad\quad \forall \bx,\by\in\R^n,
\end{equation}
for some constants $C\geq 0$ and $\delta \geq 0$. Then, the function $\valuefn$ is the unique viscosity solution to the HJ PDE~\eqref{eqt: result_HJ1_hd} in the solution set $\mathcal{G}$ defined by
\begin{equation}\label{eqt: HJ1_hd_solset_G}
\begin{split}
    \mathcal{G} := \{&W \in C(\R^n\times [0, +\infty))\colon
    \|W\|_R < \infty \ \forall R > 0, \exists \balp\in\R^n\,\,s.t.\,\, \forall T>0,\\
    &\exists C_T\in\R \,\,s.t.\,\, W(\bx,t)- \langle \balp, \bx\rangle\geq C_T \,\, \forall \bx\in\Rn,\forall t\in[0,T]\},
\end{split}
\end{equation}
where $\|W\|_R$ is defined by 
\begin{equation}\label{eqt:def_norm_R}
\|W\|_R := \sup\{|W(\bx, t)| + \|\bq\|\colon (\bx, t)\in B_R(\R^n)\times [0,R], \bq \in D_{\bx}^-W(\bx, t)\},
\end{equation}
where the set $B_R(\Rn)$ denotes the closed ball in $\Rn$ with center $\mathbf{0}$ and radius $R$ and $D_{\bx}^-W(\bx, t)$ denotes the subdifferential of $W$ with respect to $\bx$ at $(\bx,t)$.
\end{itemize}
\end{prop}
\begin{proof}
(a) 
By Lemma~\ref{lem: continuity_SJ_pde1_hd}, the function $\valuefn$ is continuous in $\Rn\times [0,+\infty)$.
By Lemma~\ref{lem: lem24_grad_SJ}, the function $\valuefn$ is continuously differentiable, and its gradient at any point $(\bx,t)\in\Rn\times (0,+\infty)$ satisfies
\begin{equation*}
\begin{split}
    \frac{\partial \valuefn(\bx,t)}{\partial t} 
    &= \sum_{i=1}^n\frac{\partial \valuefn}{\partial t}(x_i,t;p_i^*(\bx,t),a_i,b_i)\\
    &= - \sum_{i=1}^n \left(\frac{1}{2}\left(\frac{\partial \valuefn}{\partial x}(x_i,t;p_i^*(\bx,t),a_i,b_i)\right)^2 + \potentfn_i(x_i)\right)\\
    &= - \frac{1}{2}\|\nabla_{\bx} \valuefn(\bx,t)\|^2 - \potentfn(\bx),
\end{split}
\end{equation*}
where $p^*_i(\bx,t)$ denotes the $i$-th component of the unique maximizer in~\eqref{eqt: result_Hopf1_hd} at $(\bx,t)$, the first and the third equalities hold by~\eqref{eqt: lem24_grad_SJ}, and the second equality holds since each function $(x_i,t)\mapsto \valuefn(x_i,t;p^*_i,a_i,b_i)$ satisfies their corresponding one-dimensional HJ PDE by Lemma~\ref{prop: S_linearJ}. Hence, the function $\valuefn$ satisfies the differential equation in~\eqref{eqt: result_HJ1_hd}.
Also, the initial condition is satisfied according to~\eqref{eqt: lem23_SJ_initialcond} in the proof of Lemma~\ref{lem: prop_barp_hopf_hd}. Therefore, the function $\valuefn$ is a continuously differentiable solution to the HJ PDE~\eqref{eqt: result_HJ1_hd}.

(b) To prove that the function $\valuefn$ is the unique viscosity solution in the solution set $\mathcal{G}$, we first prove that the function $\valuefn$ is in $\mathcal{G}$. 
Let $\bp$ be any vector in the domain of $\initcond^*$. By~\eqref{eqt: result_Hopf1_hd} and Lemma~\ref{lem: appendix_S_bd}, we have
\begin{equation} \label{eqt:prop24_pf_S_lowerbd}
\begin{split}
    \valuefn(\bx,t)&\geq \sum_{i=1}^n \valuefn(x_i,t;p_i,a_i,b_i) - \initcond^*(\bp) \geq \sum_{i=1}^n (p_ix_i + C_{i}) - \initcond^*(\bp) \\
    &= \langle \bp,\bx\rangle + \sum_{i=1}^nC_{i} - \initcond^*(\bp),
\end{split}
\end{equation}
for all $\bx\in\Rn$ and $t\geq 0$,
where $C_{i}$ is the constant $C$ in the lower bound in Lemma~\ref{lem: appendix_S_bd} with constants $a=a_i$, $b=b_i$, and $p=p_i$.
Hence, $\valuefn$ is bounded below by an affine function.
Then, to prove $\valuefn\in\mathcal{G}$, it remains to prove that $\|\valuefn\|_R$ is finite for all $R>0$. Let $R>0$ be an arbitrary number, and let $\bp^*(\bx,t)$ denote the set of maximizers in~\eqref{eqt: result_Hopf1_hd} for any $\bx\in\Rn$ and $t\geq 0$. If the set is a singleton, by a slight abuse of notation, we denote both the set and the element in the set by $\bp^*(\bx,t)$ and we denote the $i$-th component of the element by $p_i^*(\bx,t)$.
By Lemma~\ref{lem: lem24_grad_SJ} and straightforward computation, we get
\begin{equation}\label{eqt: prop23_subdiff_SJ}
   D_{\bx}^-\valuefn(\bx, t) = \begin{dcases}
   \{\nabla_{\bx} \valuefn(\bx,t) \} & \bx\in\Rn, t>0,\\
   \partial \initcond(\bx) = \bp^*(\bx,0) & \bx\in\Rn, t=0,
   \end{dcases}
\end{equation}
where
\begin{equation*}
     \{\nabla_{\bx} \valuefn(\bx,t) \} = \left\{\left(\frac{\partial \valuefn(x_1,t;p^*_1(\bx,t),a_1,b_1)}{\partial x},\dots,   \frac{\partial \valuefn(x_n,t;p^*_n(\bx,t),a_n,b_n)}{\partial x}
   \right)\right\}.
\end{equation*}
For any $x,p\in\R$ and $a,b>0$, the function $x\mapsto \valuefn(x,0;p,a,b)$ equals $xp$, and hence we get $\frac{\partial \valuefn(x,0;p,a,b)}{\partial x}=p$. Therefore,~\eqref{eqt: prop23_subdiff_SJ} is simplified to
\begin{equation}\label{eqt: prop23_subdiff_SJ2}
\begin{split}
    & D_{\bx}^-\valuefn(\bx, t) \\
    & = \Bigg\{\left(\frac{\partial \valuefn(x_1,t;p_1,a_1,b_1)}{\partial x},\dots,
   \frac{\partial \valuefn(x_n,t;p_n,a_n,b_n)}{\partial x}
   \right)\colon (p_1,\dots,p_n)\in \bp^*(\bx,t)\Bigg\},
\end{split}
\end{equation}
for any $\bx\in\Rn, t\geq 0$.
By Lemma~\ref{lem: prop_barp_hopf_hd}(e), the set $\{\|\bp\|\colon \bx\in B_R(\Rn),\, t\in [0,R],\, \bp\in \bp^*(\bx,t)\}$ is bounded for all $R>0$, and we denote the bound by $R_p$. Then, by Lemma~\ref{lem: appendix_dSdx_bd}, for all $\bx\in B_R(\Rn)$, $t\in[0,R]$, and $\bp\in \bp^*(\bx,t)$, we get
\begin{equation}\label{eqt:prop24_pf_gradS_bd}
\begin{split}
   &\left\|\left(\frac{\partial \valuefn(x_1,t;p_1,a_1,b_1)}{\partial x},\dots, 
   \frac{\partial \valuefn(x_n,t;p_n,a_n,b_n)}{\partial x}
   \right)\right\|\\
   =&\sqrt{\sum_{i=1}^n \left|\frac{\partial \valuefn(x_i,t;p_i,a_i,b_i)}{\partial x}\right|^2}
   \leq \sqrt{\sum_{i=1}^nC_i(R,R,R_p)^2},
\end{split}
\end{equation}
where $C_i$ is the function in the upper bound defined in Lemma~\ref{lem: appendix_dSdx_bd} for the parameters $a_i$ and $b_i$. (Note that in different contexts, we may reuse the notation $C_i$ to denote different bounds if there is no ambiguity.) Combining~\eqref{eqt: prop23_subdiff_SJ2} and~\eqref{eqt:prop24_pf_gradS_bd}, we have
\begin{equation*}
    \|\valuefn\|_R \leq \sup_{\bx\in B_R(\R), t\in [0,R]}|\valuefn(\bx,t)| + \sqrt{\sum_{i=1}^nC_i(R,R,R_p)^2} < +\infty.
\end{equation*}
Therefore, $\valuefn$ is a function in $\mathcal{G}$. We have proved in (a) that $\valuefn$ is a continuously differentiable solution, and hence, $\valuefn$ is a viscosity solution in $\mathcal{G}$.

Then, we apply Lemma~\ref{lem:A10_uniqueness} to prove the uniqueness of the viscosity solution. To apply Lemma~\ref{lem:A10_uniqueness}, we need to check its assumptions. The assumption on $\initcond$ is satisfied since $\initcond$ is a convex function satisfying~\eqref{eqt: HJ1_hd_condJ}. The assumption on $M$ is satisfied since $M$ is the identity matrix in this section. The assumption on $\potentfn$ is satisfied since $\potentfn$ is a non-positive 1-homogeneous concave function, which implies that $\potentfn$ is Lipschitz continuous. Now, it remains to check the assumption on $\valuefn$ in Lemma~\ref{lem:A10_uniqueness}(b). Let $\balp\in\Rn$ be a vector such that $\bx\mapsto \initcond(\bx)-\langle \balp, \bx\rangle $ is bounded from below. The convexity of $\initcond$ implies that $\balp$ is in the domain of $\initcond^*$. Using~\eqref{eqt:prop24_pf_S_lowerbd} with $\bp=\balp$, we have that
\begin{equation*}
    \valuefn(\bx,t) - \langle \balp,\bx\rangle \geq \sum_{i=1}^nC_{i} - \initcond^*(\balp)\in\R.
\end{equation*}
In other words, the function $(\bx,t)\mapsto \valuefn(\bx,t) - \langle \balp, \bx\rangle $ is bounded from below. Therefore, the assumptions for Lemma~\ref{lem:A10_uniqueness}(b) holds, and the viscosity solution to~\eqref{eqt: result_HJ1_hd} in $\mathcal{G}$ is unique.
\end{proof}

\begin{prop} \label{prop: optctrl1_hd}
Let $\initcond\colon \Rn\to\R$ be a convex function satisfying~\eqref{eqt: HJ1_hd_condJ}. Let $a_1,\dots,a_n,b_1,\dots, b_n$ be positive constants and $\potentfn_i\colon \R\to\R$ be defined by~\eqref{eqt: result_defV} with constants $a = a_i$ and $b = b_i$ for each $i\in\{1,\dots, n\}$. 
Then, for any $\bx\in\R^n$ and $t>0$, the unique optimal trajectory for the optimal control problem~\eqref{eqt: result_optctrl1_hd} is given by the function $[0,t]\ni s\mapsto \opttrajhd(s;\bx,t)\in\R^n$ defined in~\eqref{eqt: optctrlhd_traj_J}.
Moreover, the optimal value of the optimal control problem~\eqref{eqt: result_optctrl1_hd} equals $\valuefn(\bx,t)$ defined in~\eqref{eqt: result_Hopf1_hd}.
\end{prop}

\begin{proof}
Let $\bx\in\R^n$ and $t>0$.
In this proof, we write $\opttrajhd(s)$ instead of $\opttrajhd(s;\bx,t)$ whenever there is no ambiguity.
By Lemma~\ref{lem: appendix_optval_equalS} and the definition of $\opttrajhd$, we have
\begin{equation}\label{eqt: prop27_costeq1}
\begin{split}
&\int_0^t \left(\frac{\|\dot{\opttrajhd}(s)\|^2}{2} - \potentfn(\opttrajhd(s))\right) ds + \initcond(\opttrajhd(0)) \\
=\ & \sum_{i=1}^n\int_0^t \left(\frac{(\dot{\opttrajhd}_i(s))^2}{2} - \potentfn_i(\opttrajhd_i(s))\right) ds + \initcond(\opttrajhd(0))\\
=\ & \sum_{i=1}^n \valuefn(x_i,t;p^*_i,a_i,b_i) - \langle \bp^*,\opttrajhd(0)\rangle + \initcond(\opttrajhd(0)),
\end{split}
\end{equation}
where $\opttrajhd_i$ denotes the $i$-th component of $\opttrajhd$ and $\bp^*=(p_1^*,\dots,p_n^*)$ is the unique maximizer in~\eqref{eqt: result_Hopf1_hd}.
By Lemma~\ref{lem: appendix_gradpS_gamma0}, for each $i\in\{1,\dots,n\}$, we have $\frac{\partial \valuefn}{\partial p}(x_i,t;p^*_i,a_i,b_i) = \opttraj(0;x_i,t,p^*_i,a_i,b_i)$. 
Note that~\eqref{eqt: result_Hopf1_hd} is a concave optimization problem by Lemma~\ref{lem: prop_S} and the convexity of $\initcond^*$.
Since $\bp^*=(p^*_1,\dots,p^*_n)$ is the maximizer in~\eqref{eqt: result_Hopf1_hd}, by the first order optimality condition, we have 
\begin{equation*}
    \left(\frac{\partial \valuefn}{\partial p}(x_1,t;p^*_1,a_1,b_1), \dots, \frac{\partial \valuefn}{\partial p}(x_n,t;p^*_n,a_n,b_n)\right)\in \partial \initcond^*(\bp^*),
\end{equation*}
where $\partial \initcond^*$ denotes the subdifferential operator of $\initcond^*$.
Therefore, we conclude that $\opttrajhd(0)\in \partial \initcond^*(\bp^*)$, which implies that $\initcond(\opttrajhd(0)) - \langle \bp^*,\opttrajhd(0)\rangle + \initcond^*(\bp^*) = 0$. In other words, by~\eqref{eqt: prop27_costeq1}, we have
\begin{equation*}
\begin{split}
\int_0^t \left(\frac{\|\dot{\opttrajhd}(s)\|^2}{2} - \potentfn(\opttrajhd(s))\right) ds + \initcond(\opttrajhd(0)) &= \sum_{i=1}^n \valuefn(x_i,t;p^*_i,a_i,b_i) - \initcond^*(\bp^*) \\
&= \valuefn(\bx,t).
\end{split}
\end{equation*}
Moreover, by Lemma~\ref{lem:A10_uniqueness}(a) with $M$ being the identity matrix (whose assumptions are proved in the proof of Proposition~\ref{prop: HJ1_hd_1}(b)), the value $\valuefn(\bx,t)$ is the optimal value of the optimal control problem~\eqref{eqt: result_optctrl1_hd}.
Hence, the cost of $\opttrajhd$ equals the optimal cost $\valuefn(\bx,t)$. By straightforward calculation, the function $\opttrajhd$ is Lipschitz continuous and satisfies $\opttrajhd(t) = \bx$.
Therefore, $\opttrajhd$ is an optimal trajectory, and the corresponding optimal value equals $\valuefn(\bx,t)$. The optimal trajectory is unique since the problem~\eqref{eqt: result_optctrl1_hd} is a strictly convex problem.
\end{proof}

\begin{rem}
For the special case when $\initcond$ is a linear function, i.e., $\initcond(\bx) = \langle \bp, \bx\rangle$ for all $\bx\in\R^n$ and for some constant vector $\bp =(p_1,\dots, p_n)\in \Rn$, the function $\initcond^*$ equals the indicator function of $\{\bp\}$, and hence the solution to the corresponding HJ PDE~\eqref{eqt: result_HJ1_hd} reads:
\begin{equation*}
    \valuefn(\bx,t) = \sum_{i=1}^n \valuefn(x_i, t; p_i, a_i,b_i) \quad \forall \bx\in\Rn, t\geq 0,
\end{equation*}
where the function $\valuefn$ in the summation on the right-hand side is defined by~\eqref{eqt: result_S1_1d} and~\eqref{eqt: result_S1_1d_negative}. 
In this case, the optimal trajectory of the corresponding optimal control problem~\eqref{eqt: result_optctrl1_hd} equals 
\begin{equation*}
    \opttrajhd(s; \bx,t) = \left(\opttraj(s; x_1,t,p_1, a_1,b_1), \dots, \opttraj(s; x_n,t,p_n, a_n,b_n)\right) \quad \forall\, s\in[0,t],
\end{equation*}
where each component $\opttraj(s; x_i,t,p_i, a_i,b_i)$ on the right-hand side is defined by~\eqref{eqt: optctrl_defx_1},~\eqref{eqt: optctrl_defx_2},~\eqref{eqt: optctrl_defx_3},~\eqref{eqt: optctrl_defx_4},~\eqref{eqt: optctrl_defx_5}, and~\eqref{eqt: optctrl_defx_neg} for different cases.
\end{rem}
\def\uz{\bu_0}

\subsection{The general high-dimensional case}\label{sec:HJhd_general}

In this section, we provide the analytical solution to the high-dimensional problems~\eqref{eqt: result_optctrl1_hd_general} and~\eqref{eqt: HJhd_1_general} under more general assumptions.
Assume there exists a vector $\uz\in\Rn$ such that 
the function $\bx\mapsto \potentfn(\bx + \uz)$ is 1-homogeneous and there exists an invertible matrix $\matP$ with $n$ rows and $n$ columns whose column vectors $\bu_1,\dots, \bu_n\in\Rn$ satisfy
\begin{equation}\label{eqt: prop24_levelset_V}
    \{\bx\in\Rn\colon \potentfn(\bx + \uz)\geq -1\} = \co \left(\bigcup_{j=1}^n \left\{\frac{1}{a_j}\bu_j, -\frac{1}{b_j}\bu_j\right\}\right),
\end{equation}
for some positive scalars $a_i,b_i>0, i\in\{1,\dots, n\}$, where 
$\co E$ denotes the convex hull of a set $E$.

We define the function $\valuefn\colon \Rn\times[0,+\infty)\to\R$ by the following Hopf-type formula:
\begin{equation} \label{eqt: result_HJ1_Hopf_hd_general}
    \valuefn(\bx,t) := \sup_{\bp\in\R^n} \left\{\sum_{i=1}^n \valuefn\left((\matP^{-1}\bx-\matP^{-1}\uz)_i, t; p_i, a_i,b_i\right) - \tilde{\initcond}^*(\bp)\right\},
\end{equation}
for any $\bx\in\Rn, t\geq 0$,
where the function $\valuefn$ in the summation on the right-hand side is defined by~\eqref{eqt: result_S1_1d} and~\eqref{eqt: result_S1_1d_negative} and the function $\tilde{\initcond}^*$ is the Legendre-Fenchel transform of the function $\tilde{\initcond}\colon \Rn\to\R$, which is defined by
\begin{equation}\label{eqt: J_changeofvariable}
    \tilde{\initcond}(\by):= \initcond(\matP\by+\uz)\quad \forall\,\by\in\Rn.
\end{equation}
By straightforward calculation, the function $\tilde{\initcond}^*$ equals
\begin{equation}\label{eqt: J_star_changeofvariable}
    \tilde{\initcond}^*(\bp) = \initcond^*\left((\matP^{-1})^T\bp\right) - \langle \bp, \matP^{-1}\uz\rangle \quad \forall \bp\in\Rn.
\end{equation}
Hence, for any $\bx\in\Rn$ and $t\geq 0$, the Hopf-type formula~\eqref{eqt: result_HJ1_Hopf_hd_general} can equivalently be formulated as:
\begin{equation}\label{eqt:Hopf_hd_general_ver2}
\begin{split}
    &\valuefn(\bx,t) \\
    =& \sup_{\bp\in\R^n} \left\{\sum_{i=1}^n \valuefn\left((\matP^{-1}\bx-\matP^{-1}\uz)_i, t; p_i, a_i,b_i\right) - \initcond^*\left((\matP^{-1})^T\bp\right) + \langle \bp, \matP^{-1}\uz\rangle\right\}\\
    =&\sup_{\bq\in\R^n} \left\{\sum_{i=1}^n \valuefn\left((\matP^{-1}\bx-\matP^{-1}\uz)_i, t; (\matP^T\bq)_i, a_i,b_i\right) - \initcond^*(\bq) + \langle \bq, \uz\rangle\right\},
\end{split}
\end{equation}
where the second equality holds by the change of variable $\bq = (\matP^{-1})^T\bp$.
Define the trajectory $[0,t]\ni s\mapsto \opttrajhd(s; \bx,t)\in\R^n$ by 
\begin{equation}\label{eqt: optctrlhd_traj_general}
    \opttrajhd(s; \bx,t) :=
    \matP\tilde{\opttrajhd}(s; \matP^{-1}(\bx-\bu_0),t) + \bu_0,
\end{equation}
where the function $s\mapsto \tilde{\opttrajhd}(s; \by,t)$ for any $\by=(y_1,\dots, y_n)\in\Rn$ and $t>0$ is defined by
\begin{equation}\label{eqt: optctrl1_hd_general_opty}
    \tilde{\opttrajhd}(s; \by, t) := \left(\opttraj(s; y_1, t, p_1^*, a_1, b_1),\dots, \opttraj(s; y_n, t, p_n^*, a_n, b_n)\right) \quad\forall s\in[0,t],
\end{equation}
where $\bp^*=(p_1^*,\dots,p_n^*)$ is the maximizer in~\eqref{eqt: result_HJ1_Hopf_hd_general} and the $i$-th component $\opttraj(s; y_i,t,p_i^*, a_i,b_i)$ on the right-hand side is the one-dimensional trajectory defined by~\eqref{eqt: optctrl_defx_1},~\eqref{eqt: optctrl_defx_2},~\eqref{eqt: optctrl_defx_3},~\eqref{eqt: optctrl_defx_4},~\eqref{eqt: optctrl_defx_5}, and~\eqref{eqt: optctrl_defx_neg} for different cases of $y_i,t$, and $p_i^*$.
Note that for $t>0$, by Lemma~\ref{lem: prop_S}, the negative of the objective function in~\eqref{eqt: result_HJ1_Hopf_hd_general} is 1-coercive and strictly convex. Therefore, the optimal value $\valuefn(\bx,t)$ in~\eqref{eqt: result_HJ1_Hopf_hd_general} is finite, and the maximizer $\bp^*$ exists and is unique. Hence, the functions $\valuefn$ and $\opttrajhd$ are well-defined.

The following propositions show that the functions $\opttrajhd$ and $\valuefn$ defined above do indeed solve the problems~\eqref{eqt: result_optctrl1_hd_general} and~\eqref{eqt: HJhd_1_general}. Proposition~\ref{prop: HJ1_hd_2} proves that the function $\valuefn$ is the unique viscosity solution to the HJ PDE~\eqref{eqt: HJhd_1_general} under some assumptions. Proposition~\ref{prop: HJ1_opttraj_general} shows that the function $\opttrajhd$ is the unique optimal trajectory of the optimal control problem~\eqref{eqt: result_optctrl1_hd_general} under some assumptions and that the corresponding optimal value equals $\valuefn(\bx,t)$.

\begin{prop} \label{prop: HJ1_hd_2}
Let $\initcond\colon\R^n\to\R$ be a convex function.
Let $\potentfn\colon \R^n\to (-\infty, 0]$ be a piecewise affine concave function. Assume there exists a vector $\uz\in\Rn$ such that $\bx\mapsto \potentfn(\bx + \uz)$ is 1-homogeneous and there exists an invertible matrix $\matP$ whose column vectors $\bu_1,\dots,\bu_n\in\Rn$ satisfy~\eqref{eqt: prop24_levelset_V} with the vector $\uz$ and some positive scalars $a_1,\dots, a_n,b_1,\dots,b_n>0$. 
Let $M$ be a matrix with $n$ rows and $n$ columns satisfying $M = \matP \matP^T$. Let $\valuefn\colon\R^n\times [0,+\infty)\to\R$ be the function defined in~\eqref{eqt: result_HJ1_Hopf_hd_general}. Then, the function $\valuefn$ is a continuously differentiable solution to the HJ PDE~\eqref{eqt: HJhd_1_general}. Moreover, if $\initcond$ satisfies~\eqref{eqt: HJ1_hd_condJ}, the function
$\valuefn$ is the unique viscosity solution to the HJ PDE~\eqref{eqt: HJhd_1_general} in the solution set $\mathcal{G}$ defined in~\eqref{eqt: HJ1_hd_solset_G}.
\end{prop}
\begin{proof}
Define $\tilde{\valuefn}\colon \Rn\times [0,+\infty)\to \R\cup \{+\infty\}$ by
\begin{equation}\label{eqt: prop25_def_tildeS}
    \tilde{\valuefn}(\by,t) := \valuefn(\matP\by + \bu_0,t) = \sup_{\bp\in\R^n} \left\{\sum_{i=1}^n \valuefn(y_i, t; p_i, a_i,b_i) - \tilde{\initcond}^*(\bp)\right\},
\end{equation}
for any $\by\in\R^n,t\geq 0$,
where the second equality follows directly from the definition of $\valuefn$ in~\eqref{eqt: result_HJ1_Hopf_hd_general}. By definition~\eqref{eqt: J_changeofvariable}, the function $\tilde{\initcond}$ is the composition of the convex function $\initcond$ with an affine map, and hence, $\tilde{\initcond}$ is also convex. By~\eqref{eqt: prop25_def_tildeS} and Proposition~\ref{prop: HJ1_hd_1}, $\tilde{\valuefn}$ is a continuously differentiable solution to the HJ PDE~\eqref{eqt: result_HJ1_hd} with the convex initial data $\tilde{\initcond}$. By straightforward calculation, we have
\begin{equation}\label{eqt: prop_hd2_gradtildeS}
    \nabla \tilde{\valuefn}(\by,t) = \left(\matP^T\nabla_{\bx} \valuefn(\matP\by + \bu_0,t), \frac{\partial \valuefn}{\partial t}(\matP\by + \bu_0,t)\right).
\end{equation}
Applying~\eqref{eqt: prop_hd2_gradtildeS} and the change of variable $\bx = \matP\by+\bu_0$ to the HJ PDE~\eqref{eqt: result_HJ1_hd}, we conclude that $\valuefn$ is a continuously differentiable solution of the following PDE:
\begin{equation}\label{eqt:prop25_pf_HJ}
    \begin{dcases} 
\frac{\partial \valuefn}{\partial t}(\bx,t)+\frac{\|\matP^T\nabla_{\bx}\valuefn(\bx,t)\|^2}{2} + \sum_{i=1}^n \potentfn_i((\matP^{-1}\bx-\matP^{-1}\uz)_i) = 0 & \bx\in\mathbb{R}^n, t>0,\\
\valuefn(\bx,0)=\initcond(\bx) & \bx\in\mathbb{R}^n.
\end{dcases}
\end{equation}
Since we assume $M=\matP \matP^T$, after some computation, we get that
\begin{equation*}
\begin{split}
    \frac{1}{2}\|\matP^T\nabla_{\bx}\valuefn(\bx,t)\|^2
    &= \frac{1}{2}\nabla_{\bx}\valuefn(\bx,t)^T\matP \matP^T \nabla_{\bx}\valuefn(\bx,t)
    = \frac{1}{2}\nabla_{\bx}\valuefn(\bx,t)^TM \nabla_{\bx}\valuefn(\bx,t) \\
    &= \frac{1}{2}\|\nabla_{\bx}\valuefn(\bx,t)\|^2_M.
\end{split}
\end{equation*}
Then, to prove that $\valuefn$ solves~\eqref{eqt: HJhd_1_general}, it suffices to prove that $\tilde{\potentfn}= \potentfn$, where $\tilde{\potentfn}\colon \Rn\to (-\infty, 0]$ is defined by
\begin{equation*}
    \tilde{\potentfn}(\bx):= \sum_{i=1}^n \potentfn_i((\matP^{-1}\bx-\matP^{-1}\uz)_i) \quad \forall\,\bx\in\Rn.
\end{equation*}
Recall that each function $\potentfn_i\colon \R\to(-\infty,0]$ is 1-homogeneous, and hence, the function $\bx\mapsto \tilde{\potentfn}(\bx+\uz) = \sum_{i=1}^n \potentfn_i((\matP^{-1}\bx)_i) \in (-\infty,0]$ is also 1-homogeneous. Since any non-positive 1-homogeneous function is uniquely determined by its superlevel set at $-1$, to prove $\tilde{\potentfn}= \potentfn$, it suffices to prove 
\begin{equation}\label{eqt:prop25_pf_superlevel_equal}
    \{\bx\in\Rn\colon \tilde{\potentfn}(\bx + \uz)\geq -1\}
    =\{\bx\in\Rn\colon \potentfn(\bx + \uz)\geq -1\}.
\end{equation}
Denote by $\unitsim_n$ the simplex set, i.e. define $\unitsim_n$ by 
\begin{equation*}
    \unitsim_n := \{(\alpha_1,\dots,\alpha_n)\in [0,1]^n\colon \sum_{i=1}^n\alpha_i = 1\}.
\end{equation*}
By straightforward calculation, we have
\begin{equation*}
\begin{split}
    &\{\bx\in\Rn\colon \tilde{\potentfn}(\bx + \uz)\geq -1\}\\
    =& \left\{\bx\in\Rn\colon \sum_{i=1}^n \potentfn_i((\matP^{-1}\bx)_i) \geq -1\right\}\\
    =& \bigcup_{\balp\in\unitsim_n}\{\bx\in\Rn\colon \potentfn_i((\matP^{-1}\bx)_i) \geq -\alpha_i \,\forall i\in\{1,\dots,n\}\}\\
    =& \bigcup_{\balp\in\unitsim_n} \left\{\bx\in\Rn\colon (\matP^{-1}\bx)_i\in \left[-\frac{\alpha_i}{b_i}, \frac{\alpha_i}{a_i}\right] \,\forall i\in\{1,\dots,n\}\right\}\\
    =& \left\{\bx\in\Rn\colon \matP^{-1}\bx \in \co\left(\bigcup_{j=1}^n \left\{\frac{\be_j}{a_j}, -\frac{\be_j}{b_j} \right\}\right) \right\}\\
    =&\co \left(\bigcup_{j=1}^n \left\{\frac{1}{a_j}\bu_j, -\frac{1}{b_j}\bu_j\right\}\right)\\
    =&\{\bx\in\Rn\colon \potentfn(\bx + \uz)\geq -1\},
\end{split}
\end{equation*}
where $\be_1,\dots,\be_n\in\Rn$ denote the standard basis vectors in $\Rn$.
Therefore,~\eqref{eqt:prop25_pf_superlevel_equal} holds, and we obtain $\tilde{\potentfn}= \potentfn$. As a result, the HJ PDE~\eqref{eqt:prop25_pf_HJ} coincides with~\eqref{eqt: HJhd_1_general}, and hence, the function $\valuefn$ defined in~\eqref{eqt: result_HJ1_Hopf_hd_general} is a continuously differentiable solution to the HJ PDE~\eqref{eqt: HJhd_1_general}. 

Now assume that $\initcond$ satisfies~\eqref{eqt: HJ1_hd_condJ}. Then, after straightforward calculation, we obtain that for any $\bx,\by\in\Rn$,
\begin{equation*}
\begin{split}
    \left|\tilde{\initcond}(\bx) - \tilde{\initcond}(\by)\right|
    &= |\initcond(\matP\bx + \uz) - \initcond(\matP\by + \uz)|\\
    &\leq C\|\matP\bx-\matP\by\|(1 + \|\matP\bx+\uz\|^\delta + \|\matP\by+\uz\|^\delta) \\
    &\leq C\|\matP\|\|\bx-\by\|\left(1 + 2^{\delta}(\|\matP\|^{\delta}\|\bx\|^{\delta}+\|\uz\|^\delta) + 2^{\delta}(\|\matP\|^{\delta}\|\by\|^{\delta}+\|\uz\|^\delta)\right) \\
    &\leq C\|\matP\|\max\{1+2^{\delta+1}\|\uz\|^{\delta}, 2^{\delta}\|\matP\|^{\delta}\}\|\bx-\by\|\left(1 + \|\bx\|^{\delta}+\|\by\|^{\delta}\right).
\end{split}
\end{equation*}
Hence, $\tilde{\initcond}$ is a convex function satisfying~\eqref{eqt: HJ1_hd_condJ}. By Proposition~\ref{prop: HJ1_hd_1}(b), the function $\tilde{\valuefn}$ defined in~\eqref{eqt: prop25_def_tildeS} is in the solution set $\mathcal{G}$. Since the function $\valuefn$ is the composition of $\tilde{\valuefn}$ and an affine function, it is straightforward to check that the function $\valuefn$ is also in the solution set $\mathcal{G}$.

Now, we prove the uniqueness of the viscosity solution to the HJ PDE~\eqref{eqt: HJhd_1_general} in the solution set $\mathcal{G}$. Note that for any viscosity solution $W\in\mathcal{G}$ to the HJ PDE~\eqref{eqt: HJhd_1_general}, the corresponding function $\tilde{W}\colon \Rn\times [0,+\infty)\to\R$ defined by
\begin{equation*}
\tilde{W}(\by,t) := W(\matP\by + \uz,t)\quad \forall \by\in\R^n, t\geq 0
\end{equation*}
is a viscosity solution to the HJ PDE~\eqref{eqt: result_HJ1_hd} with initial condition $\tilde{\initcond}$.
By Proposition~\ref{prop: HJ1_hd_1}(b), the function $\tilde{W}$ is the unique viscosity solution in the solution set $\mathcal{G}$ to the HJ PDE~\eqref{eqt: result_HJ1_hd} with the initial condition $\tilde{\initcond}$.
Then, the uniqueness of the viscosity solution $W$ follows since we have  $W(\bx,t)=\tilde{W}(\matP^{-1}(\bx-\uz),t)$ by definition of $\tilde{W}$. In other words, the function $\valuefn$ is the unique viscosity solution in the solution set $\mathcal{G}$ to the HJ PDE~\eqref{eqt: HJhd_1_general}.
\end{proof}

\begin{rem}
Under the assumptions of Proposition~\ref{prop: HJ1_hd_2}, the set $\{\bx\in\Rn\colon \potentfn(\bx + \uz)\geq -1\}$ (which is the shifted superlevel set of the function $\potentfn$) is a convex polyhedral with $2n$ vertices.
Moreover, the matrix $M$ in the kinetic energy term $\frac{1}{2}\|\nabla_{\bx} \valuefn(\bx,t)\|_M^2$ is related to the extreme points (and hence the shape) of the set $\{\bx\in\Rn\colon \potentfn(\bx + \uz)\geq -1\}$ in~\eqref{eqt: prop24_levelset_V}. For instance, if $M$ is a diagonal matrix, then the vectors $\bu_1,\dots, \bu_n$ are orthogonal to each other. In this case, the $2n$ vertices of $\{\bx\in\Rn\colon \potentfn(\bx + \uz)\geq -1\}$ can be grouped pairwise into $n$ groups, where the $i$-th group contains the points $\frac{\bu_i}{a_i}$ and $-\frac{\bu_i}{b_i}$. Then, by connecting the two vertices in each group, we obtain $n$ line segments that are pairwise orthogonal to each other.
\end{rem}

\begin{prop} \label{prop: HJ1_opttraj_general}
Let $\initcond\colon \Rn\to\R$ be a convex function satisfying~\eqref{eqt: HJ1_hd_condJ}. Let the function $\potentfn\colon\R^n\to(-\infty,0]$, the matrices $\matP$ and $M$, the vector $\uz\in\Rn$, and the scalars $a_1,\dots, a_n,b_1,\dots,b_n>0$ satisfy the assumptions in Proposition~\ref{prop: HJ1_hd_2}. 
Then, for any $\bx\in\R^n$ and $t>0$, the unique optimal trajectory of the optimal control problem~\eqref{eqt: result_optctrl1_hd_general} is given by the function $[0,t]\ni s\mapsto \opttrajhd(s;\bx,t)\in\R^n$ defined in~\eqref{eqt: optctrlhd_traj_general}. Moreover, the optimal value of the optimal control problem~\eqref{eqt: result_optctrl1_hd_general} equals $\valuefn(\bx,t)$ as defined in~\eqref{eqt: result_HJ1_Hopf_hd_general}.
\end{prop}

\begin{proof}
Let $\bx\in\Rn$ and $t>0$. First, we prove that the problem~\eqref{eqt: result_optctrl1_hd_general} is equivalent to another optimal control problem in the form of~\eqref{eqt: result_optctrl1_hd}. 
For any trajectory $s\mapsto\bx(s)$ satisfying the constraint in~\eqref{eqt: result_optctrl1_hd_general}, define another trajectory $\by(s):= \matP^{-1}(\bx(s) - \uz)$.
Note that the trajectory $s\mapsto\by(s)$ satisfies the constraint in the following optimal control problem:
\begin{equation}\label{eqt: prop28_equiv_optctrl}
    \inf \left\{\int_0^t \left(\frac{1}{2}\|\dot{\by}(s)\|^2 - \sum_{i=1}^n \potentfn_i(\by(s))\right) ds + \tilde{\initcond}(\by(0)) \colon \by(t) = \matP^{-1}(\bx - \uz)\right\}.
\end{equation}
Now, we prove that the cost of $\bx(\cdot)$ in~\eqref{eqt: result_optctrl1_hd_general} equals the cost of $\by(\cdot)$ in~\eqref{eqt: prop28_equiv_optctrl}.
By straightforward calculation, the cost of $\bx(\cdot)$ in the problem~\eqref{eqt: result_optctrl1_hd_general} equals
\begin{equation}\label{eqt:prop26_pf_cost1}
\begin{split}
    &\int_0^t \left(\frac{1}{2}\|\dot{\bx}(s)\|_{M^{-1}}^2 - \potentfn(\bx(s))\right) ds + \initcond(\bx(0))\\
    =\,& \int_0^t \left(\frac{1}{2}\|\matP\dot{\by}(s)\|_{M^{-1}}^2 - \potentfn(\matP\by(s) +\uz) \right) ds + \initcond(\matP\by(0) +\uz)\\
    =\,& \int_0^t \left(\frac{1}{2}\|\matP\dot{\by}(s)\|_{M^{-1}}^2 - \potentfn(\matP\by(s) +\uz) )\right) ds + \tilde{\initcond}(\by(0)),
\end{split}
\end{equation}
where the last equality holds by definition of $\tilde{\initcond}$ in~\eqref{eqt: J_changeofvariable}.
Since $M = \matP \matP^T$ holds and $\matP$ is invertible, we have
\begin{equation}\label{eqt:prop26_pf_cost2}
\begin{split}
    \|\matP\dot{\by}(s)\|_{M^{-1}}^2 &= \dot{\by}(s)^T \matP^T M^{-1}\matP\dot{\by}(s) = \dot{\by}(s)^T \matP^T (\matP^T)^{-1}\matP^{-1}\matP\dot{\by}(s) \\
    &= \dot{\by}(s)^T\dot{\by}(s) = \|\dot{\by}(s)\|^2.
\end{split}
\end{equation}
By assumption, the function $\by\mapsto \potentfn(\matP\by +\uz)$ is non-positive, concave, and 1-homogeneous and its superlevel set at the value $-1$ is calculated as follows:
\begin{equation*}
\begin{split}
    \{\by\in\Rn\colon \potentfn(\matP\by +\uz)\geq -1\}
    &= \matP^{-1}\left(\co \left(\bigcup_{j=1}^n \left\{\frac{1}{a_j}\bu_j, -\frac{1}{b_j}\bu_j\right\}\right)\right)\\
    &= \co \left(\bigcup_{j=1}^n \left\{\frac{1}{a_j}\be_j, -\frac{1}{b_j}\be_j\right\}\right),
\end{split}
\end{equation*}
where the first equality holds by~\eqref{eqt: prop24_levelset_V} and the second equality follows from straightforward calculation (recall that $\be_1,\dots,\be_n\in\Rn$ denote the standard basis vectors in $\Rn$).
As a result, we have 
\begin{equation}\label{eqt:prop26_pf_cost3}
\potentfn(\matP\by +\uz) = \sum_{i=1}^n \potentfn_i(y_i) \quad\forall \by=(y_1,\dots,y_n)\in\Rn,
\end{equation}
where each $\potentfn_i$ is the function defined in~\eqref{eqt: result_defV} with parameters $a = a_i$ and $b = b_i$.

Combining~\eqref{eqt:prop26_pf_cost1},~\eqref{eqt:prop26_pf_cost2}, and~\eqref{eqt:prop26_pf_cost3}, we conclude that the cost of $\bx(\cdot)$ in~\eqref{eqt: result_optctrl1_hd_general} equals the cost of $\by(\cdot)$ in~\eqref{eqt: prop28_equiv_optctrl}.
Moreover, this relation between the trajectory $\bx(\cdot)$ and the trajectory $\by(\cdot)$ is a bijection. Hence, the two problems~\eqref{eqt: result_optctrl1_hd_general} and~\eqref{eqt: prop28_equiv_optctrl} are equivalent. By Proposition~\ref{prop: optctrl1_hd} (whose assumptions are checked in the proof of Proposition~\ref{prop: HJ1_hd_2}), the unique optimal trajectory of the problem~\eqref{eqt: prop28_equiv_optctrl} is $s\mapsto \by(s):= \tilde{\opttrajhd}(s;\matP^{-1}(\bx-\uz),t)$ as defined in~\eqref{eqt: optctrl1_hd_general_opty}, and the optimal value of the problem~\eqref{eqt: prop28_equiv_optctrl} equals $\valuefn(\bx,t)$ in~\eqref{eqt: result_HJ1_Hopf_hd_general}. Therefore, the unique optimal trajectory of the problem~\eqref{eqt: result_optctrl1_hd_general} is the corresponding trajectory $\bx(\cdot)$ given by $\bx(s)=\matP\by(s)+\uz = \matP\tilde{\opttrajhd}(s; \matP^{-1}(\bx-\bu_0),t) + \bu_0 = \opttrajhd(s; \bx,t)$, and the optimal value of the problem~\eqref{eqt: result_optctrl1_hd_general} also equals $\valuefn(\bx,t)$.
\end{proof}

\begin{rem}
If the initial cost $\initcond$ is a linear function given by $\initcond(\bx)=\langle \bp,\bx\rangle$ for all $\bx\in\Rn$ and for some constant vector $\bp=(p_1,\dots,p_n)\in\Rn$, then the solution $\valuefn$ to the corresponding HJ PDE~\eqref{eqt: HJhd_1_general} becomes
\begin{equation*}
    \valuefn(\bx,t) = \sum_{i=1}^n \valuefn\left((\matP^{-1}\bx-\matP^{-1}\uz)_i, t; p_i, a_i,b_i\right)\quad 
    \forall\, \bx\in\Rn, t\geq 0,
\end{equation*}
where the function $\valuefn$ in the summation on the right-hand side is defined by~\eqref{eqt: result_S1_1d} and~\eqref{eqt: result_S1_1d_negative}. 
The optimal trajectory of the optimal control problem~\eqref{eqt: result_optctrl1_hd_general} is defined by~\eqref{eqt: optctrlhd_traj_general} and~\eqref{eqt: optctrl1_hd_general_opty}, where the point $\bp^*$ equals $\bp$, which is the slope of the linear initial cost $\initcond$.
\end{rem}
\subsection{An extension to certain non-convex initial costs} \label{sec:HJhd_minplus}

In this section, we solve the high-dimensional problems~\eqref{eqt: result_optctrl1_hd_general} and~\eqref{eqt: HJhd_1_general} for certain non-convex initial data $\initcond$. To be specific, we assume the function $\initcond\colon \Rn\to\R$ satisfies
\begin{equation} \label{eqt: J_minplus}
    \initcond(\bx) := \min_{j\in\{1,\dots, m\}} \initcond_j(\bx) \quad\forall \bx\in\Rn,
\end{equation}
where $\initcond_1,\dots,\initcond_m\colon\Rn\to\R$ are convex functions satisfying~\eqref{eqt: HJ1_hd_condJ}. 

The solution $\valuefn\colon \Rn\times [0,+\infty)\to\R$ is defined by 
\begin{equation}\label{eqt:Hopf_minplus}
\begin{split}
    \valuefn(\bx,t) = \min_{j\in\{1,\dots,m\}} \valuefn_{\initcond_j}(\bx,t) \quad \forall\, \bx\in\Rn,\, t\geq 0,
\end{split}
\end{equation}
where for each $j\in\{1,\dots,m\}$, the function $\valuefn_{\initcond_j}$ on the right-hand side is the solution defined by~\eqref{eqt: result_HJ1_Hopf_hd_general} (or~\eqref{eqt: result_Hopf1_1d} and~\eqref{eqt: result_Hopf1_hd} for special cases) with the initial data $\initcond_j$. Similarly, the optimal trajectory $[0,t]\ni s\mapsto \opttrajhd(s;\bx,t)\in\R^n$ is defined by
\begin{equation}\label{eqt:opttraj_minplus}
    \opttrajhd(s;\bx,t):= \opttrajhd_{\initcond_r}(s;\bx,t) \quad\forall s\in[0,t],\quad\text{ where } r\in \argmin_{j\in\{1,\dots,m\}} \valuefn_{\initcond_j}(\bx,t),
\end{equation}
and the function $s\mapsto \opttrajhd_{\initcond_r}(s;\bx,t)$ is the trajectory defined by~\eqref{eqt: optctrlhd_traj_general} (or~\eqref{eqt: optctrl_traj_J} and~\eqref{eqt: optctrlhd_traj_J} for special cases) with the initial cost $\initcond_r$.

In the following proposition, we prove that the function $\valuefn$ defined above is the viscosity solution to the HJ PDE~\eqref{eqt: HJhd_1_general} and the value function of the optimal control problem~\eqref{eqt: result_optctrl1_hd_general} with initial data $\initcond$. Moreover, we show that the trajectory $s\mapsto\opttrajhd(s;\bx,t)$ defined above is an optimal trajectory of the optimal control problem~\eqref{eqt: result_optctrl1_hd_general} with initial cost $\initcond$.

\begin{prop} \label{prop:HJ_minplus}
Let $\initcond\colon \Rn\to\R$ be a continuous function of the form of~\eqref{eqt: J_minplus} for some convex functions $\initcond_1,\dots,\initcond_m\colon\Rn\to\R$ satisfying~\eqref{eqt: HJ1_hd_condJ}. Assume the function $\initcond$ is bounded below by an affine function. Let the function $\potentfn\colon\R^n\to(-\infty,0]$, the matrices $\matP$ and $M$, the vector $\uz\in\Rn$, and the scalars $a_1,\dots, a_n,b_1,\dots,b_n>0$ satisfy the assumptions in Proposition~\ref{prop: HJ1_hd_2}. Let $\valuefn$ be the function defined in~\eqref{eqt:Hopf_minplus} and $\opttrajhd$ be the trajectory defined in~\eqref{eqt:opttraj_minplus}.
Then, the following statements hold:
\begin{itemize}
    \item[(a)] The function $\valuefn$ is the unique viscosity solution in the solution set $\mathcal{G}$ in~\eqref{eqt: HJ1_hd_solset_G} to the HJ PDE~\eqref{eqt: HJhd_1_general} with initial data $\initcond$.
    
    \item[(b)] Let $\bx$ be an arbitrary vector in $\Rn$ and $t>0$ be an arbitrary positive number. The trajectory $s\mapsto\opttrajhd(s;\bx,t)$ is an optimal trajectory of the optimal control problem~\eqref{eqt: result_optctrl1_hd_general} with initial cost $\initcond$. Moreover, the optimal value of the optimal control problem~\eqref{eqt: result_optctrl1_hd_general} equals $\valuefn(\bx,t)$.
\end{itemize}
\end{prop}

\begin{proof} 
We prove (b) first.
Since each $\initcond_j$ is a convex function satisfying~\eqref{eqt: HJ1_hd_condJ}, by Proposition~\ref{prop: HJ1_opttraj_general}, the value $\valuefn_{\initcond_j}(\bx,t)$ is the optimal value of the optimal control problem~\eqref{eqt: result_optctrl1_hd_general} with initial cost $\initcond_j$. In other words, we have that
\begin{equation*}
    \valuefn_{\initcond_j}(\bx,t) = \inf \left\{\int_0^t \left(\frac{1}{2}\|\dot{\bx}(s)\|_{M^{-1}}^2 - \potentfn(\bx(s))\right) ds + \initcond_j(\bx(0)) \colon \bx(t) = \bx\right\},
\end{equation*}
for any $ j\in\{1,\dots,m\}$.
After some calculations, we obtain
\begin{equation*}
\begin{split}
    \valuefn(\bx,t) &= \min_{j\in\{1,\dots,m\}} \valuefn_{\initcond_j}(\bx,t)\\
    &= \min_{j\in\{1,\dots,m\}}\inf \left\{\int_0^t \left(\frac{1}{2}\|\dot{\bx}(s)\|_{M^{-1}}^2 - \potentfn(\bx(s))\right) ds + \initcond_j(\bx(0)) \colon \bx(t) = \bx\right\}\\
    &= \inf \left\{\int_0^t \left(\frac{1}{2}\|\dot{\bx}(s)\|_{M^{-1}}^2 - \potentfn(\bx(s))\right) ds + \min_{j\in\{1,\dots,m\}}\initcond_j(\bx(0)) \colon \bx(t) = \bx\right\}\\
    &= \inf \left\{\int_0^t \left(\frac{1}{2}\|\dot{\bx}(s)\|_{M^{-1}}^2 - \potentfn(\bx(s))\right) ds + \initcond(\bx(0)) \colon \bx(t) = \bx\right\}.
\end{split}
\end{equation*}
Therefore, $\valuefn(\bx,t)$ is the optimal value of the problem~\eqref{eqt: result_optctrl1_hd_general} with initial cost $\initcond$. 

Now, we show that the trajectory $s\mapsto\opttrajhd(s;\bx,t)$ is an optimal trajectory. We abuse notation and use $\opttrajhd(s;\bx,t)$ and $\opttrajhd(s)$ interchangeably whenever there is no ambiguity.
Let $r$ be the index in~\eqref{eqt:opttraj_minplus}. By Proposition~\ref{prop: HJ1_opttraj_general} and~\eqref{eqt:opttraj_minplus}, $\opttrajhd$ is the optimal trajectory of the problem~\eqref{eqt: result_optctrl1_hd_general} with initial cost $\initcond_r$. As a result, its cost equals the optimal value $\valuefn_{\initcond_r}(\bx,t)$, which equals $\valuefn(\bx,t)$ since $r$ is a minimizer in~\eqref{eqt: J_minplus}. Hence, we get
\begin{equation*}
\begin{split}
    \valuefn(\bx,t) &= \int_0^t \left(\frac{1}{2}\|\dot{\opttrajhd}(s)\|_{M^{-1}}^2 - \potentfn(\opttrajhd(s))\right) ds + \initcond_r(\opttrajhd(0))\\
    &\geq \int_0^t \left(\frac{1}{2}\|\dot{\opttrajhd}(s)\|_{M^{-1}}^2 - \potentfn(\opttrajhd(s))\right) ds + \min_{j\in\{1,\dots,m\}}\initcond_j(\opttrajhd(0))
    \geq \valuefn(\bx,t).
\end{split}
\end{equation*}
Therefore, the inequalities above are equalities. In other words, the cost of $\opttrajhd$ in the optimal control problem~\eqref{eqt: result_optctrl1_hd_general} with initial cost $\initcond$ equals the optimal value $\valuefn(\bx,t)$, and hence, $\opttrajhd$ is an optimal trajectory of~\eqref{eqt: result_optctrl1_hd_general} with initial cost $\initcond$.

It remains to prove (a). We will prove (a) by applying Lemma~\ref{lem:A10_uniqueness}. We need to check the assumptions in Lemma~\ref{lem:A10_uniqueness}. Note that each $\initcond_j$ satisfies~\eqref{eqt: HJ1_hd_condJ} for some constants $C_j,\delta_j\geq 0$. Choose $C:=3\max_{j\in\{1,\dots,m\}}C_j\geq 0$ and $\delta:=\max_{j\in\{1,\dots,m\}}\delta_j\geq 0$. For each $j\in\{1,\dots,m\}$, we have
\begin{equation*}
\begin{split}
    |\initcond_j(\bz) - \initcond_j(\by)| &\leq C_j\|\bz-\by\|(1 + \|\bz\|^{\delta_j} + \|\by\|^{\delta_j})\\
    &\leq C_j\|\bz-\by\|\left(1 + (1+\|\bz\|^{\delta}) + (1+\|\by\|^{\delta})\right)\\
    &\leq 3C_j\|\bz-\by\|(1 + \|\bz\|^{\delta} + \|\by\|^{\delta})\\
    &\leq C\|\bz-\by\|(1 + \|\bz\|^{\delta} + \|\by\|^{\delta})
    \quad\quad \forall \bz,\by\in\R^n.
\end{split}
\end{equation*}
For any $\bz,\by\in\Rn$, letting $k\in\argmin_{j\in\{1,\dots,m\}}\initcond_j(\by)$, there holds
\begin{equation*}
    \initcond(\bz)-\initcond(\by)= \initcond(\bz) - \initcond_k(\by)\leq \initcond_k(\bz)- \initcond_k(\by) \leq C\|\bz-\by\|(1 + \|\bz\|^{\delta} + \|\by\|^{\delta}).
\end{equation*}
Similarly, letting $k\in\argmin_{j\in\{1,\dots,m\}}\initcond_j(\bz)$, we have 
\begin{equation*}
    \initcond(\bz)-\initcond(\by)= \initcond_k(\bz) - \initcond(\by)\geq \initcond_k(\bz)- \initcond_k(\by) \geq -C\|\bz-\by\|(1 + \|\bz\|^{\delta} + \|\by\|^{\delta}).
\end{equation*}
Therefore, the function $\initcond$ also satisfies~\eqref{eqt: HJ1_hd_condJ}. Recall that by assumption, $\initcond$ is continuous and bounded below by an affine function. Thus, the assumptions of Lemma~\ref{lem:A10_uniqueness} on $\initcond$ hold. The assumptions of Lemma~\ref{lem:A10_uniqueness} on $\potentfn$ and $M$ also hold by straightforward reasoning. It remains to show that the assumptions on $\valuefn$ in Lemma~\ref{lem:A10_uniqueness}(b) hold. Let $\balp\in\Rn$ be a vector such that $\bx\mapsto \initcond(\bx)-\langle \balp, \bx\rangle $ is bounded from below, and denote the lower bound by $\beta\in\R$. By~\eqref{eqt: J_minplus}, for each $j\in\{1,\dots,m\}$, there holds
\begin{equation*}
    \initcond_j(\bx) - \langle \balp, \bx\rangle\geq \initcond(\bx)-\langle \balp, \bx\rangle \geq \beta \quad \forall \bx\in\Rn,
\end{equation*}
which implies that $\balp$ is in the domain of $\initcond_j^*$.
Recall that the function $\valuefn_{\initcond_j}$ is defined by~\eqref{eqt: result_HJ1_Hopf_hd_general} with initial data $\initcond_j$, and hence, $\valuefn_{\initcond_j}$ satisfies~\eqref{eqt:Hopf_hd_general_ver2} with initial condition $\initcond_j$. 
By Lemma~\ref{lem: appendix_S_bd}, we have
\begin{equation*}
\begin{split}
    \valuefn_{\initcond_j}(\bx,t) &\geq \sum_{i=1}^n \valuefn\left((\matP^{-1}\bx-\matP^{-1}\uz)_i, t; (\matP^T\balp)_i, a_i,b_i\right) - \initcond_j^*(\balp) + \langle \balp, \uz\rangle \\
    &\geq \sum_{i=1}^n \left((\matP^{-1}\bx-\matP^{-1}\uz)_i (\matP^T\balp)_i + C_i\right) - \initcond_j^*(\balp) + \langle \balp, \uz\rangle \\
    &= \langle \matP^{-1}\bx-\matP^{-1}\uz, \matP^T\balp\rangle + \sum_{i=1}^n C_i - \initcond_j^*(\balp) + \langle \balp, \uz\rangle\\
    &= \langle \balp,\bx\rangle + \sum_{i=1}^n C_i - \initcond_j^*(\balp),
\end{split}
\end{equation*}
where each $C_i$ is the constant in the lower bound in Lemma~\ref{lem: appendix_S_bd} with constants $a=a_i$, $b=b_i$, and $p = (\matP^T\balp)_i$.
Then, by~\eqref{eqt:Hopf_minplus}, there holds
\begin{equation*}
\begin{split}
    \valuefn(\bx,t) &= \min_{j\in\{1,\dots,m\}}\valuefn_{\initcond_j}(\bx,t)\geq \min_{j\in\{1,\dots,m\}}\left\{\langle \balp,\bx\rangle + \sum_{i=1}^n C_i - \initcond_j^*(\balp)\right\}\\
    &= \langle \balp,\bx\rangle + \sum_{i=1}^n C_i - \max_{j\in\{1,\dots,m\}}\initcond_j^*(\balp).
\end{split}
\end{equation*}
Since $\initcond_j^*(\balp)$ is a finite number for each $j\in\{1,\dots,m\}$, the function $(\bx,t)\mapsto \valuefn(\bx,t) - \langle \balp,\bx\rangle$ is bounded below by $\sum_{i=1}^n C_i - \max_{j\in\{1,\dots,m\}}\initcond_j^*(\balp)\in\R$, and hence, the assumption on $\valuefn$ in Lemma~\ref{lem:A10_uniqueness}(b) is satisfied. Therefore, all the assumptions in Lemma~\ref{lem:A10_uniqueness}(a)-(b) are satisfied. Applying Lemma~\ref{lem:A10_uniqueness}, we get that the value function in~\eqref{eqt: result_optctrl1_hd_general} is the unique viscosity solution to the HJ PDE~\eqref{eqt: HJhd_1_general} in the solution set $\mathcal{G}$. Moreover, by (b), which we proved earlier, the function $\valuefn$ is the value function in~\eqref{eqt: result_optctrl1_hd_general}. Therefore, the function $\valuefn$ is the unique viscosity solution to the HJ PDE~\eqref{eqt: HJhd_1_general} in the solution set $\mathcal{G}$.
\end{proof}

\begin{rem}
The technique used in~\eqref{eqt:Hopf_minplus} is called the min-plus technique (or max-plus if the optimal control problem is formulated as a maximization problem). For more details, see~\cite{McEneaney2006maxplus,Kolokoltsov1997Idempotent}, for instance. The min-plus technique can also be applied to solve the problems in Sections~\ref{sec:HJ1d} or~\ref{sec:HJhd_separable}. The unique viscosity solution and the optimal value are given by~\eqref{eqt:Hopf_minplus}, where each $\valuefn_{\initcond_j}$ is the solution to the corresponding HJ PDE in Sections~\ref{sec:HJ1d} or~\ref{sec:HJhd_separable} with initial data $\initcond_j$. An optimal trajectory is given by~\eqref{eqt:opttraj_minplus}, where the trajectory $s\mapsto\opttrajhd_{\initcond_r}(s;\bx,t)$ is the optimal trajectory of the corresponding optimal control problem in Sections~\ref{sec:HJ1d} or~\ref{sec:HJhd_separable} with initial cost $\initcond_r$.
\end{rem}

\begin{rem}\label{rmk:minplus_nonunique}
Note that the minimizer $r$ in~\eqref{eqt:opttraj_minplus} may be not unique. Whenever there are multiple minimizers, each minimizer $r$ gives an optimal trajectory $s\mapsto\opttrajhd_{\initcond_r}(s;\bx,t)$. The optimal trajectory of the optimal control problem~\eqref{eqt: result_optctrl1_hd_general} may be non-unique since~\eqref{eqt: result_optctrl1_hd_general} is not a convex optimization problem in this case (the initial cost $\initcond$ is non-convex).
\end{rem}
\def\admmb{w}
\def\admmbb{\boldsymbol{\admmb}}
\def\Sanaly{\valuefn}
\def\gmanaly{\opttrajhd}
\def\Snum{\hat{\valuefn}}
\def\gmnum{\hat{\opttrajhd}}

\section{Efficient numerical algorithms}\label{sec:admm}

In this section, we present efficient numerical solvers based on some optimization methods that evaluate the optimal trajectory of the high-dimensional optimal control problem \eqref{eqt: result_optctrl1_hd} as well as the solution of the corresponding high-dimensional HJ PDE \eqref{eqt: result_HJ1_hd}. 
We note that the algorithms we present in this section can be extended to solve \eqref{eqt: result_optctrl1_hd_general} and \eqref{eqt: HJhd_1_general}, instead. To solve the more general problems in~\eqref{eqt: result_optctrl1_hd_general} and~\eqref{eqt: HJhd_1_general}, we apply our algorithms to compute the optimizer $\bp^*$ with the terminal position $\bx$ replaced by $\matP^{-1}\bx - \matP^{-1}\bu_0$ and the Legendre transform of the initial cost $\initcond^*$ replaced by $\tilde \initcond^*$ as defined by \eqref{eqt: J_star_changeofvariable}. Then, we compute the optimal values and optimal trajectories using~\eqref{eqt: result_HJ1_Hopf_hd_general} and~\eqref{eqt: optctrlhd_traj_general}, respectively.

Recall that in Section~\ref{sec: theory}, we provided representation formulas for the problems \eqref{eqt: result_optctrl1_hd} and \eqref{eqt: result_HJ1_hd}. Thus, we can numerically solve these problems using these representation formulas if the optimization problem in \eqref{eqt: result_Hopf1_hd} is numerically solvable. In this section, we provide different methods to solve \eqref{eqt: result_Hopf1_hd} for different classes of initial costs $\initcond$. More specifically, in Section~\ref{subsec:numerical_quad}, we present efficient numerical solvers for quadratic initial costs based on explicit formulas for solving \eqref{eqt: result_Hopf1_hd} exactly.
In Section~\ref{sec:numerical_convex}, we solve \eqref{eqt: result_Hopf1_hd} with more general convex initial costs using optimization methods that utilize the numerical solver presented in Section~\ref{subsec:numerical_quad} as a building block.
For illustrative purposes, the ADMM algorithm (see~\cite{Glowinski2014Alternating,Boyd2011Distributed}) is applied. However, we note that ADMM can be replaced by any other appropriate convex optimization algorithm. In Section~\ref{sec: ADMM_nonconvex}, we extend our numerical methods to address the class of non-convex initial costs that are of the form of~\eqref{eqt: J_minplus}. In each of these three sections, we also present high-dimensional numerical examples and timing results, which demonstrate the efficiency of our proposed methods in each of these cases.
All of the numerical examples in Sections~\ref{subsec:numerical_quad}-~\ref{sec: ADMM_nonconvex} are run using a C++ implementation on an 8th Gen Intel Laptop Core i5-8250U with a 1.60GHz processor. Finally, in Section \ref{sec:fpga}, we describe a high throughput FPGA implementation of our building block from Section \ref{subsec:numerical_quad} and present some numerical results demonstrating the performance boost that can be obtained using FPGAs. 

To avoid confusion, we use $\valuefn$ and $\opttrajhd$ to denote the analytical solutions to the HJ PDEs and optimal control problems, while we use $\Snum$ and $\gmnum$ to denote their numerical approximations as obtained by our proposed methods.
\subsection{Quadratic initial costs} \label{subsec:numerical_quad}
In this section, we solve the optimal control problem~\eqref{eqt: result_optctrl1_hd} and the HJ PDE~\eqref{eqt: result_HJ1_hd} with quadratic initial cost defined by
\begin{equation*}
    \initcond(\bx) = \frac{1}{2\lambda}\|\bx - \by\|^2 + \alpha \quad \forall \bx\in\Rn,
\end{equation*}
where $\by\in\Rn$, $\lambda>0$, and $\alpha\in\R$ are some parameters. Recall that $\|\cdot\|$ denotes the $\ell^2$-norm in $\Rn$. 
To solve these problems, we need to solve the optimization problem in~\eqref{eqt: result_Hopf1_hd}. Then, the solution is given by the explicit formulas~\eqref{eqt: result_Hopf1_hd} and~\eqref{eqt: optctrlhd_traj_J}.
By straightforward computation, the Legendre-Fenchel transform of $\initcond$ is
\begin{equation*}
    \initcond^*(\bp) = \frac{\lambda}{2}\left\|\bp + \frac{\by}{\lambda}\right\|^2 - \frac{\|\by\|^2}{2\lambda} - \alpha \quad \forall \bp\in\Rn.
\end{equation*}
Therefore, for quadratic initial costs, the optimization problem in~\eqref{eqt: result_Hopf1_hd} is equivalent to 
\begin{equation}\label{eqt:opt_quad_hd}
    \min_{\bp\in\R^n} \left\{\sum_{i=1}^n\left(- \valuefn(x_i, t; p_i, a_i,b_i) + \frac{\lambda}{2}\left(p_i + \frac{y_i}{\lambda}\right)^2\right)\right\}.
\end{equation}
Note that the optimization problem~\eqref{eqt:opt_quad_hd} can be divided into $n$ one-dimensional subproblems since each term in the summation, which corresponds to each state dimension, is independent from each other. The $i$th subproblem amounts to computing
\begin{equation}\label{eqt: proxptofminusV}
    p_i^* = \argmin_{p\in\R} \left\{- \valuefn(x_i, t; p, a_i,b_i) + \frac{\lambda}{2}\left(p + \frac{y_i}{\lambda}\right)^2\right\},
\end{equation}
which can be calculated using the numerical solver described in Appendix~\ref{sec:appendix_numerical_prox} with parameters $x=x_i$, $a=a_i$, $b=b_i,$ and $c =-\frac{y_i}{\lambda}$.
Thus, solving~\eqref{eqt:opt_quad_hd} is embarrassingly parallel. We also note that the minimizer in~\eqref{eqt: proxptofminusV} is the proximal point of $p\mapsto -\frac{1}{\lambda}\valuefn(x_i, t; p, a_i,b_i)$ at $-\frac{y_i}{\lambda}$. This relation suggests that our proposed numerical solver for quadratic initial costs may be a useful building block in proximal point-based methods for solving the problems with more general initial costs.

\begin{figure}[htbp]
    \centering
    \begin{subfigure}{0.49\textwidth}
        \centering \includegraphics[width=\textwidth]{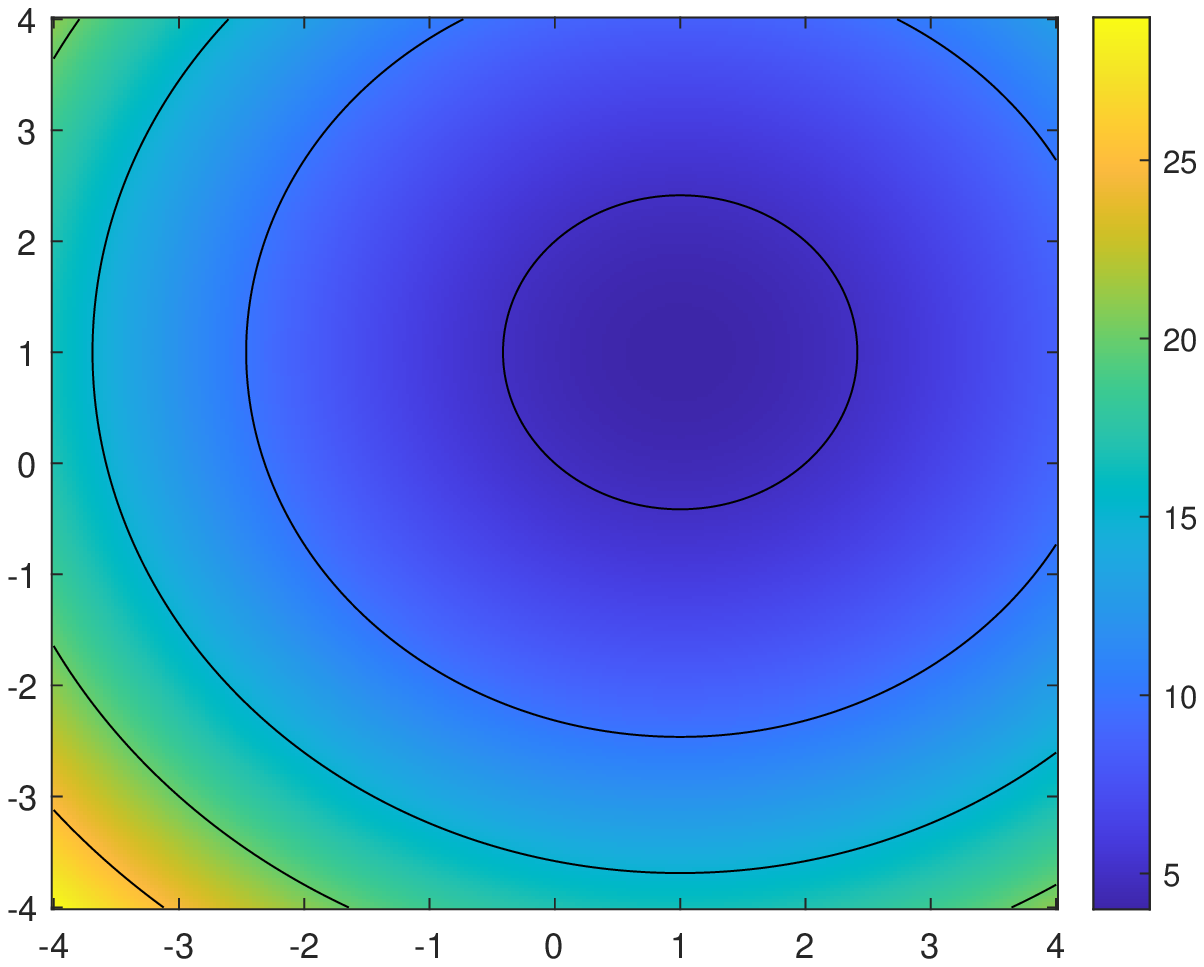}
        \caption{$t=0$}
    \end{subfigure}
    \hfill
    \begin{subfigure}{0.49\textwidth}
        \centering \includegraphics[width=\textwidth]{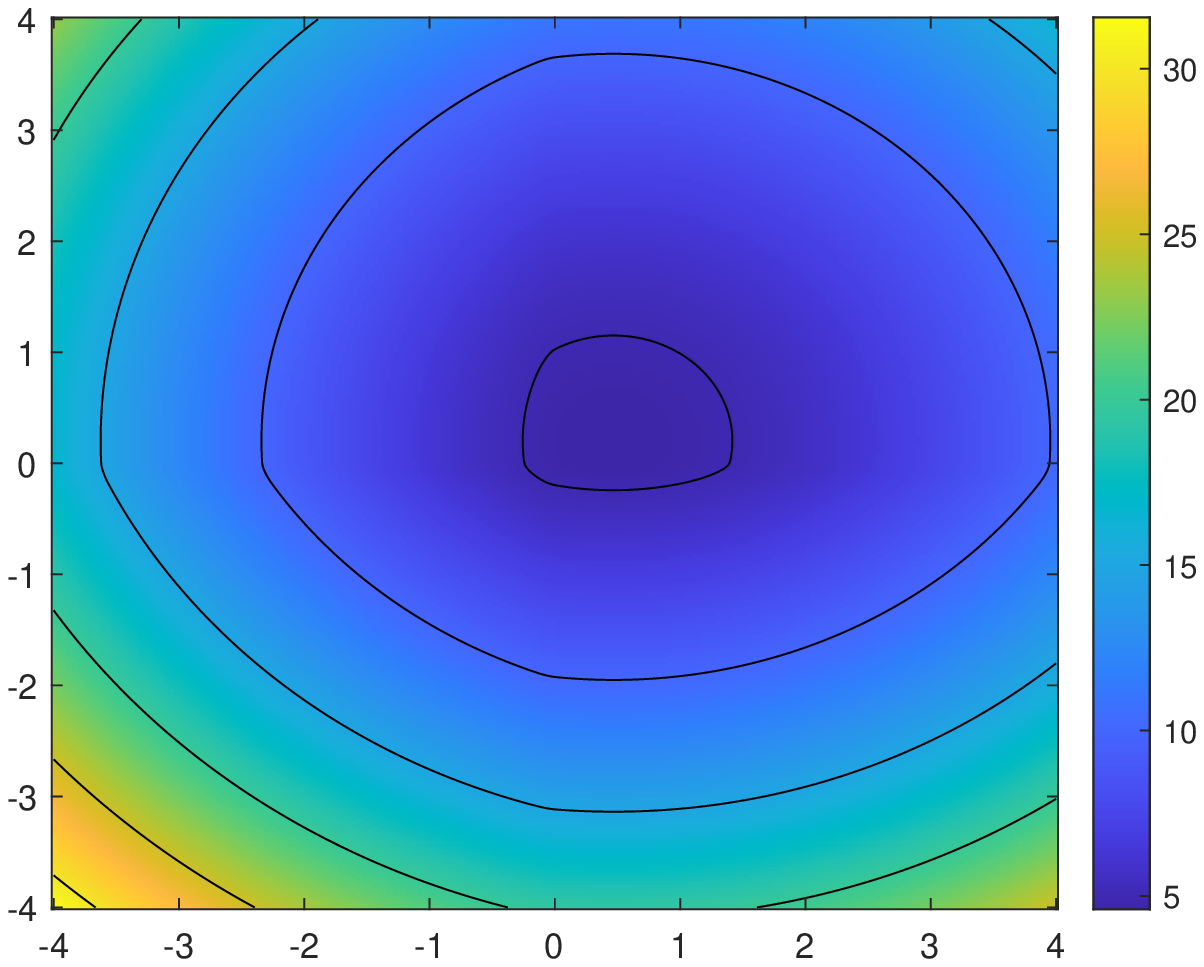}
        \caption{$t=0.125$}
    \end{subfigure}
    
    \begin{subfigure}{0.49\textwidth}
        \centering \includegraphics[width=\textwidth]{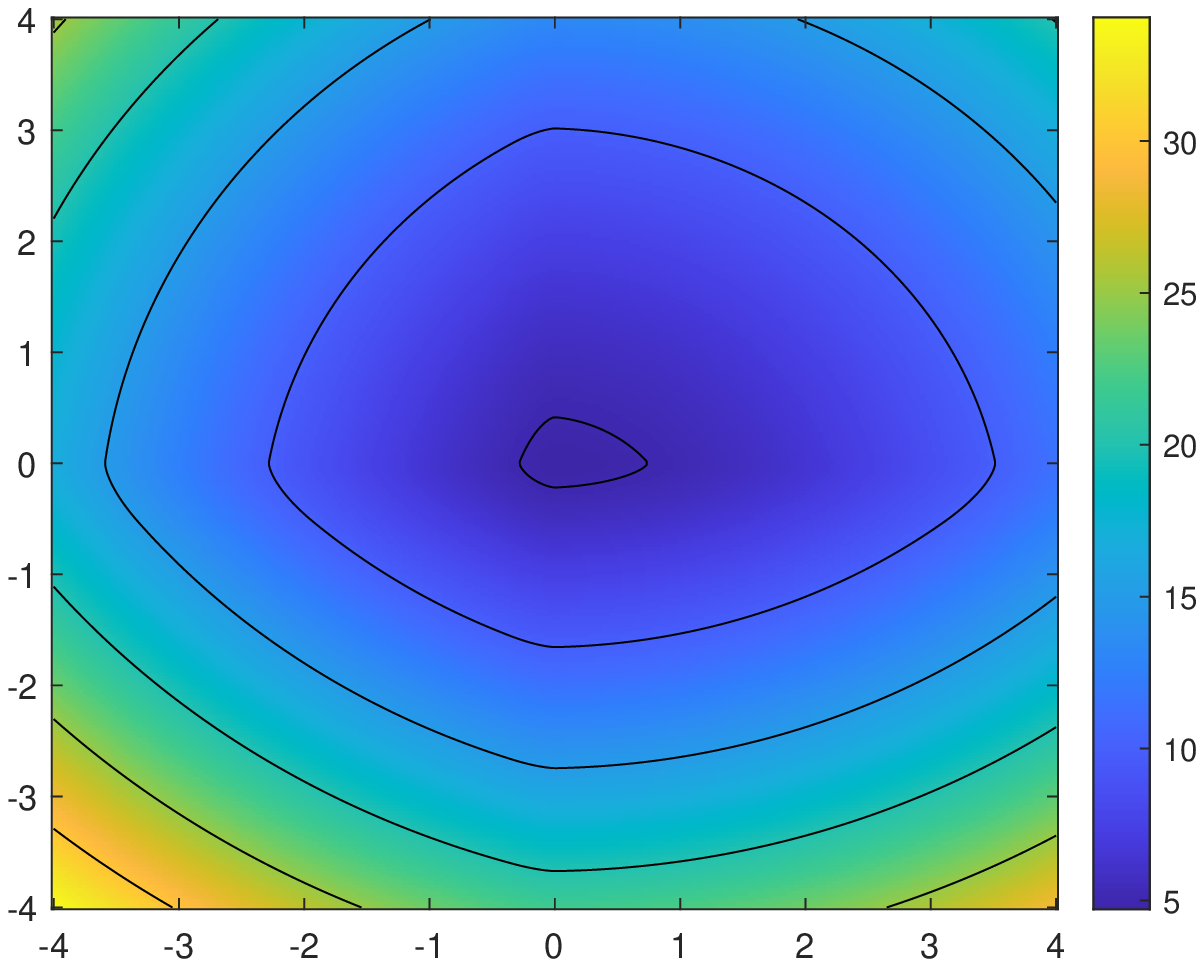}
        \caption{$t=0.25$}
    \end{subfigure}
    \hfill
    \begin{subfigure}{0.49\textwidth}
        \centering \includegraphics[width=\textwidth]{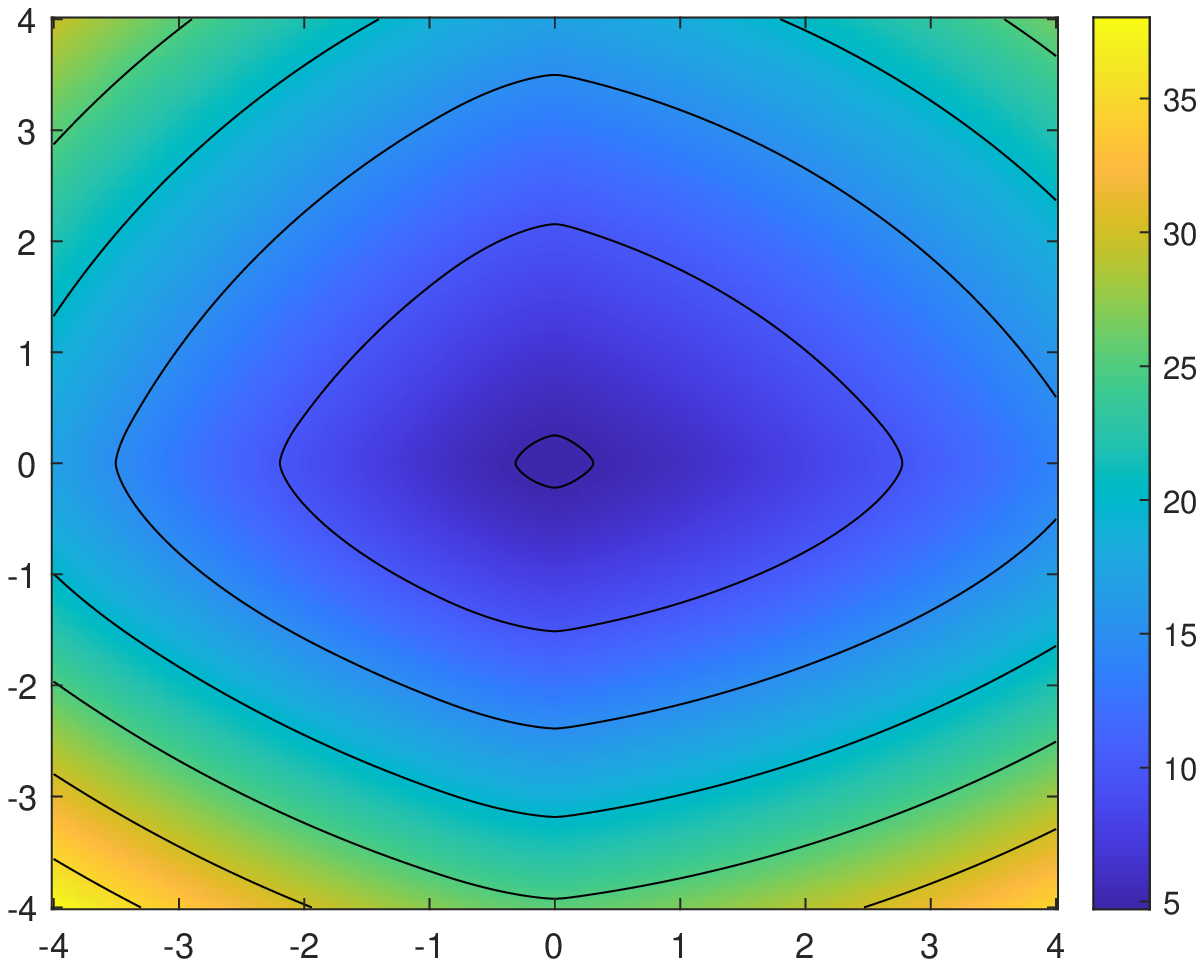}
        \caption{$t=0.5$}
    \end{subfigure}
    
    \caption{Evaluation of the solution $\Snum(\bx,t)$ of the HJ PDE (\ref{eqt: result_HJ1_hd}) with $\ba$ and $\bb$ defined in~\eqref{eqt:numerical_def_ab} and initial data $\initcond(\bx) = \frac{1}{2} \|\bx -\mathbf{1}\|^2$ for $\bx=(x_1,x_2,0,\dots,0)\in \R^{10}$ and different times $t$. Plots for $t =0$, $0.125$, $0.25$, and $0.5$ are depicted in (a)-(d), respectively. Level lines are superimposed on the plots.}
    \label{fig: HJ1_2d_quadratic_contour}
\end{figure}

Now, we apply our proposed numerical solver to solve the HJ PDE~\eqref{eqt: result_HJ1_hd} with the following quadratic initial cost:
\begin{equation}\label{eqt: quadratic_IC}
    \initcond(\bx) = \frac{1}{2}\|\bx - \mathbf{1}\|^2 \quad \forall \bx\in\Rn,
\end{equation}
i.e., we set $\lambda = 1$, $\alpha = 0$, and $\by=\mathbf{1}$, where $\mathbf{1}$ denotes the vector in $\R^n$ whose elements are all one. We also define the parameters $\ba=(a_1,\dots, a_n)\in \R^n$ and $\bb = (b_1,\dots, b_n)\in \R^{n}$ by
\begin{equation} \label{eqt:numerical_def_ab}
    a_i = \begin{dcases}
    4 & \text{if } i=1,\\
    6 &\text{if } i=2,\\  
    5 &\text{if } i>2,
    \end{dcases}
    \quad\quad\text{ and }\quad\quad 
    b_i = \begin{dcases}
    3 & \text{if } i=1,\\
    9 &\text{if } i=2,\\  
    6 &\text{if } i>2.
    \end{dcases}
\end{equation}
In Figure~\ref{fig: HJ1_2d_quadratic_contour}, we show the numerical solution to the HJ PDE~\eqref{eqt: result_HJ1_hd} in dimension $n=10$. 
More specifically, Figure~\ref{fig: HJ1_2d_quadratic_contour} depicts two-dimensional contour plots of the solution $\Snum(\bx,t)$ 
for $\bx=(x_1,x_2,0,\dots,0)$ and different times $t$. 
The running time for this example in different dimensions is shown in Table~\ref{tab:timing_quad}. To compute the running time, we first compute the overall running time for computing the solution at $102,400$ random points $(\bx,t)\in [-4,4]^n\times [0,0.5]$ and then report the average running time for computing the solution $\Snum(\bx,t)$ at one point $(\bx,t)$ over these $102,400$ trials.
From Table~\ref{tab:timing_quad}, we see that on average, it takes less than $2\times 10^{-6}$ seconds to compute the solution at one point in a $16$-dimensional problem, which demonstrates the efficiency of our proposed solver even in high dimensions.

\begin{table}[htbp]
    \centering
    \begin{tabular}{c|c|c|c|c}
        \hline
        \textbf{$\mathbf{n}$} & 4 & 8 & 12 & 16 \\
        \hline
        \textbf{running time (s)} & 4.2330e-07 & 9.2605e-07 & 1.4404e-06 & 1.9732e-06  \\
        \hline
    \end{tabular}
    \hfill 
    
    \caption{Time results in seconds for the average time per call over $102,400$ trials for evaluating the solution of the HJ PDE (\ref{eqt: result_HJ1_hd}) with quadratic initial condition~\eqref{eqt: quadratic_IC} for various dimensions $n$.}
    \label{tab:timing_quad}
\end{table}

\subsection{Convex initial costs}\label{sec:numerical_convex}

\begin{algorithm}[htbp]
\SetAlgoLined
\SetKwInOut{Input}{Inputs}
\SetKwInOut{Output}{Outputs}
\Input{Parameters in the problem: $\ba,\bb\in\Rn$, terminal position $\bx\in\Rn$, time horizon $t>0$, running time $s>0$ of the trajectory, and Legendre transform $\initcond^*$ of the convex initial cost $\initcond$. Parameters for ADMM: $\lambda>0$, initialization $\bd^0\in\Rn$, $\admmbb^0\in\Rn$, error tolerance $\epsilon>0$.}
\Output{The optimal trajectory $\gmnum(s;\bx,t)$ in the optimal control problem~\eqref{eqt: result_optctrl1_hd} and the solution value $\Snum(\bx,t)$ to the corresponding HJ PDE~\eqref{eqt: result_HJ1_hd}.}
 \For{$k = 1,2,\dots$}{
 Update $\bv^{k+1}\in\Rn$ by 
 \begin{equation}\label{eqt: ADMM_HJ1_vupdate}
    \bv^{k+1} = \argmin_{\bv\in \Rn} \left\{ \initcond^*(\bv) + \frac{\lambda}{2} \left\| \bv - \bd^k + \admmbb^k\right\|^2\right\}.
 \end{equation}
 \\
 Update $\bd^{k+1}\in\Rn$, where the $i$-th element $d_i^{k+1}$ is updated by 
 \begin{equation}\label{eqt: ADMM_HJ1_dupdate}
    d_i^{k+1} = \argmin_{d_i\in \R} \left\{ -\valuefn(x_i, t; d_i, a_i,b_i) + \frac{\lambda}{2} (v_i^{k+1} - d_i + \admmb_i^k)^2\right\}.
 \end{equation}
 \\
 Update $\admmbb^{k+1}\in\Rn$ by 
 \begin{equation*}
     \admmbb^{k+1} = \admmbb^k + \bv^{k+1} - \bd^{k+1}.
 \end{equation*}
 \\
 \If{$\|\bv^{k+1} - \bv^{k}\|^2 \leq \epsilon$, $\|\bd^{k+1} - \bd^{k}\|^2 \leq \epsilon$, and $\|\bv^{k+1} - \bd^{k+1}\|^2 \leq \epsilon$}{
   set $N=k+1$ and $\initposhd^N = \bv^N$\;
   break\;
   }
 }
 Output the optimal trajectory by
 \begin{equation}\label{eqt: ADMM_HJ1_opttraj}
     \gmnum(s;\bx,t) = (\opttraj(s; x_1,t,\initpos_1^N,a_1,b_1), \dots, \opttraj(s; x_n,t,\initpos_n^N,a_n,b_n)), 
 \end{equation}
 where the $i$-th component $\opttraj(s; x_i,t,\initpos_i^N, a_i,b_i)$ is defined in~\eqref{eqt: optctrl_defx_1},~\eqref{eqt: optctrl_defx_2},~\eqref{eqt: optctrl_defx_3},~\eqref{eqt: optctrl_defx_4},~\eqref{eqt: optctrl_defx_5}, and~\eqref{eqt: optctrl_defx_neg}.
 Also, output the solution to the HJ PDE by
 \begin{equation}\label{eqt: ADMM_HJ1_solution}
     \Snum(\bx,t)  = \sum_{i=1}^n \valuefn\left(x_i,t; \initpos_i^N, a_i,b_i\right) - \initcond^*\left(\initposhd^N\right),
 \end{equation}
 where the $i$-th component $\valuefn\left(x_i,t; \initpos_i^N, a_i,b_i\right)$ in the summation is defined in~\eqref{eqt: result_S1_1d},~\eqref{eqt: domain_subregions},~\eqref{eqt: case1_deff}, and~\eqref{eqt: result_S1_1d_negative}.
 \caption{An ADMM algorithm for solving the optimal control problem~\eqref{eqt: result_optctrl1_hd} and the corresponding HJ PDE~\eqref{eqt: result_HJ1_hd} with convex initial cost. \label{alg:admm_ver1}}
\end{algorithm}

In this section, we solve the optimal control problem \eqref{eqt: result_optctrl1_hd} and the corresponding HJ PDE \eqref{eqt: result_HJ1_hd} with convex initial costs. To solve these problems, we need to solve the optimization problem in the representation formula~\eqref{eqt: result_Hopf1_hd}, which can be rewritten as
\begin{equation}\label{eqt: ADMM_Hopfhd_rewritten}
    \valuefn(\bx,t) = -\inf_{\bp\in\R^n} \left\{-\sum_{i=1}^n \valuefn(x_i, t; p_i, a_i,b_i) + \initcond^*(\bp)\right\}.
\end{equation}  
By Lemma~\ref{lem: prop_S}, if $\initcond$ is convex and $t > 0$, then the optimization problem in \eqref{eqt: ADMM_Hopfhd_rewritten} is a convex optimization problem with strictly convex, 1-coercive objective function. Thus, the optimal value in \eqref{eqt: ADMM_Hopfhd_rewritten} is finite, and the minimizer exists and is unique. Furthermore, since the objective function is convex, \eqref{eqt: ADMM_Hopfhd_rewritten} can be solved numerically using convex optimization algorithms. Following the discussion in Section \ref{subsec:numerical_quad}, proximal point-based methods would be a reasonable approach. For illustrative purposes, we demonstrate how ADMM can be applied to solve \eqref{eqt: ADMM_Hopfhd_rewritten} when $\initcond^*$ has numerically-computable proximal point. The details of applying ADMM to this problem is described in Algorithm \ref{alg:admm_ver1}. 

In each iteration of ADMM in Algorithm \ref{alg:admm_ver1}, we first update $\bv^{k+1}$ using \eqref{eqt: ADMM_HJ1_vupdate}. By definition, $\bv^{k+1}$ is the proximal point of the function $\frac{\initcond^*}{\lambda}$ at the point $\bd^k-\admmbb^k$, which we denote by $\prox_{\frac{\initcond^*}{\lambda}}(\bd^k-\admmbb^k)$. In some cases, this proximal point may not be easy to compute, but the proximal point of $\bx\mapsto \frac{1}{\lambda}\initcond(\lambda \bx)$ is easy to compute. In these cases, we can use Moreau's identity~\cite{Moreau1965Proximite} to obtain $\bv^{k+1}$ as follows:
\begin{equation}\label{eqt: ADMM_vupdate_moreau}
\begin{split}
    \bv^{k+1} &= \bd^k-\admmbb^k - \prox_{\bx\mapsto \frac{1}{\lambda}\initcond(\lambda \bx)} (\bd^k-\admmbb^k)\\
    &= \bd^k-\admmbb^k - \argmin_{\bv\in\Rn}\left\{ \initcond(\lambda\bv) + \frac{\lambda}{2} \left\| \bv - \bd^k + \admmbb^k\right\|^2\right\}.
\end{split}
\end{equation}
Thus, we compute $\bv^{k+1}$ either via \eqref{eqt: ADMM_HJ1_vupdate} or via \eqref{eqt: ADMM_vupdate_moreau} using any appropriate optimization algorithm, where the choice of optimization algorithm may depend on the form of $\initcond^*$ or $\initcond$, respectively.

In the second step of each iteration, we update $\bd^{k+1}$ componentwise, as in \eqref{eqt: ADMM_HJ1_dupdate}, which can be done in parallel. More specifically, we 
compute \eqref{eqt: ADMM_HJ1_dupdate} using the numerical solver in Appendix~\ref{sec:appendix_numerical_prox} with parameters $x=x_i$, $a=a_i$, $b=b_i$, and $c = v_i^{k+1} +\admmb_i^k$. 
We note that updating $\bd^{k+1}$ in Algorithm \ref{alg:admm_ver1} is equivalent to solving~\eqref{eqt: result_HJ1_hd} with quadratic initial cost $\initcond(\bx) = \frac{1}{2\lambda} \| \bx + \lambda(\bv^{k+1} + \admmbb^{k})\|^2$. Hence, the solver proposed in Section \ref{subsec:numerical_quad} serves as a building block for ADMM in Algorithm \ref{alg:admm_ver1}.

Finally, to recover the optimal trajectory in \eqref{eqt: result_optctrl1_hd} and the viscosity solution to \eqref{eqt: result_HJ1_hd}, we compute their approximations using~\eqref{eqt: ADMM_HJ1_opttraj} and~\eqref{eqt: ADMM_HJ1_solution}, respectively. In these two formulas, we use $\initposhd^N = \bv^N$ to approximate the maximizer $\initposhd^*$ in the original formulas~\eqref{eqt: result_Hopf1_hd} and~\eqref{eqt: optctrlhd_traj_J}. 

In the following proposition, we prove that our numerical solutions $\gmnum$ and $\Snum$ converge to their respective analytical solutions as the number of ADMM iterates approaches infinity.

\begin{prop}\label{prop:convergence_ADMM_convexJ}
Let $\initcond\colon\Rn\to\R$ be a convex function and $\ba,\bb$ be two vectors in $(0,+\infty)^n$. Let $\bx$ be any vector in $\Rn$ and $t>0$ be any scalar. Let $\Sanaly$ and $\gmanaly$ be the functions defined in~\eqref{eqt: result_Hopf1_hd} and~\eqref{eqt: optctrlhd_traj_J}, respectively. Let $\lambda>0$ and the initializations $\bd^0,\admmbb^0\in\Rn$ be arbitrary parameters for Algorithm~\ref{alg:admm_ver1}. Let $\Snum^N$ and $\gmnum^N$ be the output solution and trajectory, respectively, from Algorithm \ref{alg:admm_ver1} with iteration number $N$. Then, we have
\begin{equation}\label{eqt:prop31_conv}
    \lim_{N\to\infty} \Snum^N(\bx,t) = \Sanaly(\bx,t)  \text{ and }  \lim_{N\to\infty} \sup_{s\in[0,t]}\|\gmnum^N(s;\bx,t) - \gmanaly(s;\bx,t)\| = 0.
\end{equation}
\end{prop}
\begin{proof}
The proof is provided in Appendix~\ref{subsec:appendix_pf_prop31}.
\end{proof}

Now, we present two numerical results for the optimal control problem~\eqref{eqt: result_optctrl1_hd} and the HJ PDE~\eqref{eqt: result_HJ1_hd} with convex initial cost $\initcond$.
For simplicity, we set the parameters in Algorithm~\ref{alg:admm_ver1} to be $\lambda = 1$, $\epsilon = 10^{-8}$, $\bd^0 = \bx$, and $\admmbb^0 = \mathbf{0}$ in all of our numerical experiments.

\begin{figure}[htbp]
    \centering
    \begin{subfigure}{0.45\textwidth}
        \centering \includegraphics[width=\textwidth]{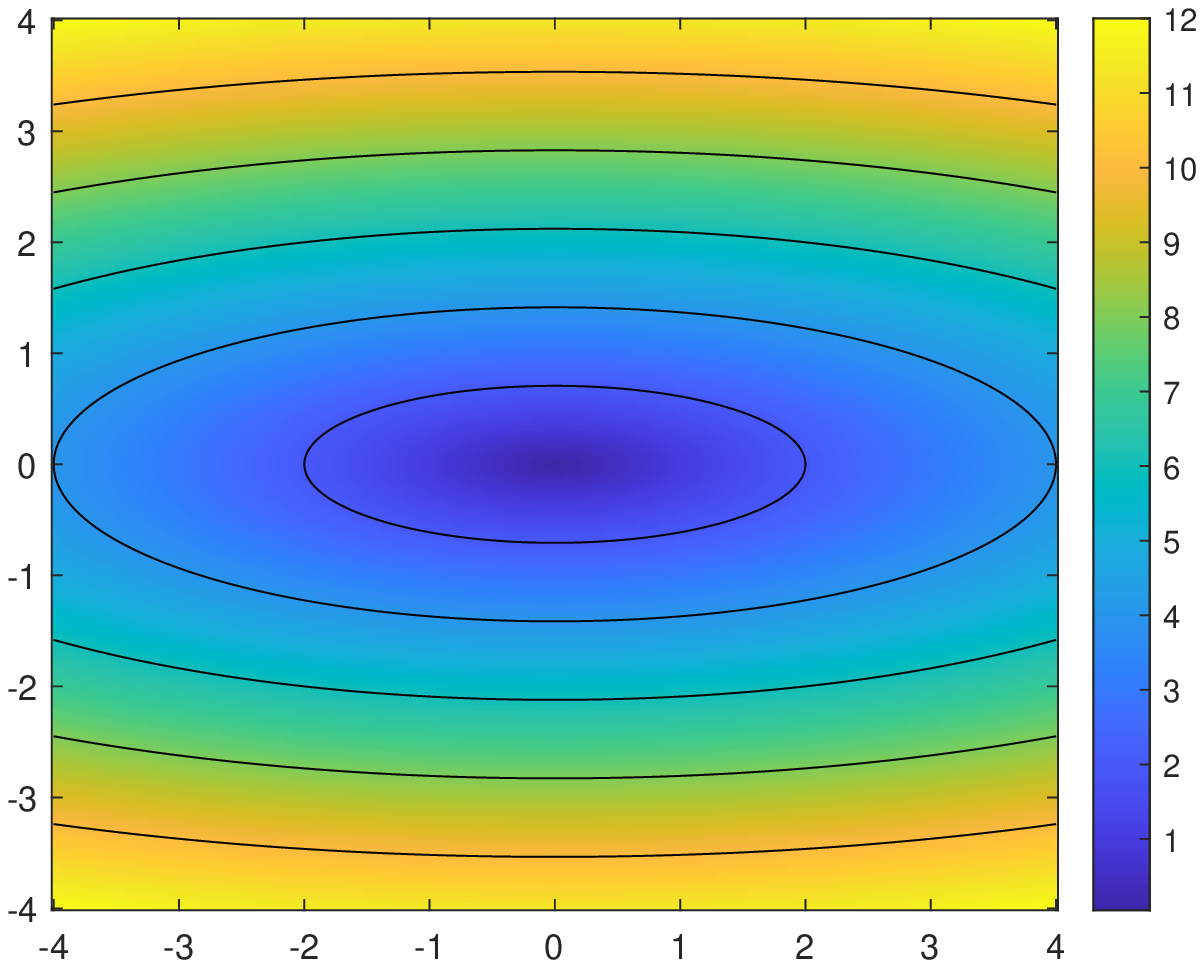}
        \caption{$t=0$}
    \end{subfigure}
    \hfill
    \begin{subfigure}{0.45\textwidth}
        \centering \includegraphics[width=\textwidth]{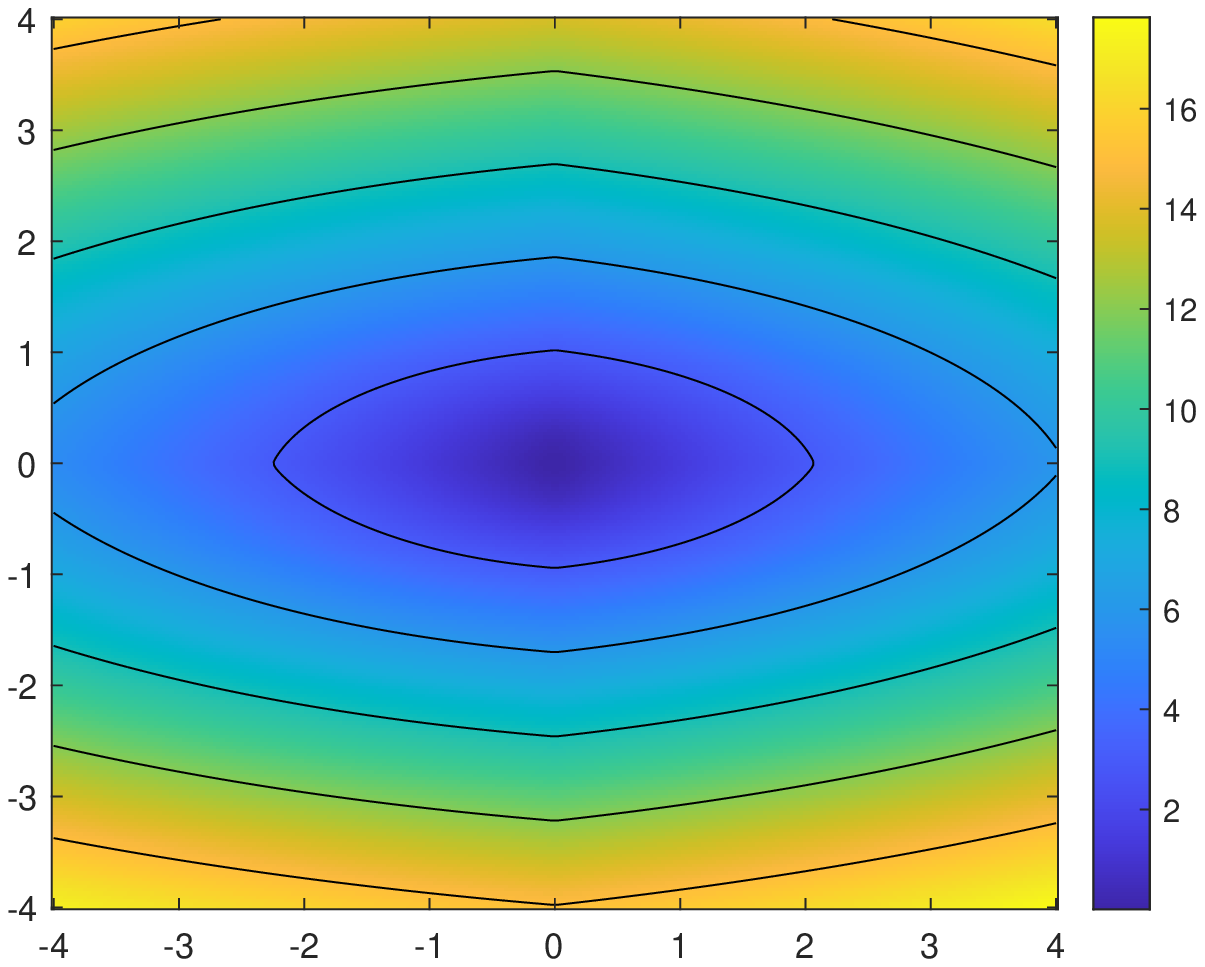}
        \caption{$t=0.125$}
    \end{subfigure}
    
    \begin{subfigure}{0.45\textwidth}
        \centering \includegraphics[width=\textwidth]{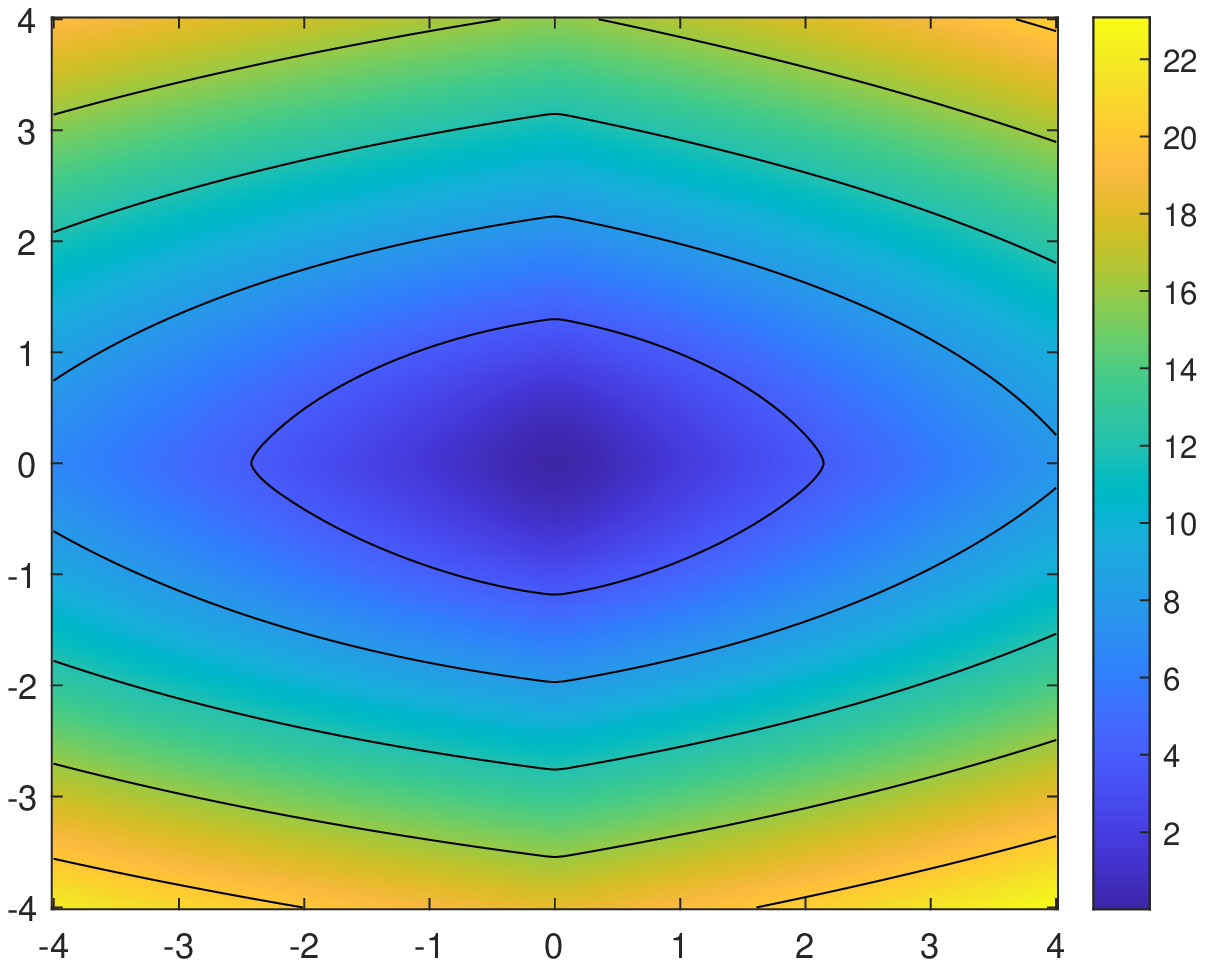}
        \caption{$t=0.25$}
    \end{subfigure}
    \hfill
    \begin{subfigure}{0.45\textwidth}
       \centering \includegraphics[width=\textwidth]{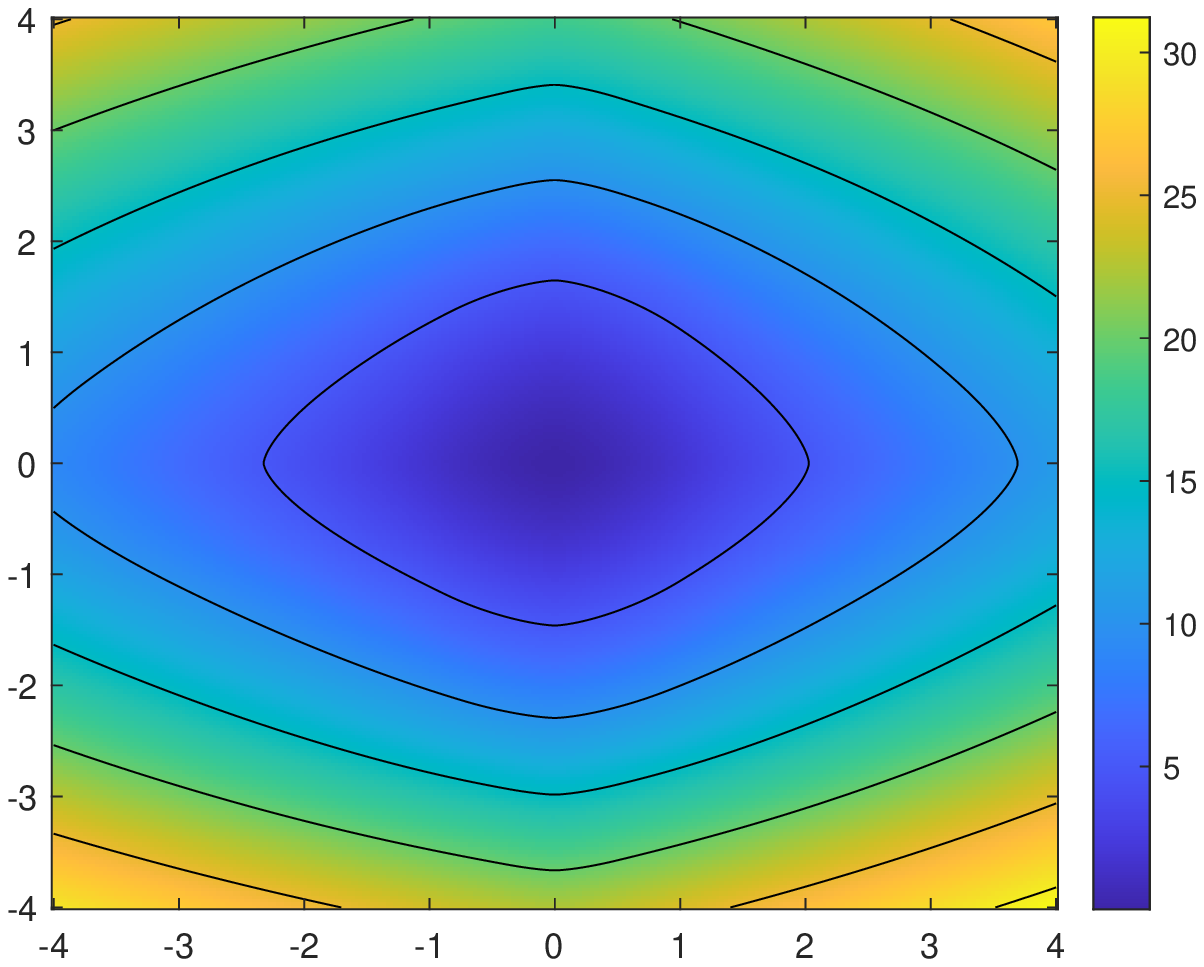}
        \caption{$t=0.5$}
    \end{subfigure}
    
    \caption{Evaluation of the solution $\Snum(\bx,t)$ of the high-dimensional HJ PDE (\ref{eqt: result_HJ1_hd}) with $\ba = ( 4, 6, 5, \dots, 5)\in\R^{10}$, $\bb=( 3, 9, 6, \dots, 6)\in\R^{10}$, and initial condition $\initcond(\bx) = \sqrt{\langle \bx, M\bx \rangle}$, where $M$ is a diagonal matrix with diagonal elements $(m_{ii}: i = 1, \dots, 10) = (1, 8, 3, 5, 1, 1,\dots, 1)$, for $\bx\in [-4,4]^2\times\{0\}^8$ and different times $t$. Plots for $t =0$, $0.125$, $0.25$, and $0.5$ are depicted in (a)-(d), respectively. Level lines are superimposed on the plots. The two axes in each figure correspond to the first and second components of the spatial variable $\bx\in\R^{10}$.}
    \label{fig: HJ1_HD_matrixnorm}
\end{figure}

We first consider the following initial cost
\begin{equation}\label{eqt: initial_cost_matrixnorm}
    \initcond(\bx) = \sqrt{\langle \bx, M \bx \rangle} \quad \forall \bx\in\Rn,
\end{equation}
where $M$ is a symmetric, positive definite matrix in $\R^{n\times n}$. Then, $\initcond^*(\bp) = \ind_{\varepsilon_M}(\bp)$, where $\ind_C$ is the indicator function defined by $\ind_C(\bx) = 0$ if $\bx\in C$ and $\ind_C(\bx) = + \infty$ otherwise and $\varepsilon_M = \{\bx\in \Rn: \langle\bx, M^{-1}\bx\rangle \leq 1\}$ is the ellipsoid associated with $M$. Thus, for this initial cost, $\bv^{k+1}$ as defined in \eqref{eqt: ADMM_HJ1_vupdate} is the projection of $\bd^k - \admmbb^k$ onto $\varepsilon_A$, which can be computed efficiently using the method described in~\cite[Section~4.3]{Darbon2016Algorithms}.
We set the parameters $\ba$ and $\bb$ to be the values defined in~\eqref{eqt:numerical_def_ab}.
For illustrative purposes, we set $M$ to be a diagonal matrix with diagonal elements $(m_{ii}: i = 1, \dots, n) = (1, 8, 3, 5, 1, 1,\dots, 1)$. Figure \ref{fig: HJ1_HD_matrixnorm} depicts two-dimensional slices of the numerical solution $\Snum(\bx,t)$ to the $10$-dimensional HJ PDE~\eqref{eqt: result_HJ1_hd} as computed using Algorithm~\ref{alg:admm_ver1}
at different positions $\bx=(x_1, x_2, 0, \dots, 0)$ and at different times $t$.

\begin{figure}[htbp]
    \centering
    \begin{subfigure}{0.32\textwidth}
        \centering \includegraphics[width=\textwidth]{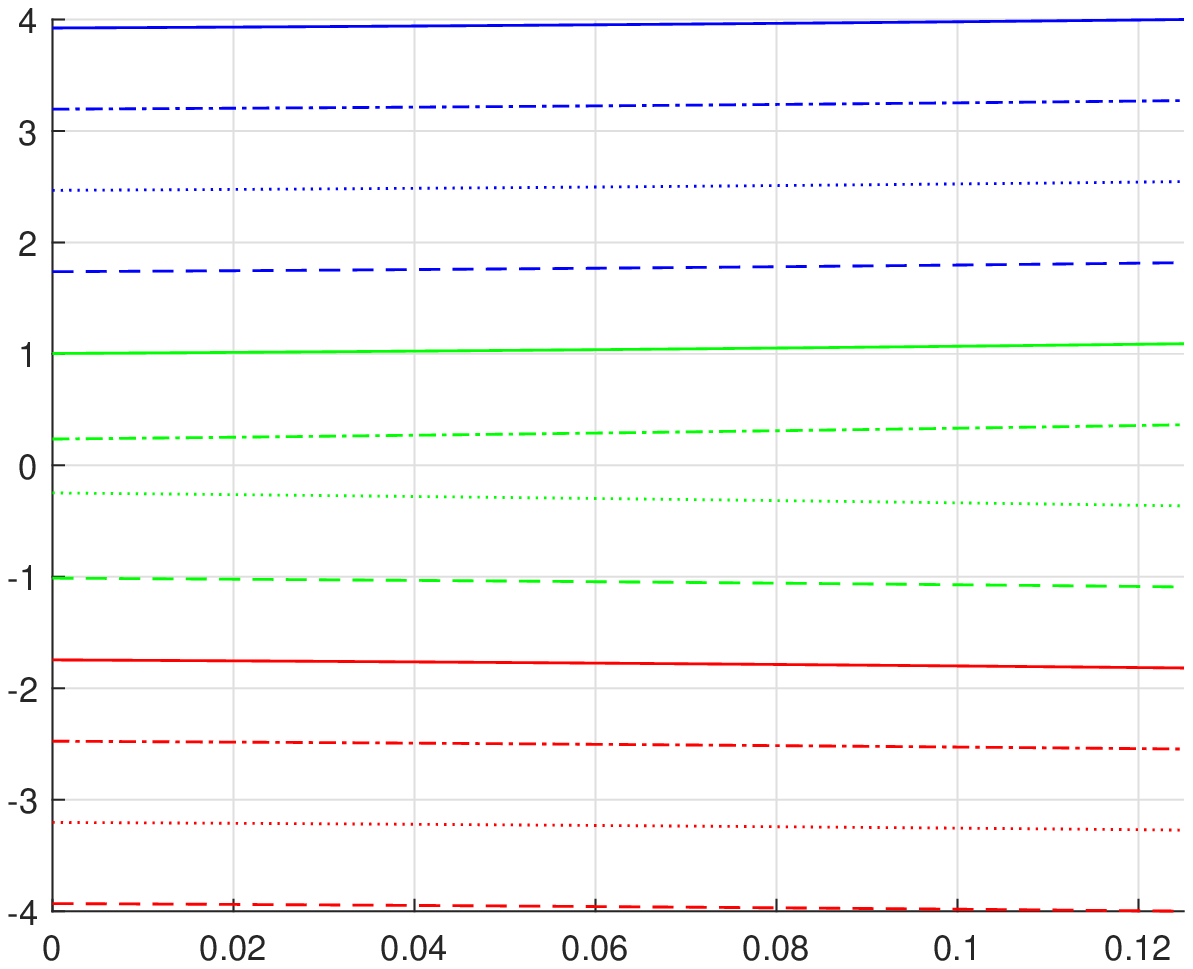}
        \caption{First component, $t=0.125$}
    \end{subfigure}
    \hfill
    \begin{subfigure}{0.32\textwidth}
        \centering \includegraphics[width=\textwidth]{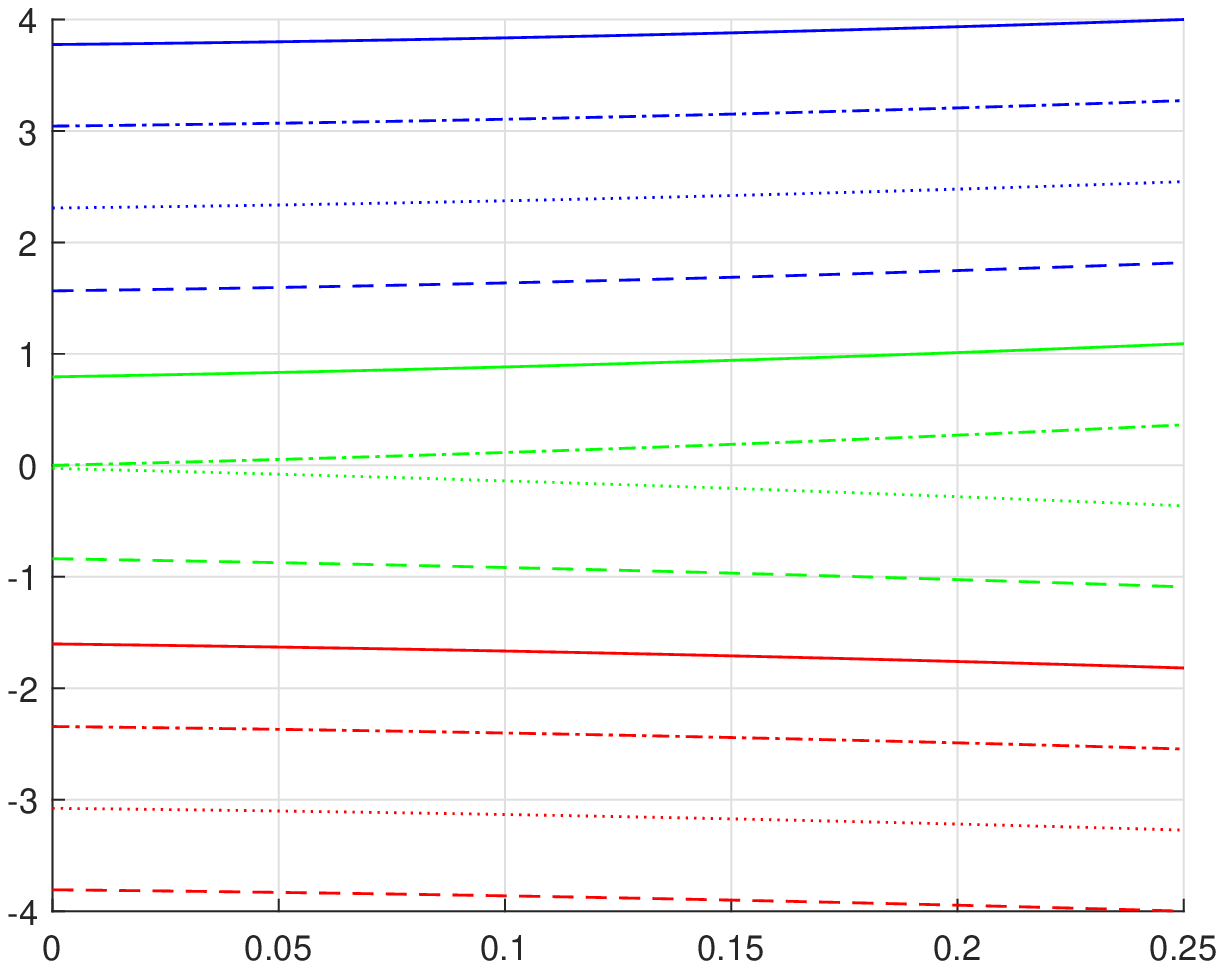}
        \caption{First component, $t=0.25$}
    \end{subfigure}
    \hfill
    \begin{subfigure}{0.32\textwidth}
        \centering \includegraphics[width=\textwidth]{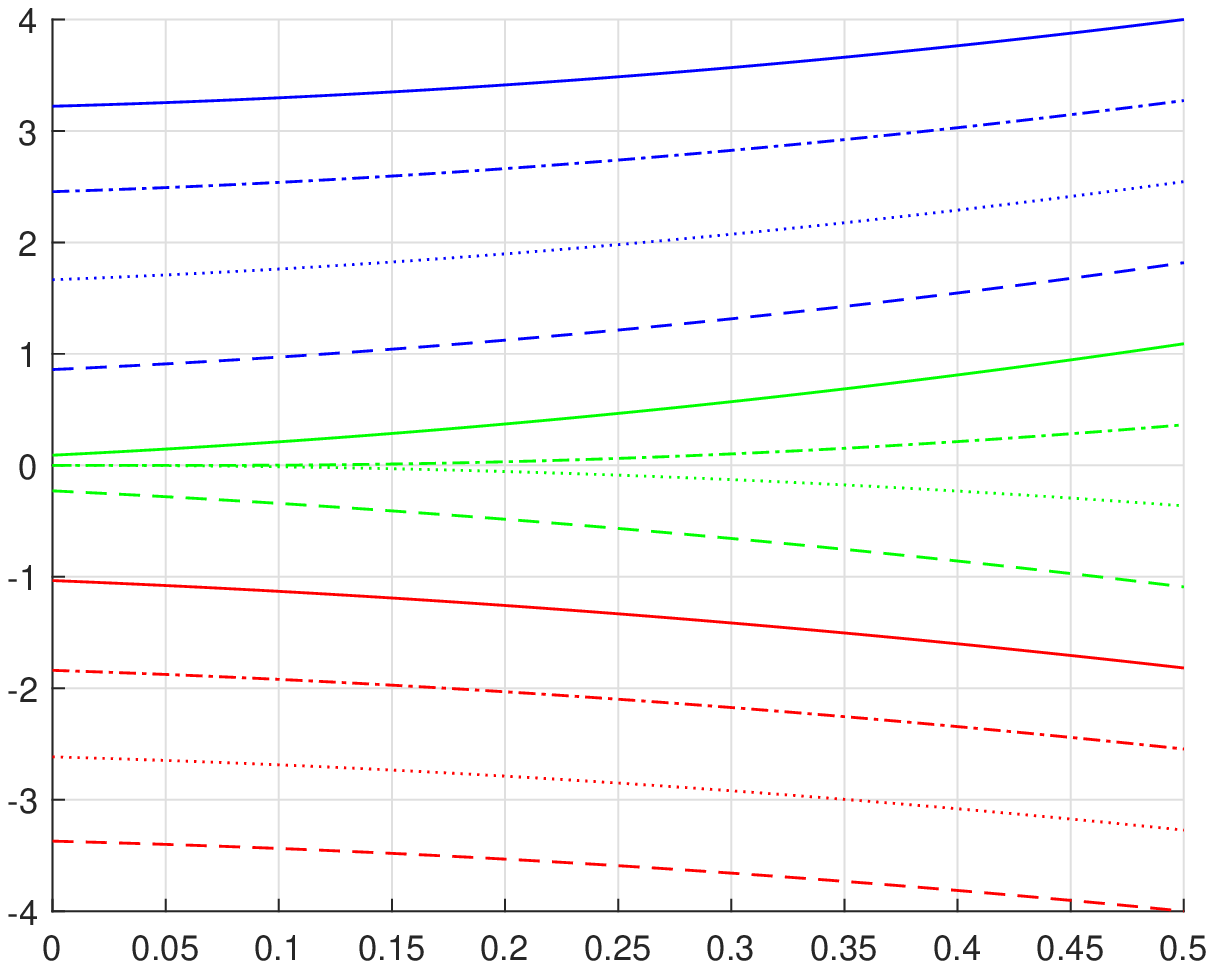}
        \caption{First component, $t=0.5$}
    \end{subfigure}
    
    \begin{subfigure}{0.32\textwidth}
        \centering \includegraphics[width=\textwidth]{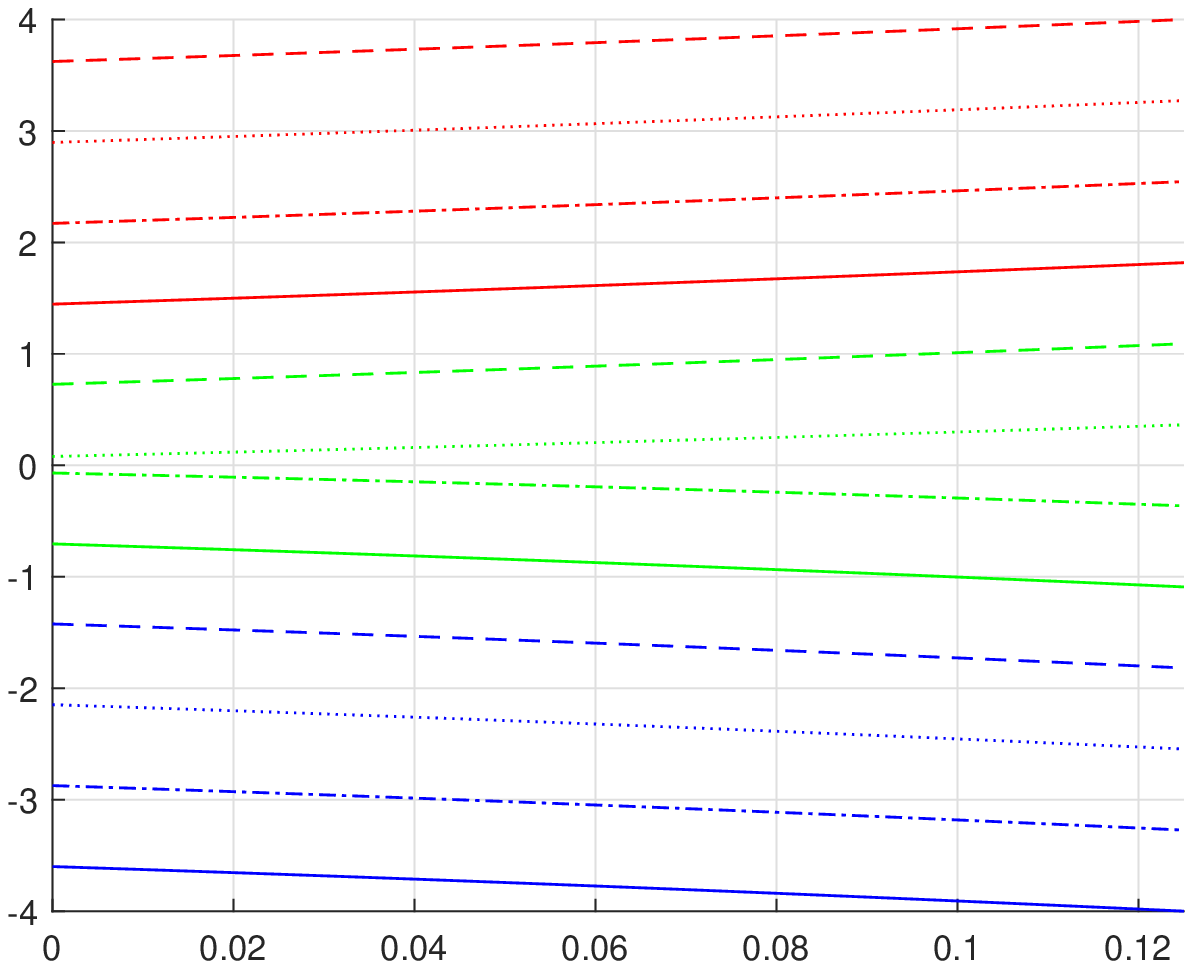}
        \caption{Second component, $t=0.125$}
    \end{subfigure}
    \hfill
    \begin{subfigure}{0.32\textwidth}
        \centering \includegraphics[width=\textwidth]{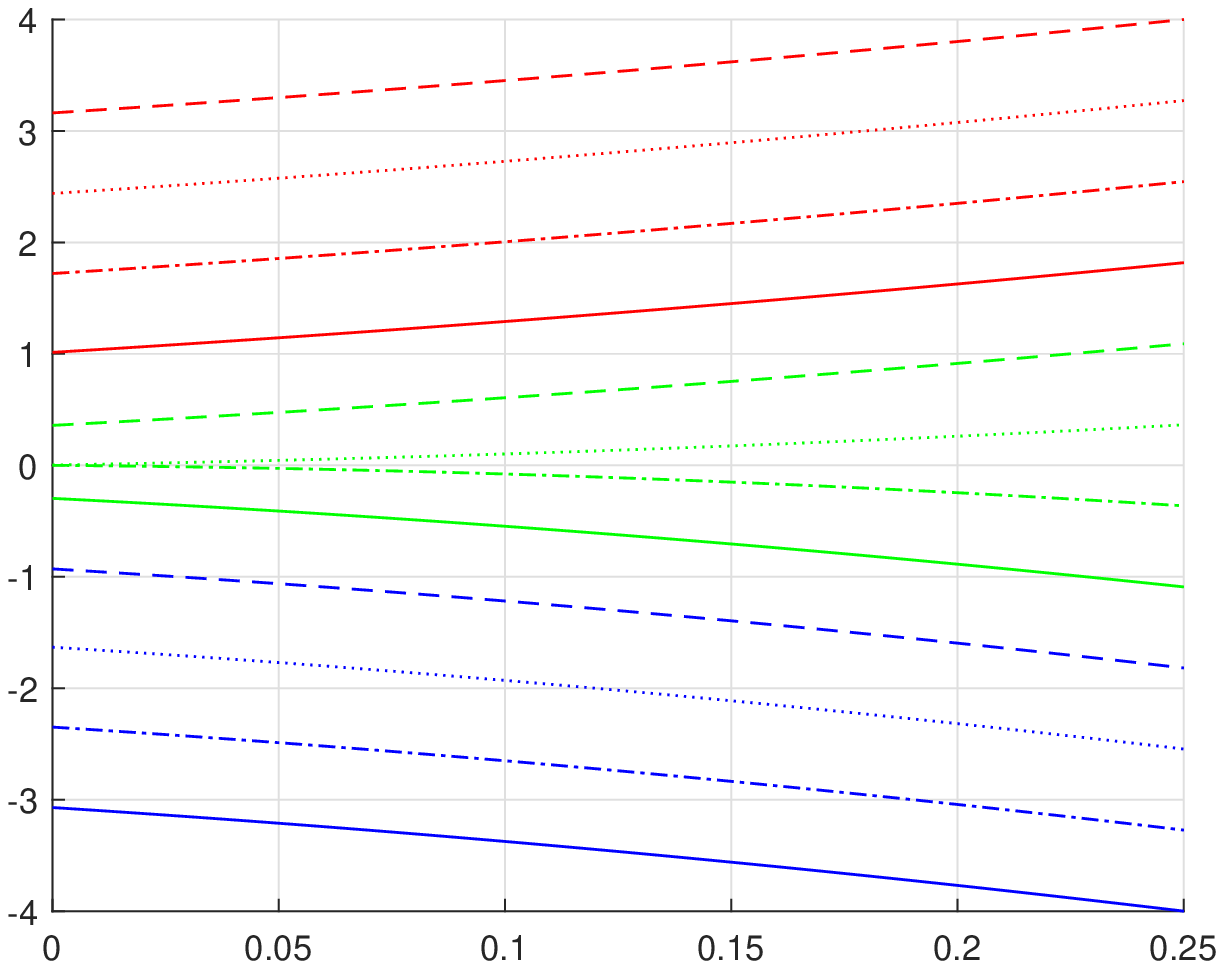}
        \caption{Second component, $t=0.25$}
    \end{subfigure}
    \hfill
    \begin{subfigure}{0.32\textwidth}
        \centering \includegraphics[width=\textwidth]{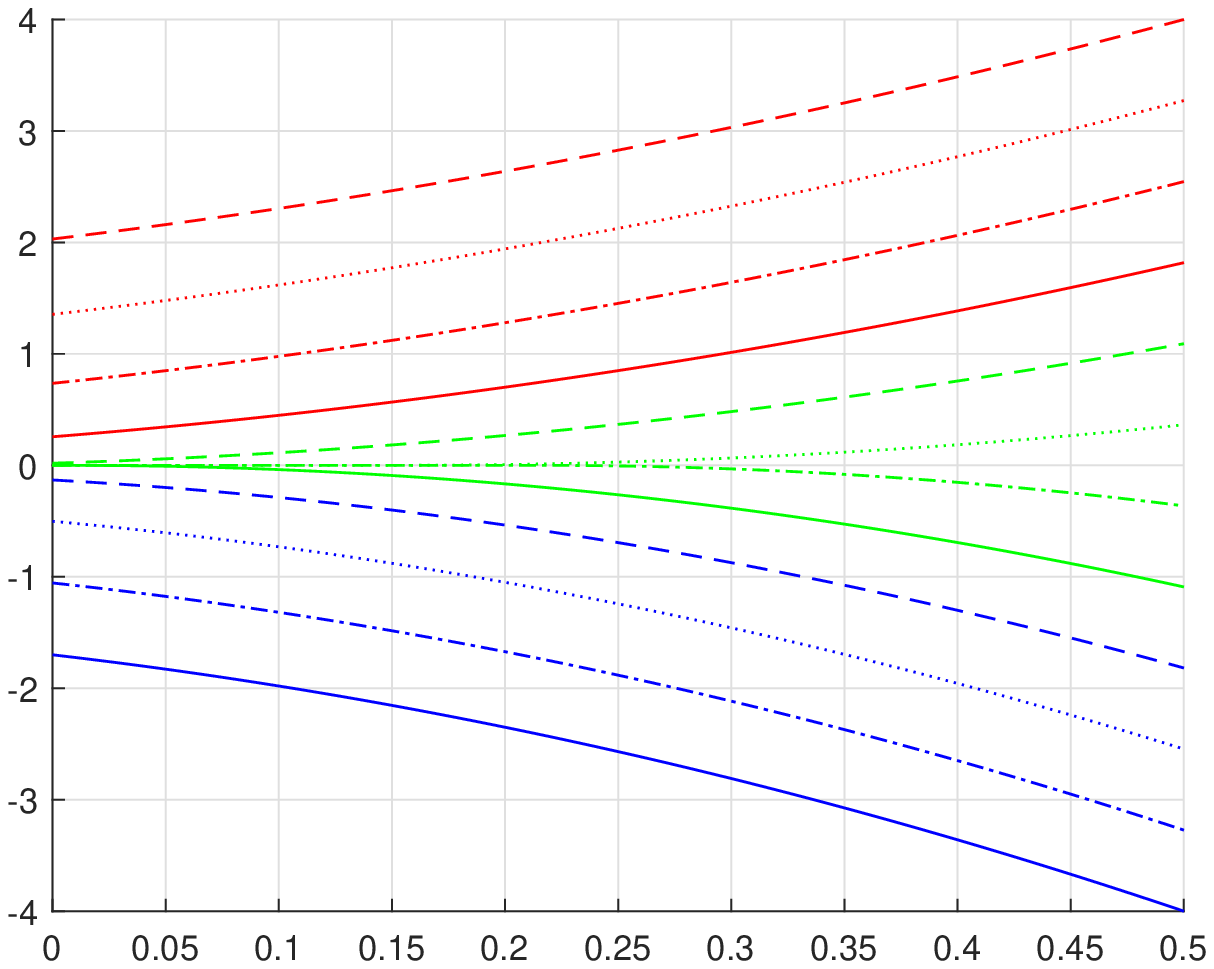}
        \caption{Second component, $t=0.5$}
    \end{subfigure}

    \caption{Evaluation of the optimal trajectory $\gmnum(s;(x,-x,0,\dots,0), t)$ of the optimal control problem (\ref{eqt: result_optctrl1_hd}) with $\ba = ( 4, 6, 5, \dots, 5)\in\R^{10}$, $\bb = ( 3, 9, 6, \dots, 6)\in\R^{10}$, and initial cost $\initcond(\bx) = \sqrt{\langle \bx, M\bx\rangle}$, where $M$ is a diagonal matrix with diagonal elements $(m_{ii}: i = 1, \dots, 10) = (1, 8, 3, 5,1,1,\dots,1) $ versus $s\in[0,t]$ for different terminal positions $(x,-x,0,\dots,0)$ ($x\in [-4,4]$) and different time horizons $t$. The different colors and line markers simply differentiate between the different trajectories. Figures (a)-(c) depict the first component of the trajectory versus $s\in[0,t]$ with different time horizons $t$, while (d)-(f) depict the second component of the trajectory versus $s\in[0,t]$ with different time horizons $t$. Plots for time horizons $t = 0.125$, $0.25$, and $0.5$ are depicted in (a)/(d), (b)/(e), and (c)/(f), respectively.}
    \label{fig: HJ1_HD_matrixnorm_trajectory}
\end{figure}

Figure \ref{fig: HJ1_HD_matrixnorm_trajectory} depicts one-dimensional slices of the optimal trajectory $\gmnum(s;\bx,t)$ of the corresponding optimal control problem (\ref{eqt: result_optctrl1_hd}) using different terminal positions $(x,-x,0,\dots,0)\in\R^{10}$ and different time horizons $t$.
In each subfigure, the time horizon $t$ is fixed, and the different trajectories correspond to different terminal positions.
We observe that the one-dimensional slices are piecewise-quadratic and continuous in $s$, which is consistent with our formulas for $\opttraj$ as defined in~\eqref{eqt: optctrl_defx_1},~\eqref{eqt: optctrl_defx_2},~\eqref{eqt: optctrl_defx_3},~\eqref{eqt: optctrl_defx_4}, and~\eqref{eqt: optctrl_defx_5}. 

In Table~\ref{tab:timing_matrixnorm}, we show the running time of this example for different dimensions $n$. We use the same method to compute the running time as in Section~\ref{subsec:numerical_quad}. From Table~\ref{tab:timing_matrixnorm}, we see that it takes, on average, less than $3\times 10^{-5}$ seconds to compute the solution at one point in a 16-dimensional problem, which demonstrates the efficiency of our proposed algorithm even in high dimensions.

\begin{table}[htbp]
    \centering
    \begin{tabular}{c|c|c|c|c}
        \hline
        \textbf{$\mathbf{n}$} & 4 & 8 & 12 & 16 \\
        \hline
        \textbf{running time (s)} & 6.9100e-06 & 9.7660e-06 & 1.5178e-05 & 2.1040e-05 \\
        \hline
    \end{tabular}
    \hfill 
    
    \caption{Time results in seconds for the average time per call over $102,400$ trials for evaluating the solution of the HJ PDE (\ref{eqt: result_HJ1_hd}) with initial cost $\initcond(\bx) = \sqrt{\langle \bx, M\bx\rangle}$, where $M$ is a diagonal matrix with diagonal elements $(m_{ii}: i = 1, \dots, n) = (1, 8, 3, 5, 1, 1,\dots, 1)$, for various dimensions $n$.}
    \label{tab:timing_matrixnorm}
\end{table}

\bigbreak
In the second example, we consider the nonsmooth convex initial cost
\begin{equation}\label{eqt: hd_initial_data_l12_shifted}
    \initcond(\bx) = \frac{1}{2} \|\bx - \boldsymbol{1}\|_1^2 \quad \forall \bx\in\Rn,
\end{equation}
where $\|\cdot\|_1$ denotes the $\ell^1$-norm in $\Rn$ and $\boldsymbol{1}$ is a vector in $\Rn$ whose components are all ones.
For this example, we update $\bv^{k+1}$ using \eqref{eqt: ADMM_vupdate_moreau}, where the proximal point of the mapping $\bx\mapsto \frac{1}{\lambda} \initcond(\lambda \bx)$ is computed efficiently using the method described in~\cite[Section~4.4]{Darbon2016Algorithms}. 
We set the parameters $\ba$ and $\bb$ to the values defined in~\eqref{eqt:numerical_def_ab}. 

\begin{figure}[htbp]
    \centering
    \begin{subfigure}{0.45\textwidth}
        \centering \includegraphics[width=\textwidth]{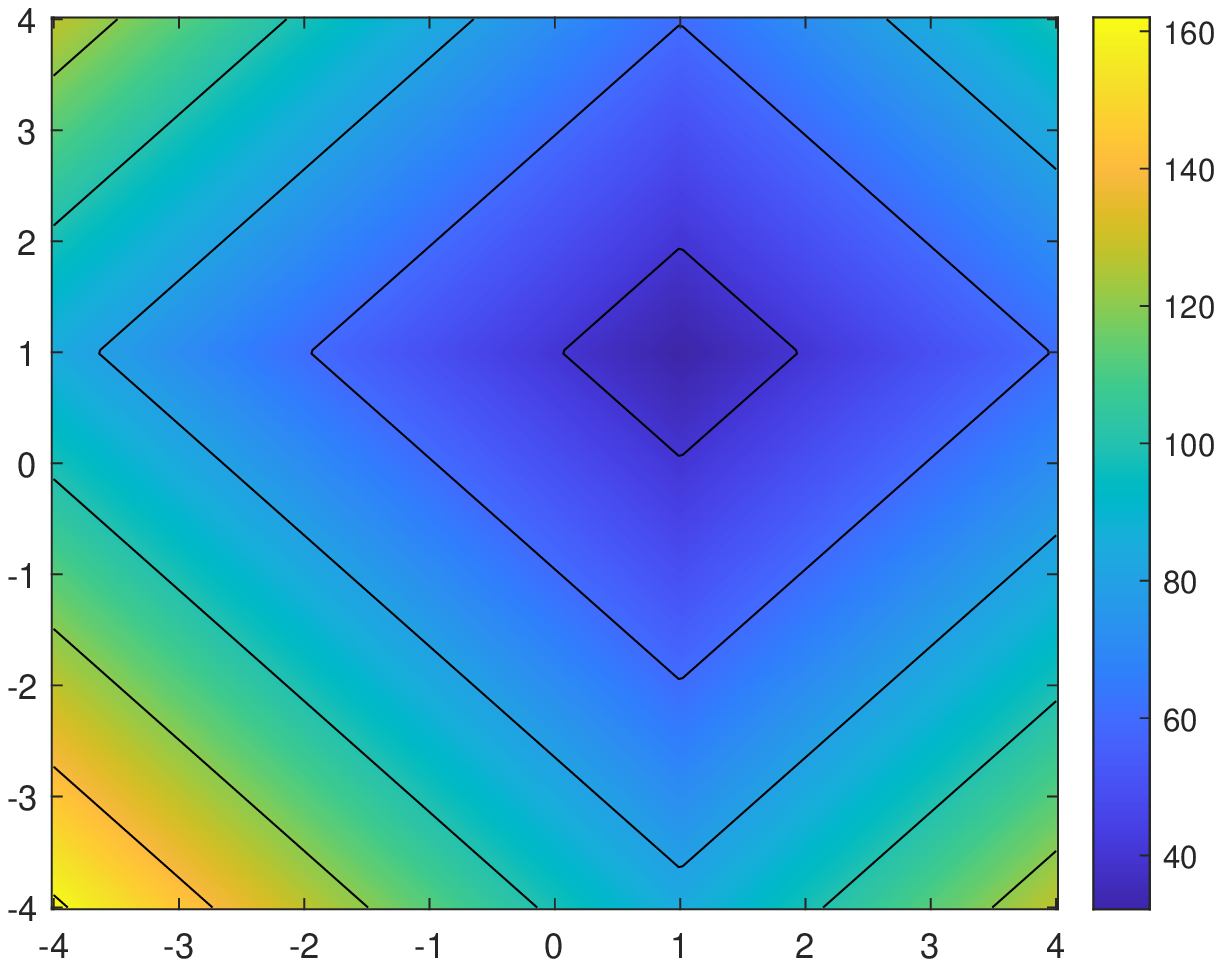}
        \caption{$t=0$}
    \end{subfigure}
    \hfill
    \begin{subfigure}{0.45\textwidth}
        \centering \includegraphics[width=\textwidth]{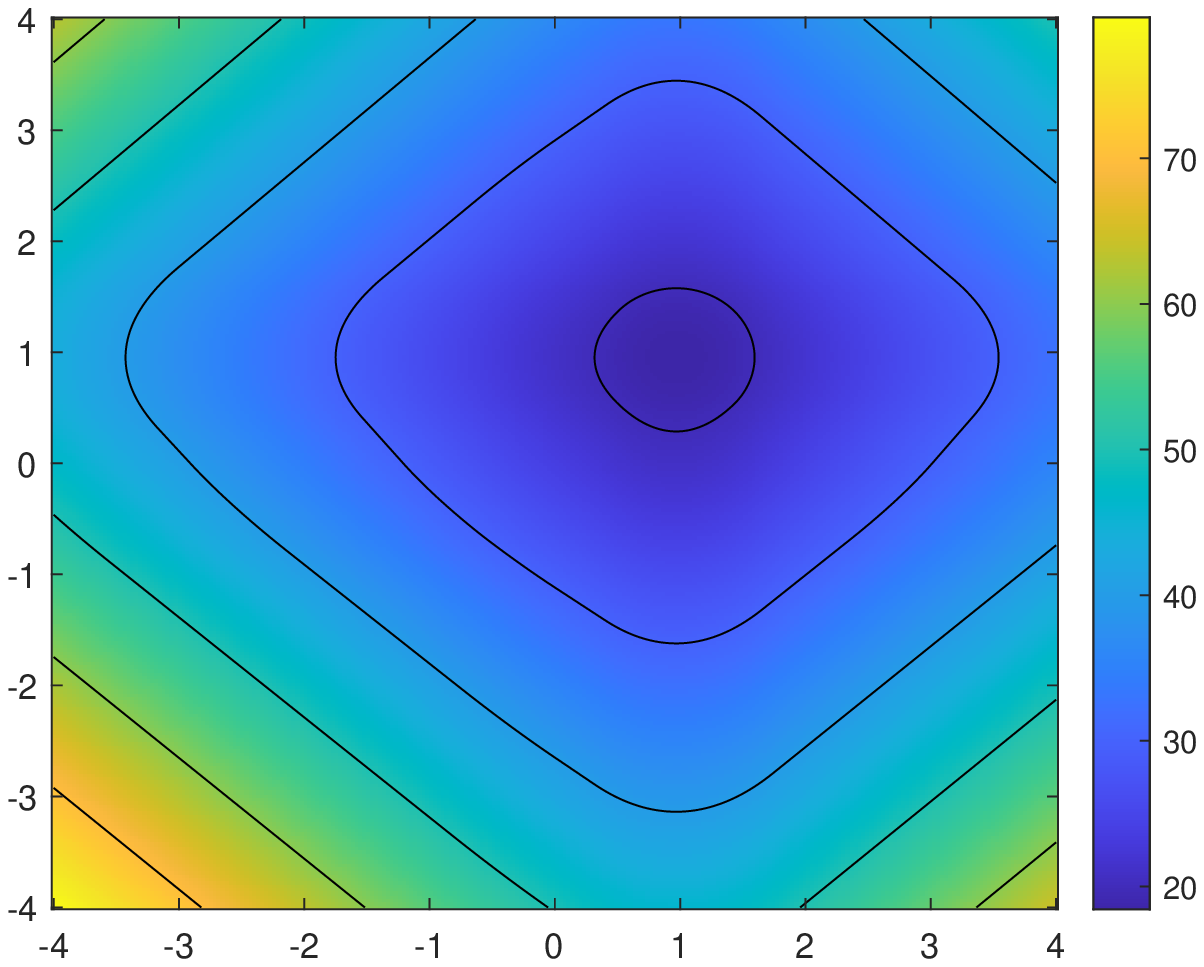}
        \caption{$t=0.125$}
    \end{subfigure}
    
    \begin{subfigure}{0.45\textwidth}
        \centering \includegraphics[width=\textwidth]{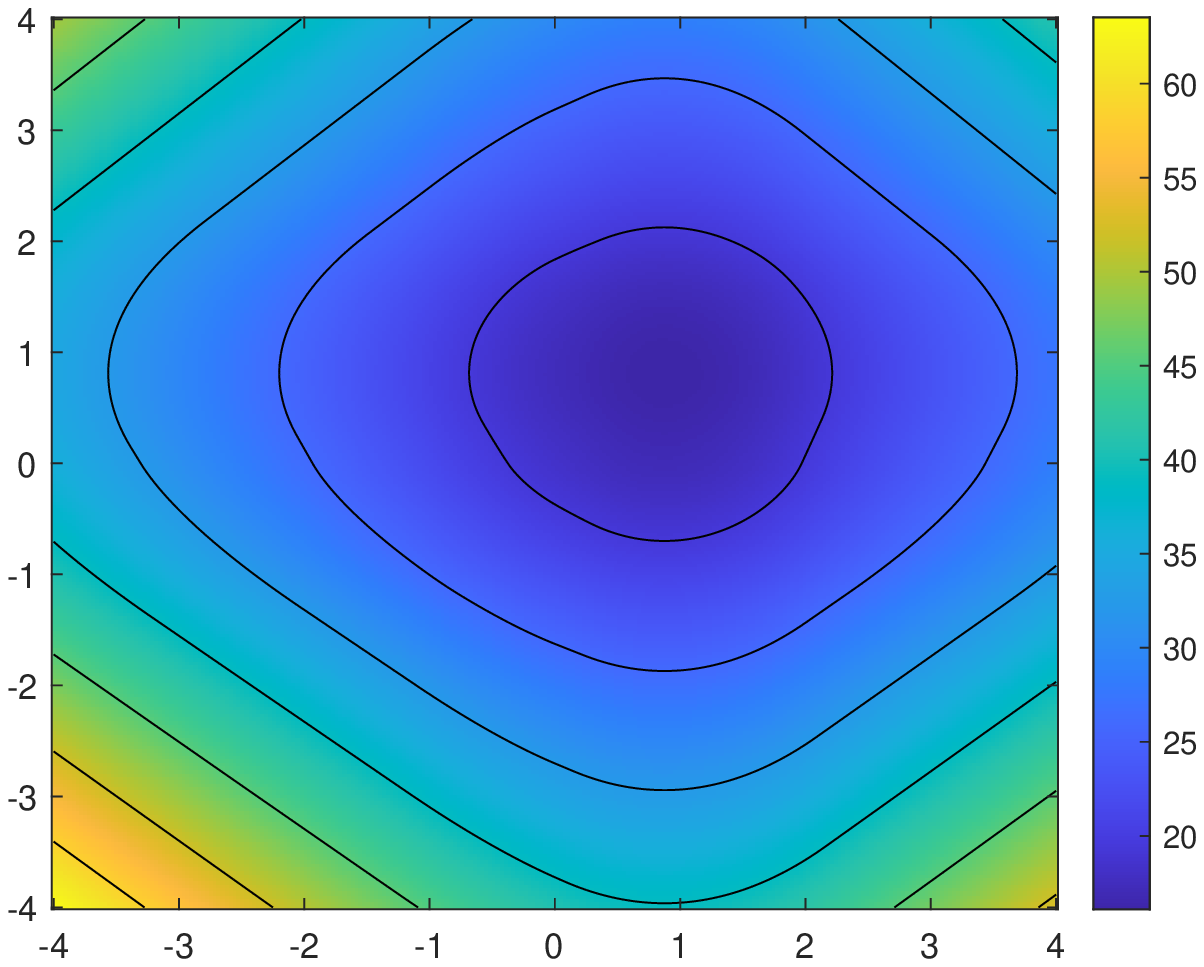}
        \caption{$t=0.25$}
    \end{subfigure}
    \hfill
    \begin{subfigure}{0.45\textwidth}
        \centering \includegraphics[width=\textwidth]{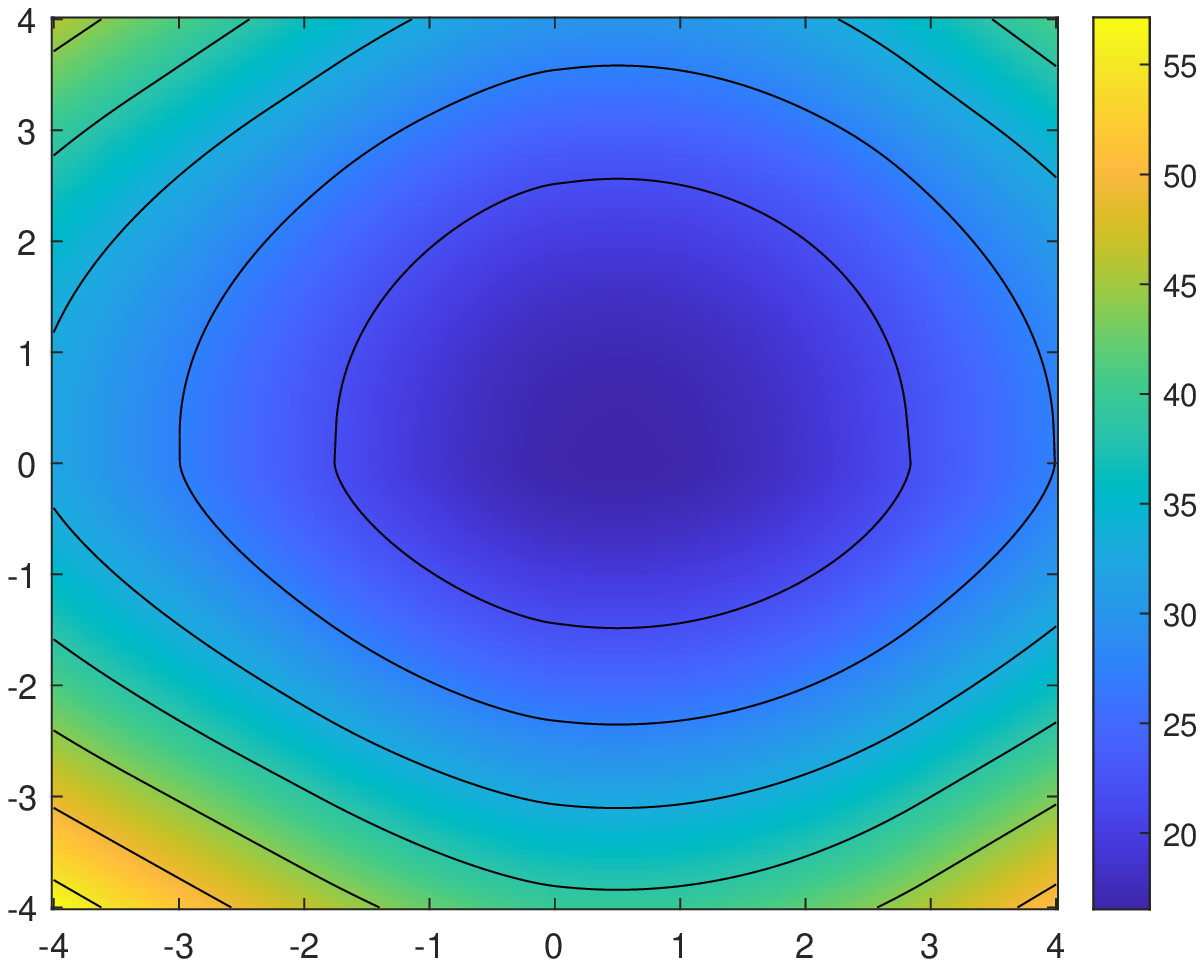}
        \caption{$t=0.5$}
    \end{subfigure}
    
    \caption{Evaluation of the solution $\Snum(\bx,t)$ of the high-dimensional HJ PDE (\ref{eqt: result_HJ1_hd}) with $\ba = ( 4, 6, 5, \dots, 5)$, $\bb=( 3, 9, 6, \dots, 6)$, and initial condition $\initcond(\bx) = \frac{1}{2}\|\bx - \boldsymbol{1}\|_1^2$ for $\bx=(x_1,x_2,0,\dots,0)\in\R^{10}$ and different times $t$. Plots for $t =0$, $0.125$, $0.25$, and $0.5$ are depicted in (a)-(d), respectively. Level lines are superimposed on the plots.}
    \label{fig: HJ1_HD_l12_shifted}
\end{figure}

\begin{figure}[htbp]
    \centering
    \begin{subfigure}{0.32\textwidth}
        \centering \includegraphics[width=\textwidth]{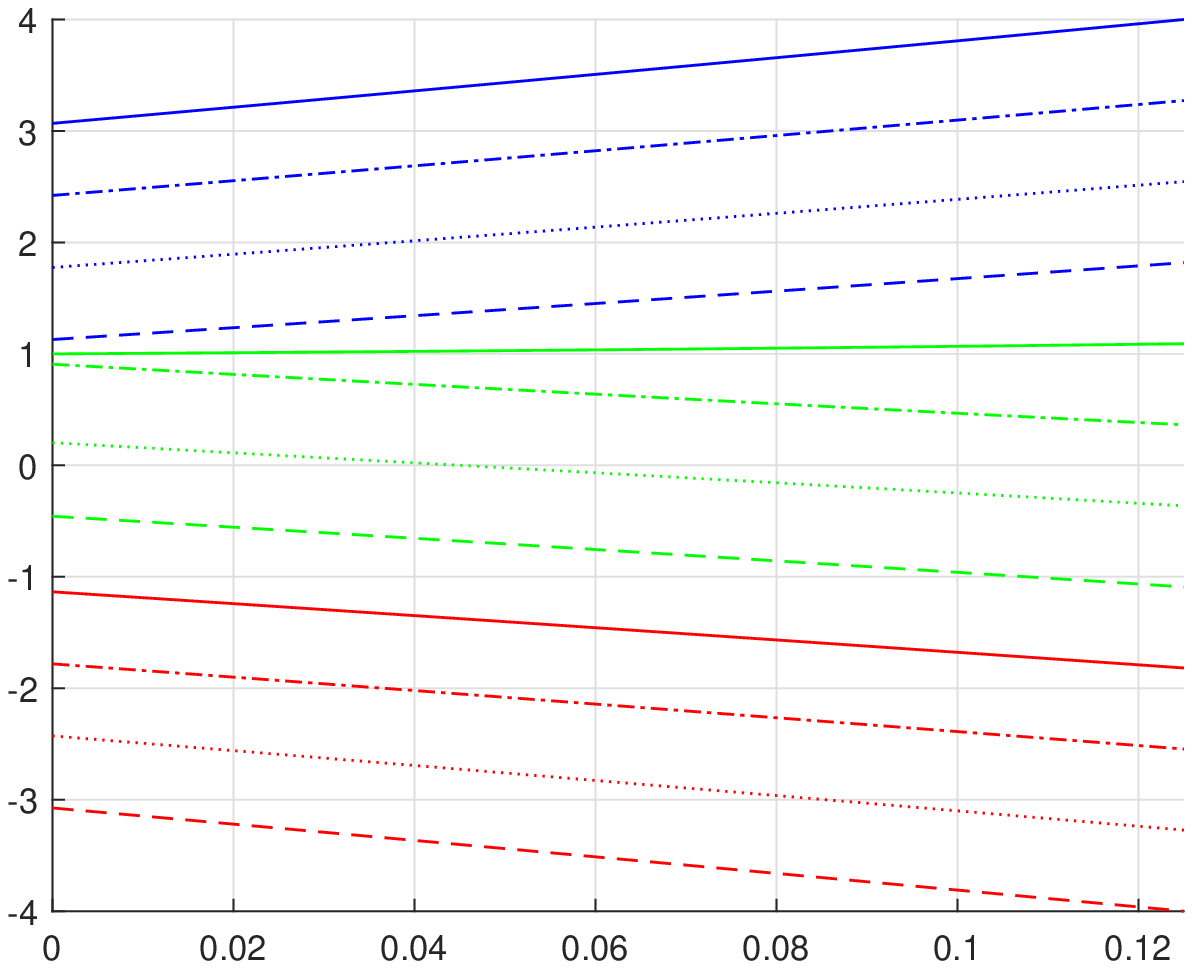}
        \caption{First component, $t=0.125$}
    \end{subfigure}
    \hfill
    \begin{subfigure}{0.32\textwidth}
        \centering \includegraphics[width=\textwidth]{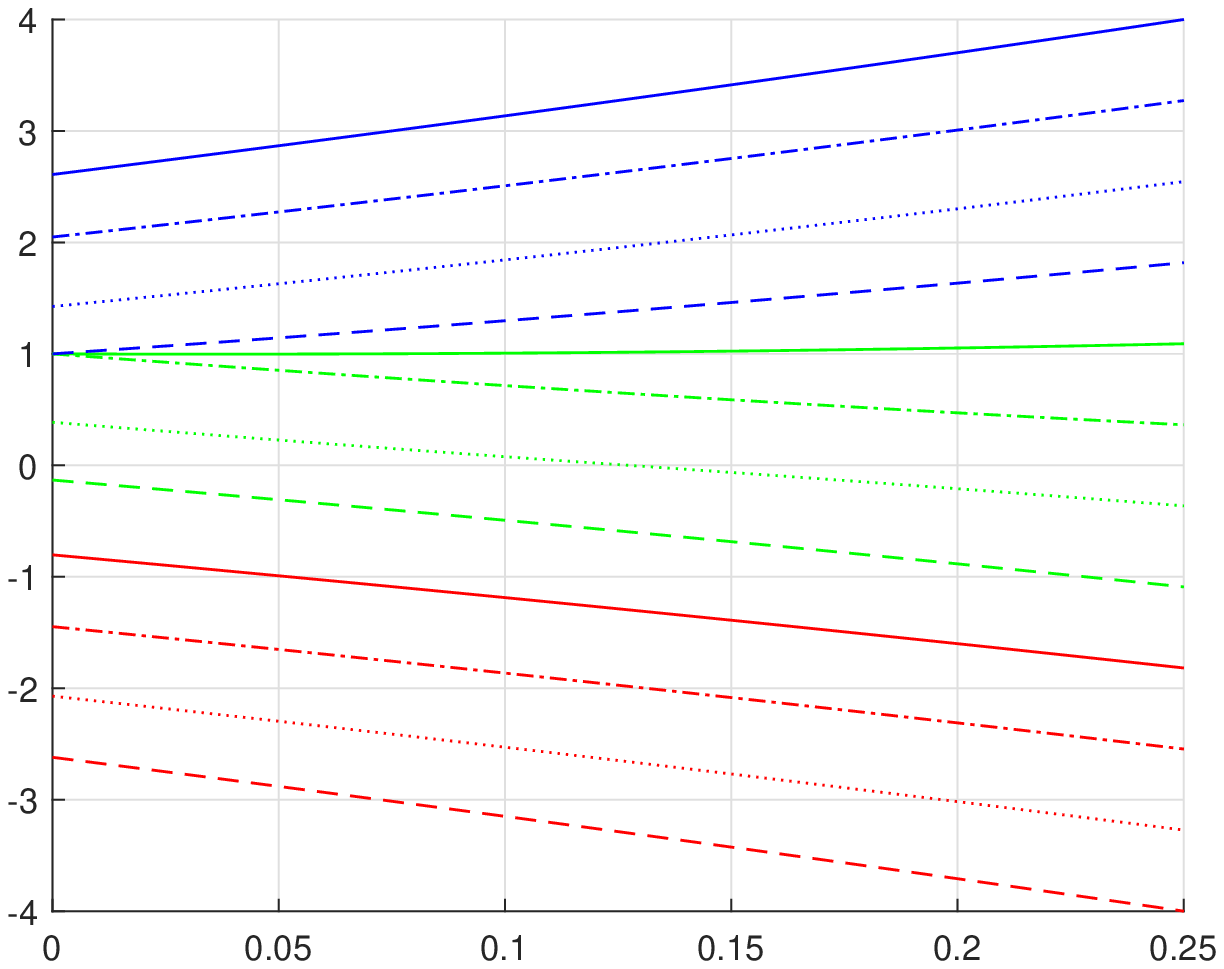}
        \caption{First component, $t=0.25$}
    \end{subfigure}
    \hfill
    \begin{subfigure}{0.32\textwidth}
        \centering \includegraphics[width=\textwidth]{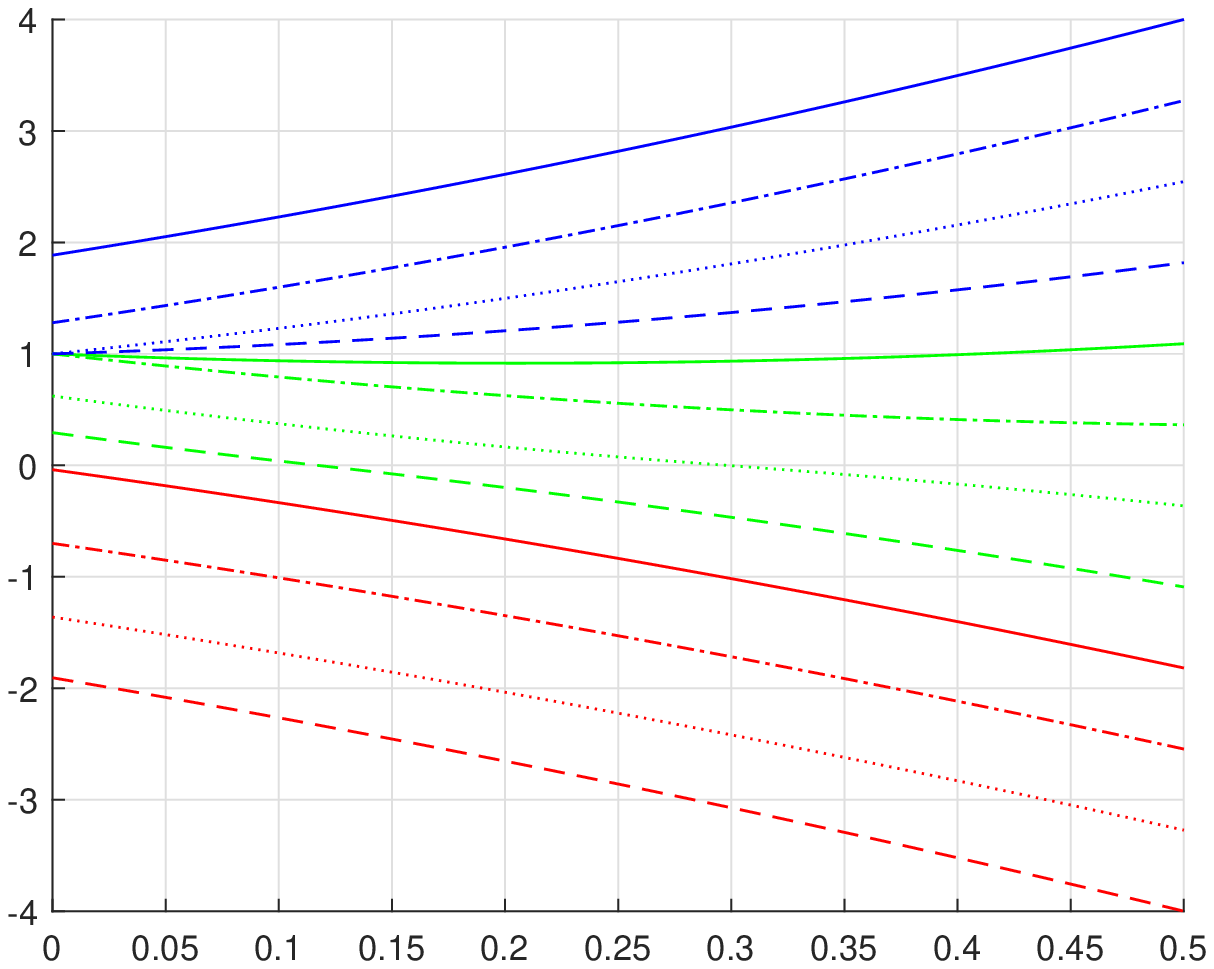}
        \caption{First component, $t=0.5$}
    \end{subfigure}
    
    \begin{subfigure}{0.32\textwidth}
        \centering \includegraphics[width=\textwidth]{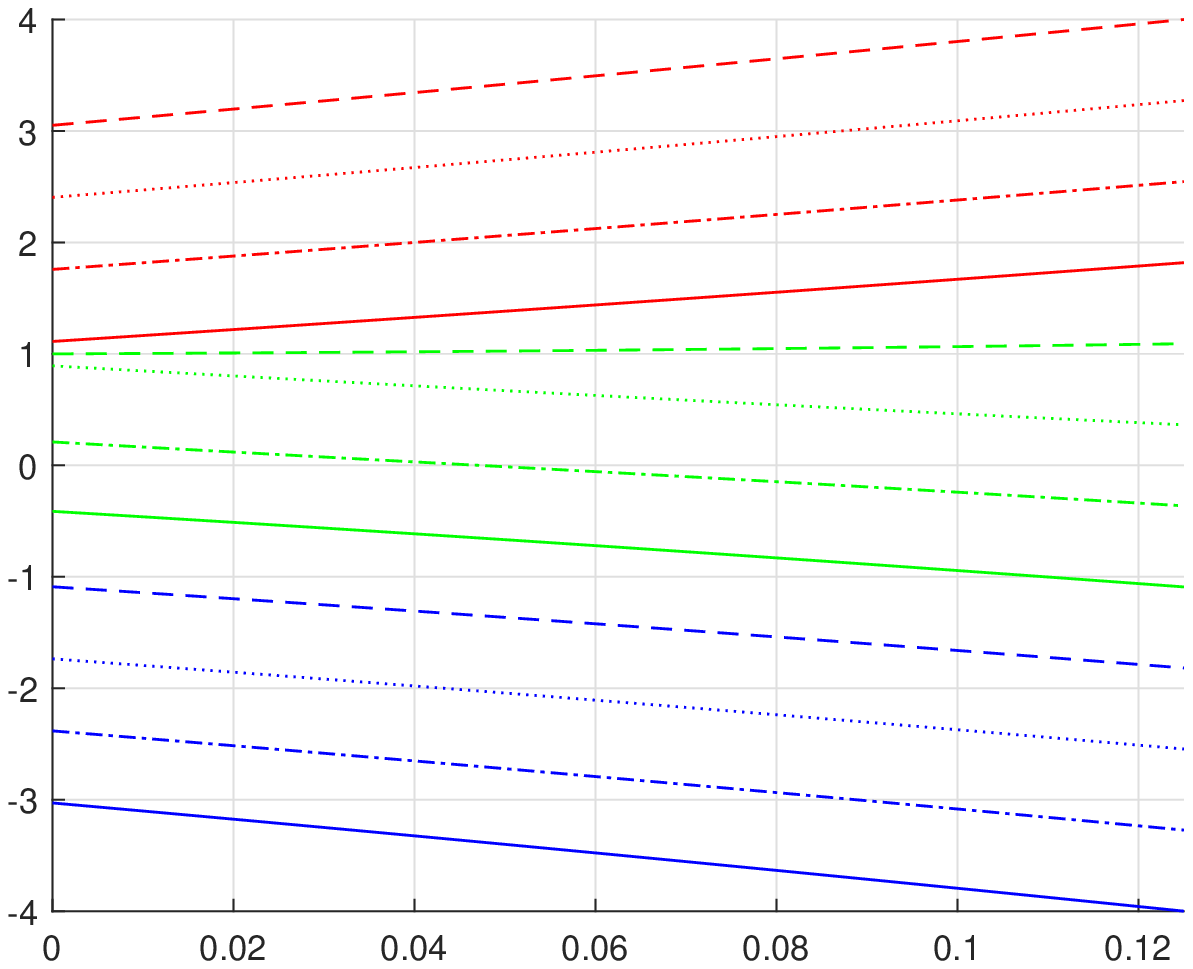}
        \caption{Second component, $t=0.125$}
    \end{subfigure}
    \hfill
    \begin{subfigure}{0.32\textwidth}
        \centering \includegraphics[width=\textwidth]{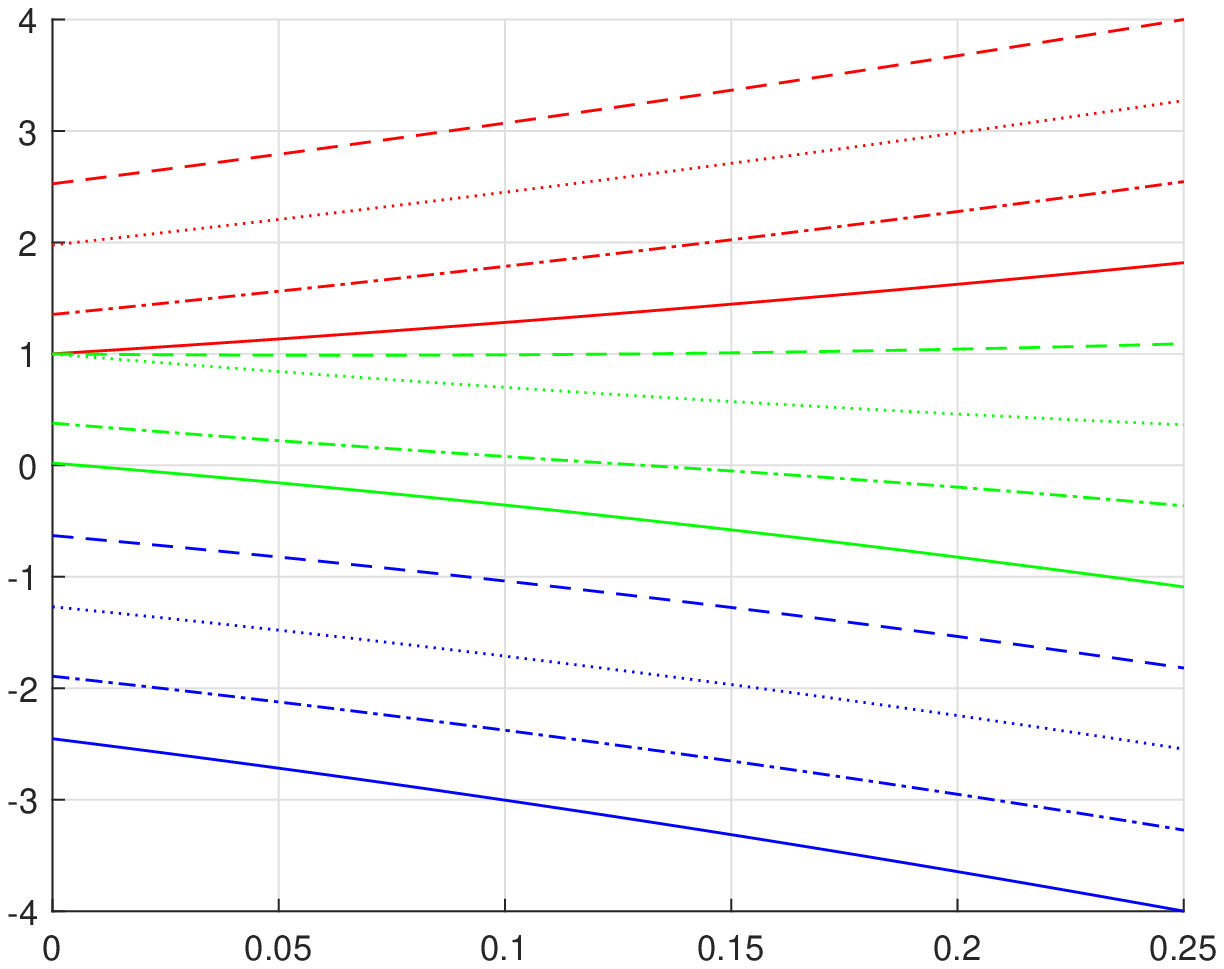}
        \caption{Second component, $t=0.25$}
    \end{subfigure}
    \hfill
    \begin{subfigure}{0.32\textwidth}
        \centering \includegraphics[width=\textwidth]{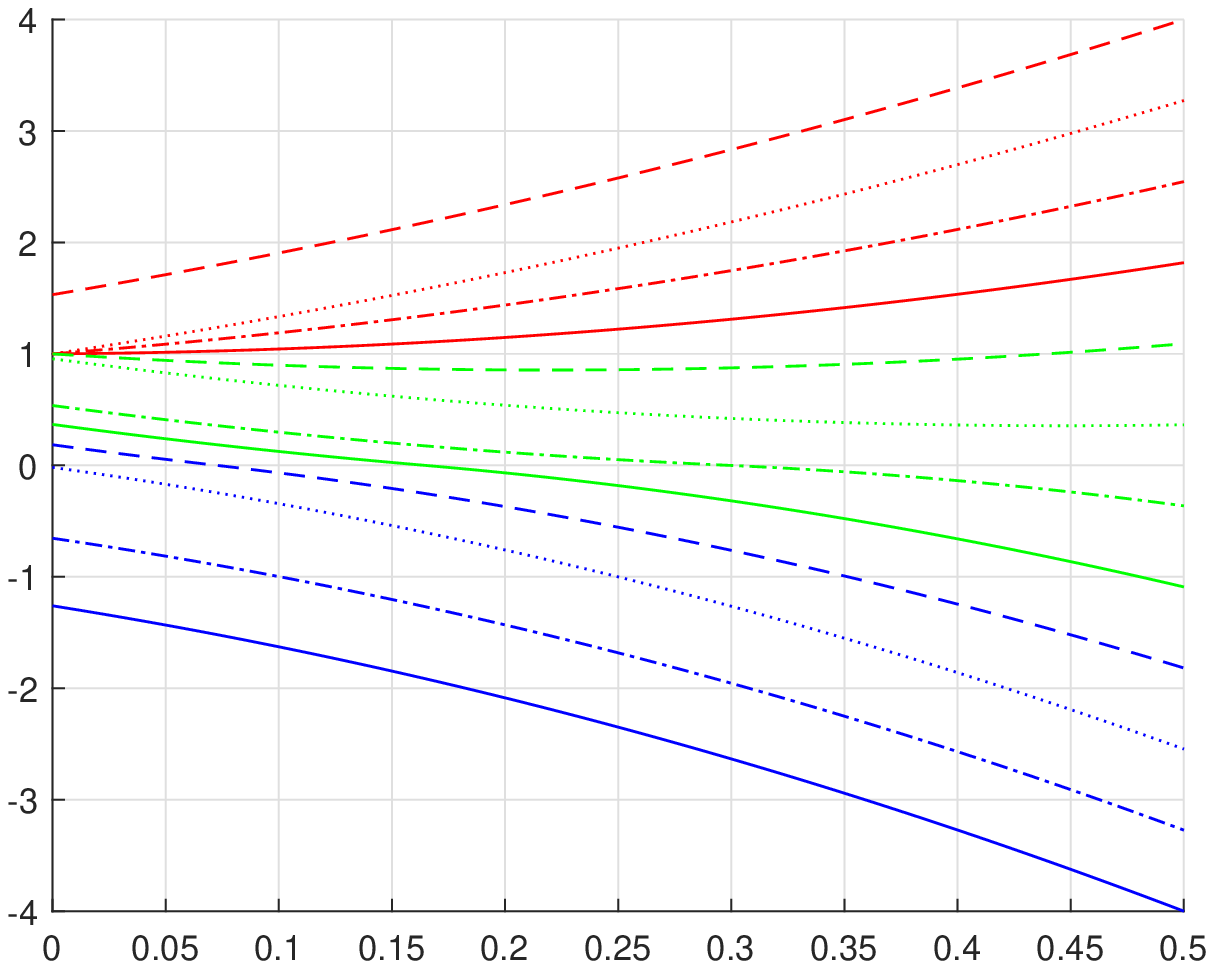}
        \caption{Second component, $t=0.5$}
    \end{subfigure}

    \caption{Evaluation of the optimal trajectory $\gmnum(s;(x,-x,0,\dots,0), t)$ of the optimal control problem (\ref{eqt: result_optctrl1_hd}) with $\ba = ( 4, 6, 5, \dots, 5)$, $\bb = ( 3, 9, 6, \dots, 6)$, and initial cost $\initcond(\bx) = \frac{1}{2}\|\bx - \mathbf{1}\|_1^2$ versus $s\in[0,t]$ for different terminal positions $(x,-x,0,\dots,0)$ ($x\in [-4,4]$) and different time horizons $t$. The different colors and line markers simply differentiate between the different trajectories. Figures (a)-(c) depict the first component of the trajectory versus $s\in[0,t]$ with different time horizons $t$, while (d)-(f) depict the second component of the trajectory versus $s\in[0,t]$ with different time horizons $t$. Plots for time horizons $t = 0.125$, $0.25$, and $0.5$ are depicted in (a)/(d), (b)/(e), and (c)/(f), respectively.}
    \label{fig: HJ1_HD_l12_shifted_trajectory}
\end{figure}

Figure \ref{fig: HJ1_HD_l12_shifted} depicts two-dimensional slices of the numerical solution $\Snum(\bx,t)$, as computed using Algorithm~\ref{alg:admm_ver1}, of the $10$-dimensional HJ PDE \eqref{eqt: result_HJ1_hd} 
at different positions $\bx=(x_1, x_2, 0, \dots, 0)$ and at different times $t$. 
Figure \ref{fig: HJ1_HD_l12_shifted_trajectory} depicts one-dimensional slices of the optimal trajectory $\gmnum(s;\bx,t)$ of the corresponding optimal control problem (\ref{eqt: result_optctrl1_hd}), using different terminal positions $\bx=(x,-x,0,\dots,0)\in\R^{10}$ and different time horizons $t$. We observe that the one-dimensional slices are piecewise quadratic and continuous in $s$, which is consistent with our formulas for $\opttraj$, given by~\eqref{eqt: optctrl_defx_1},~\eqref{eqt: optctrl_defx_2},~\eqref{eqt: optctrl_defx_3},~\eqref{eqt: optctrl_defx_4}, and~\eqref{eqt: optctrl_defx_5}.

In Table \ref{tab:timing_L1sq}, we present the running time of this example for different dimensions $n$. From Table~\ref{tab:timing_L1sq}, we see that it takes less than $4\times 10^{-4}$ seconds on average to compute the solution at one point in a 16-dimensional problem, which demonstrates the efficiency of our proposed algorithm even in high dimensions.

\begin{table}[htbp]
    \centering
    \begin{tabular}{c|c|c|c|c}
        \hline
        \textbf{$\mathbf{n}$} & 4 & 8 & 12 & 16 \\
        \hline
        \textbf{running time (s)} & 2.1192e-05 & 9.4819e-05 & 2.0531e-04 & 3.2751e-04 \\
        \hline
    \end{tabular}
    \hfill 
    
    \caption{Time results in seconds for the average time per call over $102,400$ trials for evaluating the solution of the HJ PDE (\ref{eqt: result_HJ1_hd}) with initial cost $\initcond(\bx) = \frac{1}{2} \|\bx - \mathbf{1}\|_1^2$ for various dimensions $n$.}
    \label{tab:timing_L1sq}
\end{table}
 
\subsection{A class of non-convex initial costs}\label{sec: ADMM_nonconvex}

\begin{algorithm}[htbp]
\SetAlgoLined
\SetKwInOut{Input}{Inputs}
\SetKwInOut{Output}{Outputs}
\Input{Parameters $\ba,\bb\in\Rn$, terminal position $\bx\in\Rn$, time horizon $t>0$, running time $s>0$ of the trajectory, and the convex functions $\initcond_1,\dots,\initcond_m$ in~\eqref{eqt: J_minplus}.
}
\Output{An optimal trajectory $\gmnum(s;\bx,t)$ in the optimal control problem~\eqref{eqt: result_optctrl1_hd} and the solution value $\Snum(\bx,t)$ to the corresponding HJ PDE~\eqref{eqt: result_HJ1_hd} with nonconvex initial cost $\initcond$ of the form~\eqref{eqt: J_minplus}.}
 \For{$j = 1,2,\dots,m$}{
 Numerically solve the $j$-th subproblem~\eqref{eqt:numerical_minplus_subproblem} using any appropriate method (e.g., the solver in Section~\ref{subsec:numerical_quad} or Algorithm~\ref{alg:admm_ver1} in Section~\ref{sec:numerical_convex}), and get
 the optimal trajectory $\gmnum_j(s;\bx,t)$ of optimal control problem~\eqref{eqt: result_optctrl1_hd} and the solution $\Snum_j(\bx,t)$ to the corresponding HJ PDE~\eqref{eqt: result_HJ1_hd} with initial data $\initcond_j$\;
 }
 Compute the index $r$ by
 \begin{equation}\label{eqt:minplus_defr}
    r \in \argmin_{j\in\{1, \dots, m\}} \Snum_j (\bx,t).
\end{equation}
\\
 Output an optimal trajectory $\gmnum(s;\bx,t)$ and the solution value $\Snum(\bx,t)$ using
 \begin{equation}\label{eqt:alg_minplus_output}
 \begin{split}
    \gmnum(s;\bx,t) & = \gmnum_r(s;\bx,t), \quad\quad    \Snum(\bx,t)  = \Snum_r(\bx,t) = \min_{j\in\{1, \dots, m\}} \Snum_j (\bx,t).
\end{split}
\end{equation}
 \caption{An optimization algorithm for solving the optimal control problem~\eqref{eqt: result_optctrl1_hd} and the corresponding HJ PDE~\eqref{eqt: result_HJ1_hd} with nonconvex initial cost $\initcond$ of the form~\eqref{eqt: J_minplus}. \label{alg:admm_minplus}}
\end{algorithm}

In this section, we provide an algorithm based on the min-plus technique to solve the high-dimensional problems~\eqref{eqt: result_optctrl1_hd} and~\eqref{eqt: result_HJ1_hd} with certain nonconvex initial data of the form~\eqref{eqt: J_minplus}. 
The algorithm is summarized in Algorithm \ref{alg:admm_minplus}.

Recall that the solution $\valuefn$ is given by \eqref{eqt:Hopf_minplus} and the optimal trajectory $\opttrajhd$ is given by \eqref{eqt:opttraj_minplus}. Thus, to solve~\eqref{eqt: result_optctrl1_hd} and~\eqref{eqt: result_HJ1_hd} with initial cost $\initcond$ of the form~\eqref{eqt: J_minplus}, we first divide the problem into $m$ subproblems. In the $j$-th subproblem, which corresponds to the initial cost $\initcond_j$, we must solve the following optimization problem:
\begin{equation} \label{eqt:numerical_minplus_subproblem}
    \min_{\bp\in\R^n} \left\{-\sum_{i=1}^n \valuefn(x_i, t; p_i, a_i,b_i) + \initcond_j^*(\bp)\right\}.
\end{equation}
We can apply any appropriate algorithm to solve this convex optimization problem, such as the methods in Sections~\ref{subsec:numerical_quad} and~\ref{sec:numerical_convex}. Note that the subproblems can be solved in parallel and possibly using different algorithms, the choice of which depends on the properties of $\initcond_j$.
We then compute the final solution to the overall problem using the solutions to each of the subproblems (i.e. the optimal cost $\Snum_j(\bx,t)$ and the optimal trajectory $\gmnum_j(s)$ of the $j$-th subproblem for $j = 1, \dots, m$) and~\eqref{eqt:alg_minplus_output}.
As defined in~\eqref{eqt:minplus_defr}, the index $r$ in~\eqref{eqt:alg_minplus_output} is the index of the minimal cost $\Snum_r(\bx,t)$ among all possible costs $\Snum_1(\bx,t), \dots, \Snum_m(\bx,t)$.
As noted in Remark \ref{rmk:minplus_nonunique}, the index $r$ and hence the optimal trajectory $\opttrajhd(s;\bx,t)$ may be non-unique due to the nonconvexity of the initial data.

Note that~\cite[Proposition~8]{Chen2021Lax} still holds (since the proof only relies on the min-plus technique), and hence, the convergence of Algorithm~\ref{alg:admm_minplus} is guaranteed.
In other words, as long as the algorithm for each subproblem converges, the numerical solution $\Snum$ given by Algorithm~\ref{alg:admm_minplus} converges pointwise to the analytical solution. Moreover, any cluster point of the numerical optimal trajectory $\gmnum$ yields an optimal trajectory in~\eqref{eqt: result_optctrl1_hd}. As noted previously, since the initial cost $\initcond$ is nonconvex, the optimal trajectory of the optimal control problem~\eqref{eqt: result_optctrl1_hd} may be non-unique, and thus, the output trajectory $s\mapsto\gmnum(s;\bx,t)$ may have multiple cluster points. Therefore, the conclusion holds only for the cluster points of $\gmnum$, and the convergence of $\gmnum$ is not guaranteed.
Specifically, whenever the optimal trajectory is unique, the output trajectories of Algorithm~\ref{alg:admm_minplus} converge as the error in each subproblem converges to zero.

Next, we present a high-dimensional numerical example using nonconvex $\initcond$ of the form~\eqref{eqt: J_minplus}. More specifically, we consider the following nonconvex initial cost:
\begin{equation}\label{eqt: hd_initial_data_minofquad}
    \initcond(\bx) = \min_{j\in\{1,2,3\}} \initcond_j(\bx) = \min_{j\in\{1,2,3\}} \left\{\frac{1}{2}\|\bx - \by_j\|^2 + \alpha_j\right\},
\end{equation}
where $\by_1 = (-2,0,\dots,0)$, $\by_2 = (2,-2,-1,0, \dots, 0)$, and $\by_3 = (0, 2, 0, \dots, 0)$ are vectors in $\R^n$ and $\alpha_1 = -0.5$, $\alpha_2 = 0$, $\alpha_3 = -1$ are scalars in $\R$. Recall that $\|\cdot\|$ denotes the $\ell^2$-norm in the Euclidean space $\Rn$. We also set $\ba$ and $\bb$ to be the vectors defined in~\eqref{eqt:numerical_def_ab}.

\begin{figure}[htbp]
    \centering
    \begin{subfigure}{0.45\textwidth}
        \centering \includegraphics[width=\textwidth]{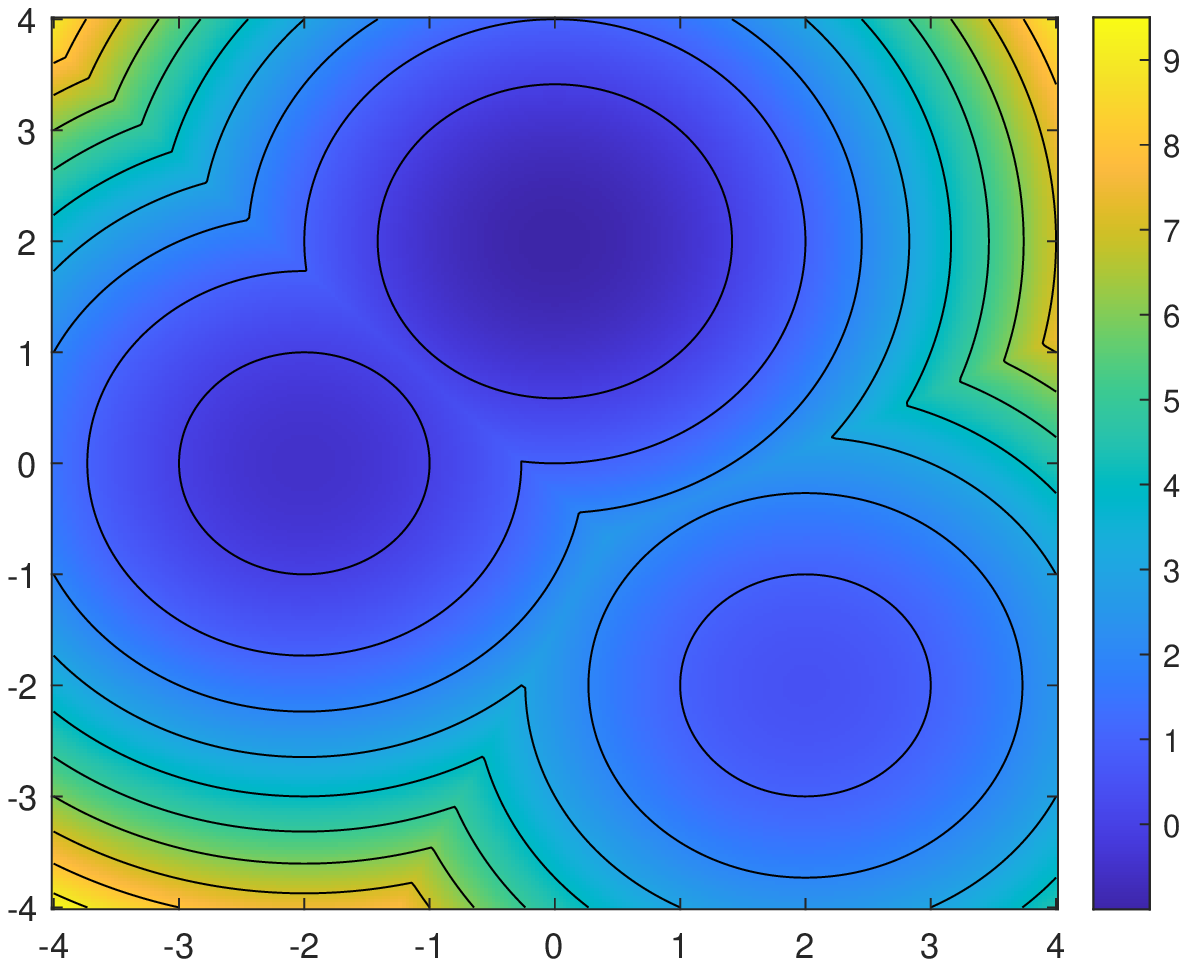}
        \caption{$t=0$}
    \end{subfigure}
    \hfill
    \begin{subfigure}{0.45\textwidth}
        \centering \includegraphics[width=\textwidth]{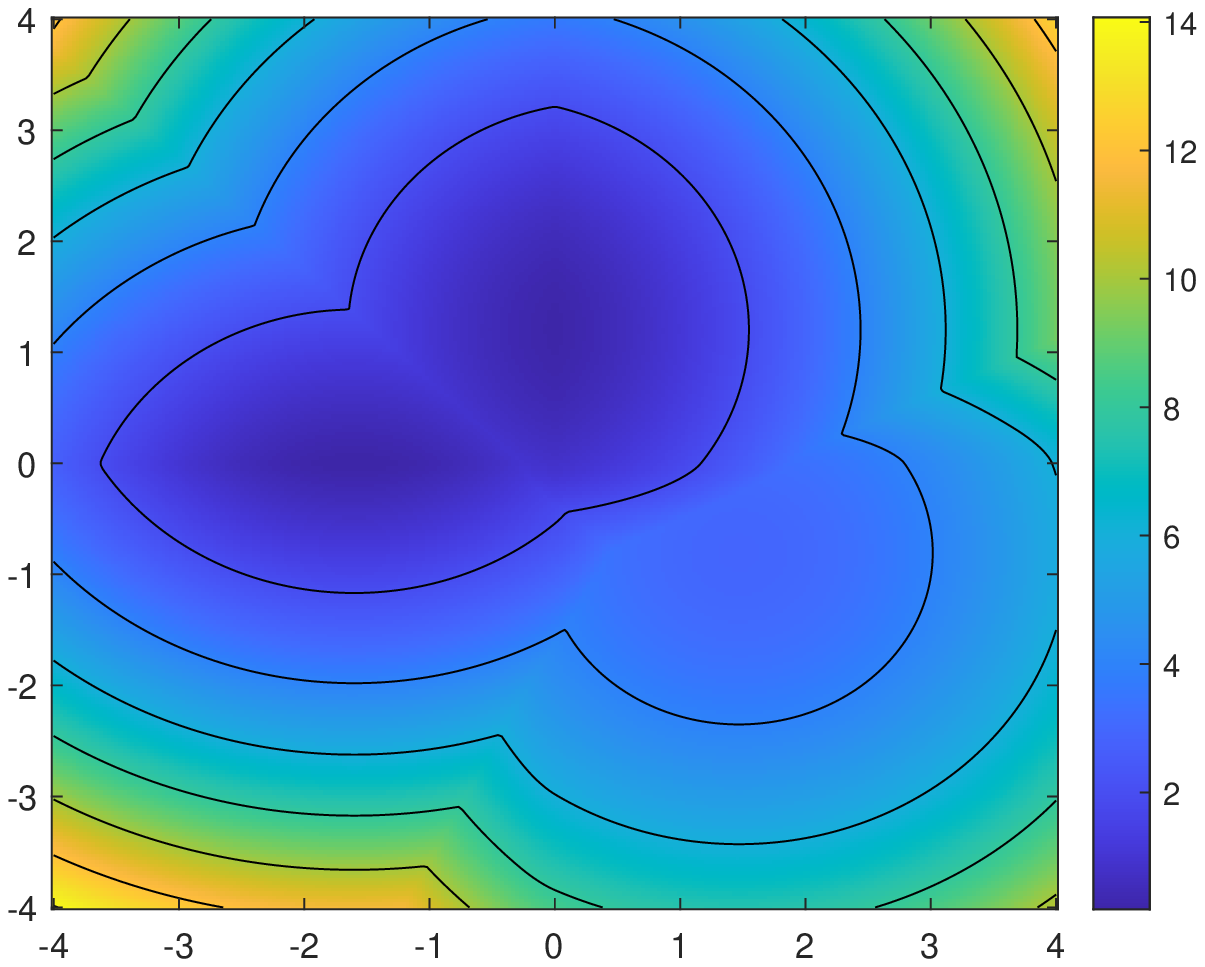}
        \caption{$t=0.125$}
    \end{subfigure}
    
    \begin{subfigure}{0.45\textwidth}
        \centering \includegraphics[width=\textwidth]{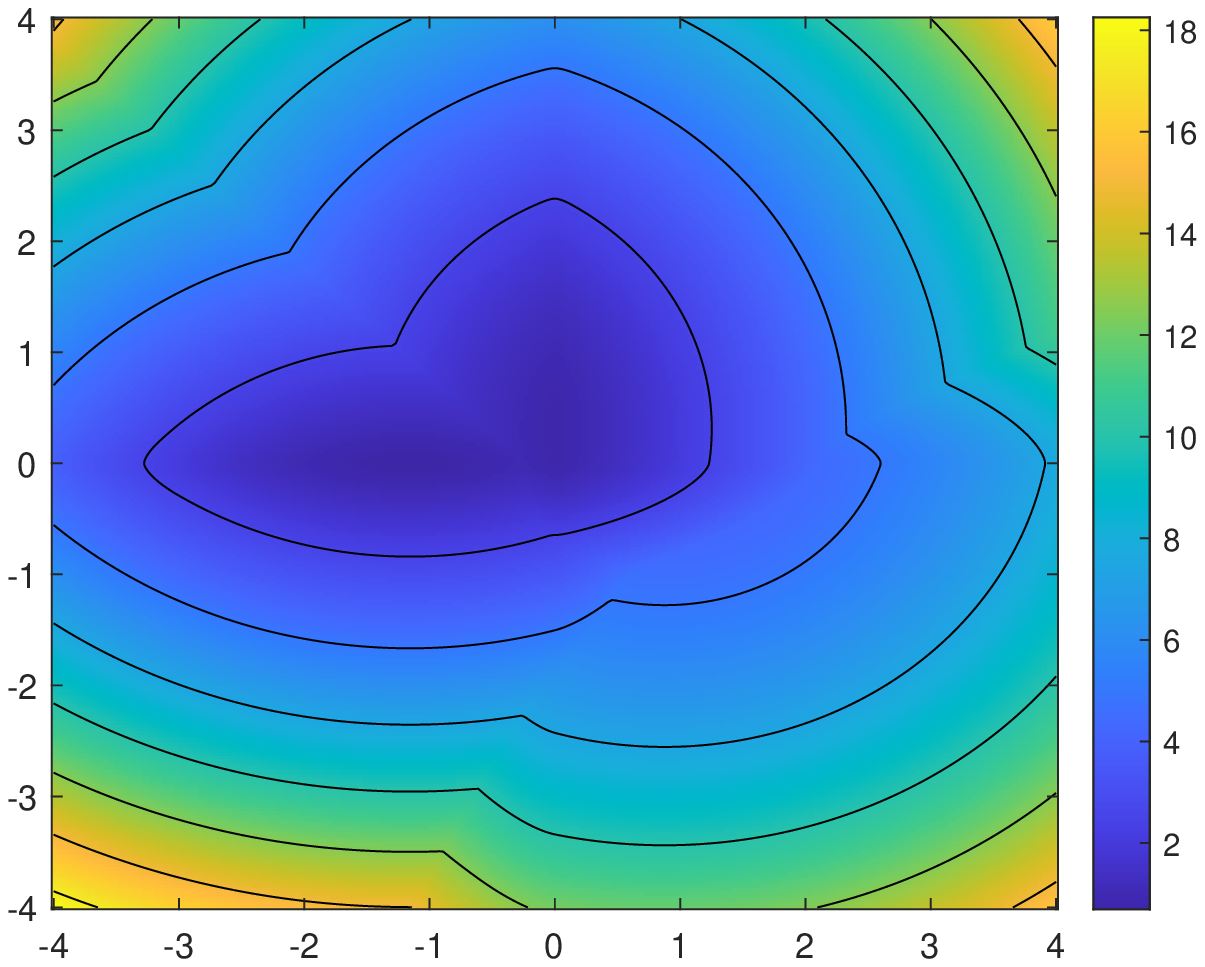}
        \caption{$t=0.25$}
    \end{subfigure}
    \hfill
    \begin{subfigure}{0.45\textwidth}
        \centering \includegraphics[width=\textwidth]{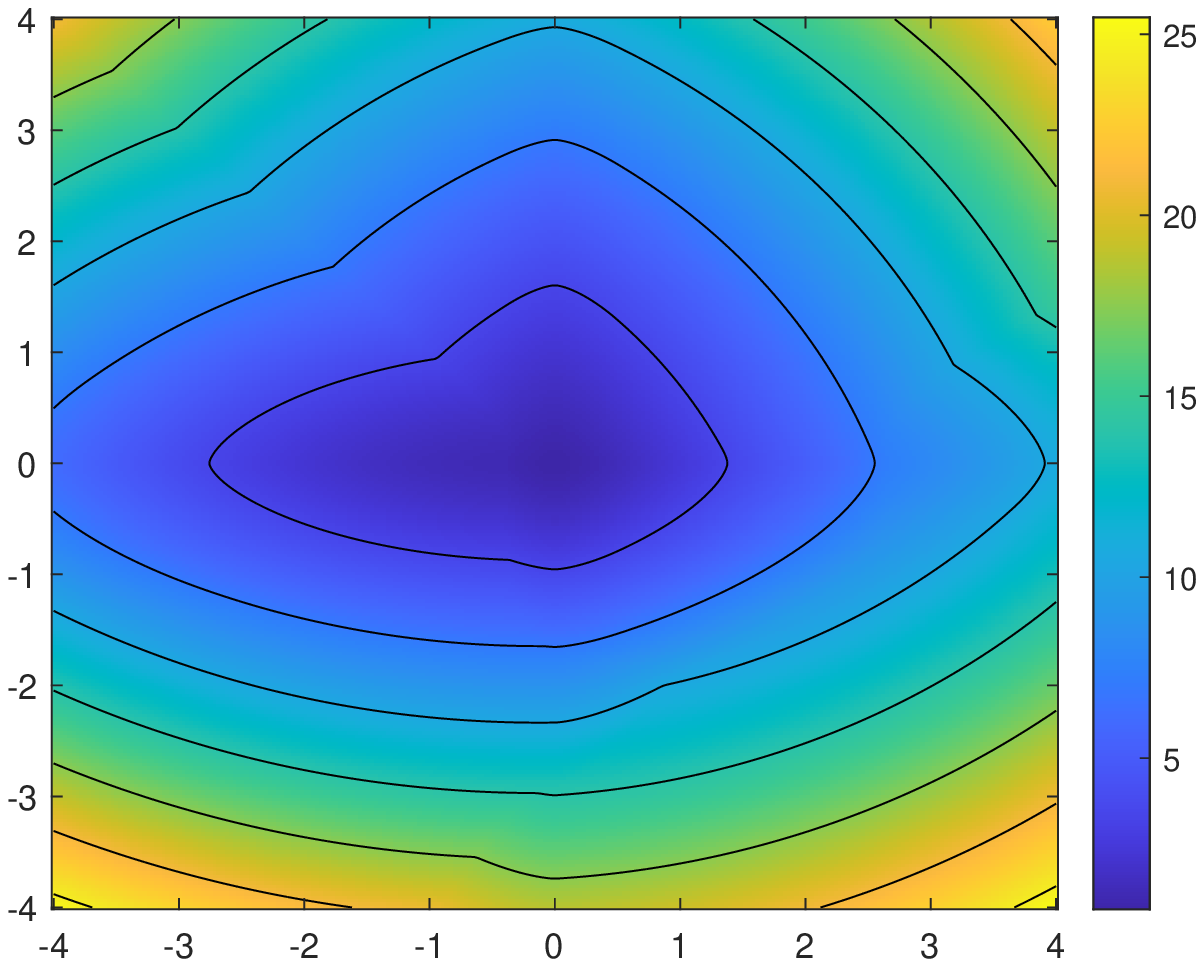}
        \caption{$t=0.5$}
    \end{subfigure}
    
    \caption{Evaluation of the solution $\Snum(\bx,t)$ of the $10$-dimensional HJ PDE (\ref{eqt: result_HJ1_hd}) with $\ba = ( 4, 6, 5, \dots, 5)$, $\bb=( 3, 9, 6, \dots, 6)$, and initial condition $\initcond(\bx) = \min_{j\in\{1,2,3\}} \initcond_j(\bx) = \min_{j\in\{1,2,3\}} \{\frac{1}{2}\|\bx - \by_j\|^2 + \alpha_j\}$, where $\by_1 = (-2,0,\dots,0)$, $\by_2 = (2,-2,-1,0, \dots, 0)$, $\by_3 = (0, 2, 0, \dots, 0)$, $\alpha_1 = -0.5$, $\alpha_2 = 0$, and $\alpha_3 = -1$, for $\bx = (x_1,x_2,0,\dots,0)$ and different times $t$. Plots for $t =0$, $0.125$, $0.25$, and $0.5$ are depicted in (a)-(d), respectively. Level lines are superimposed on the plots.}
    \label{fig:minplus_quad}
\end{figure}

Figure~\ref{fig:minplus_quad} depicts two-dimensional slices of the solution $\Snum(\bx,t)$, as computed using Algorithm \ref{alg:admm_minplus}, to the $10$-dimensional HJ PDE~\eqref{eqt: result_HJ1_hd} for different positions $\bx=(x_1, x_2, 0, \dots, 0)$ and different times $t$. In Figure \ref{fig:minplus_quad}(a), we clearly see that the initial condition $\initcond$ is not smooth at the interfaces of the quadratics $\initcond_i$, e.g., there are obvious kinks near $(x_1,x_2) = (0,0)$, $(-2,2)$, $(0,-2)$, $(2.5,0)$. We see that over time, the solution also evolves with several kinks (Figures \ref{fig:minplus_quad}(b)-(d)). These kinks provide numerical validation that the algorithm does indeed provide the non-smooth viscosity solution to the corresponding HJ PDE. 

\begin{figure}[htbp]
    \centering
    \begin{subfigure}{0.32\textwidth}
        \centering \includegraphics[width=\textwidth]{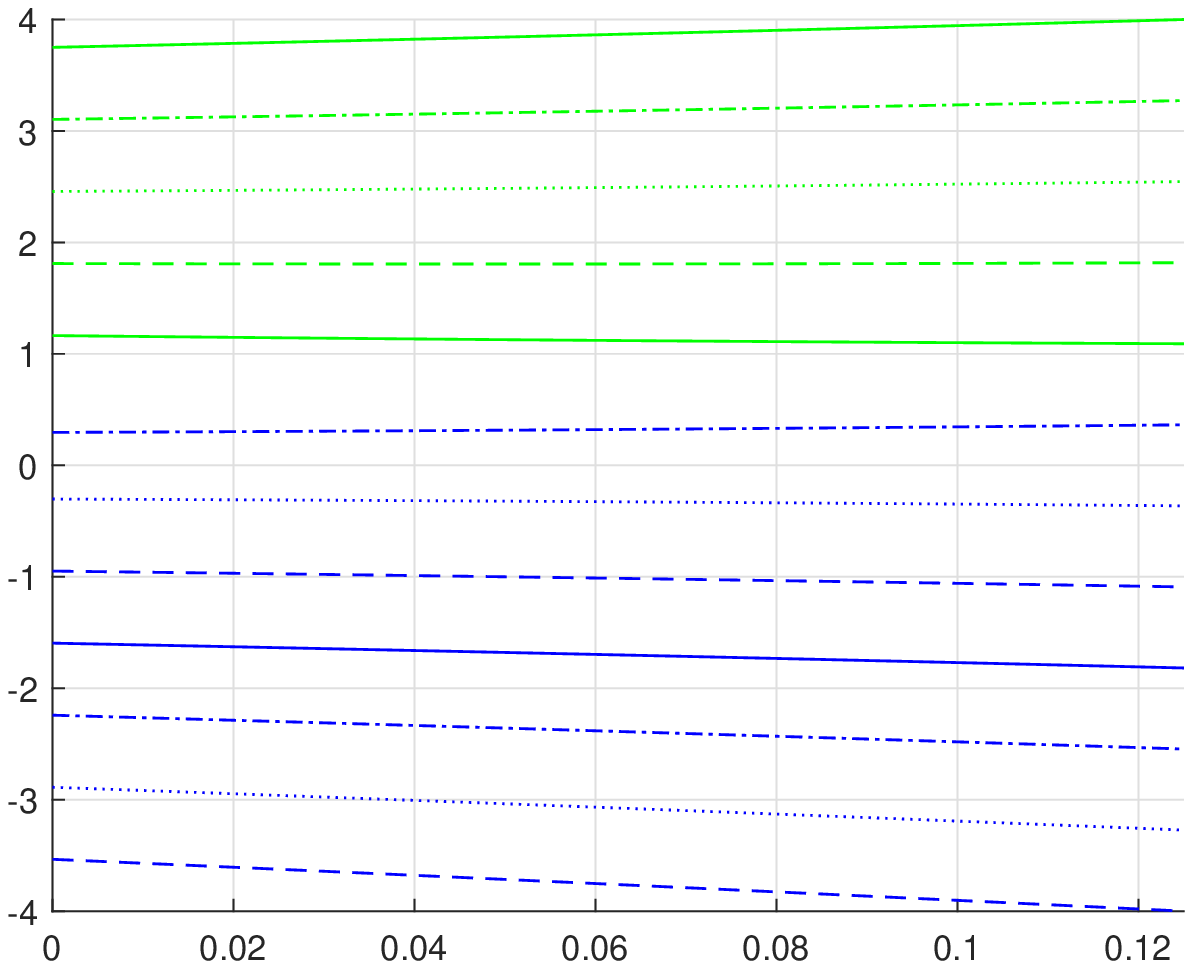}
        \caption{First component, $t=0.125$}
    \end{subfigure}
    \hfill
    \begin{subfigure}{0.32\textwidth}
        \centering \includegraphics[width=\textwidth]{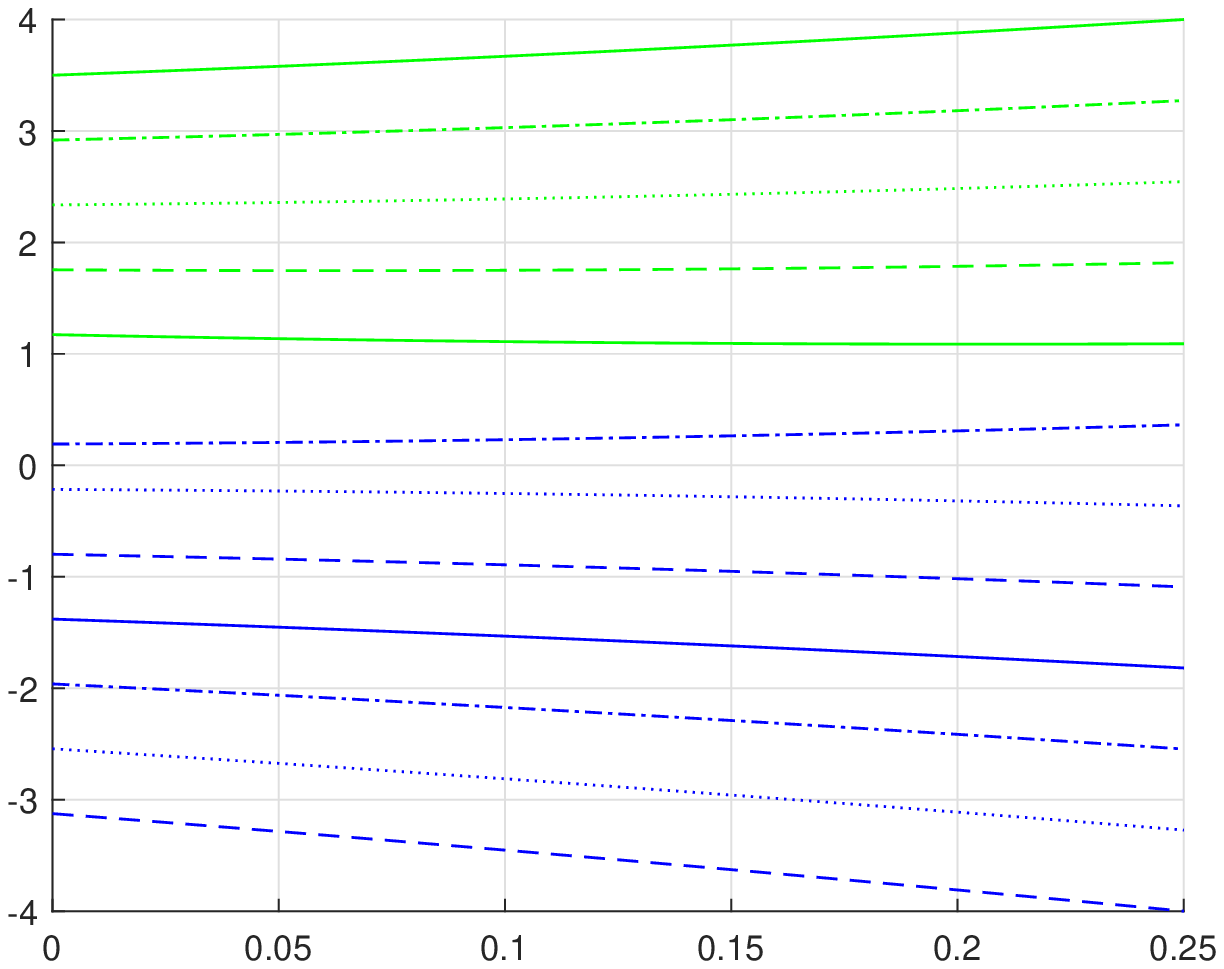}
        \caption{First component, $t=0.25$}
    \end{subfigure}
    \hfill
    \begin{subfigure}{0.32\textwidth}
        \centering \includegraphics[width=\textwidth]{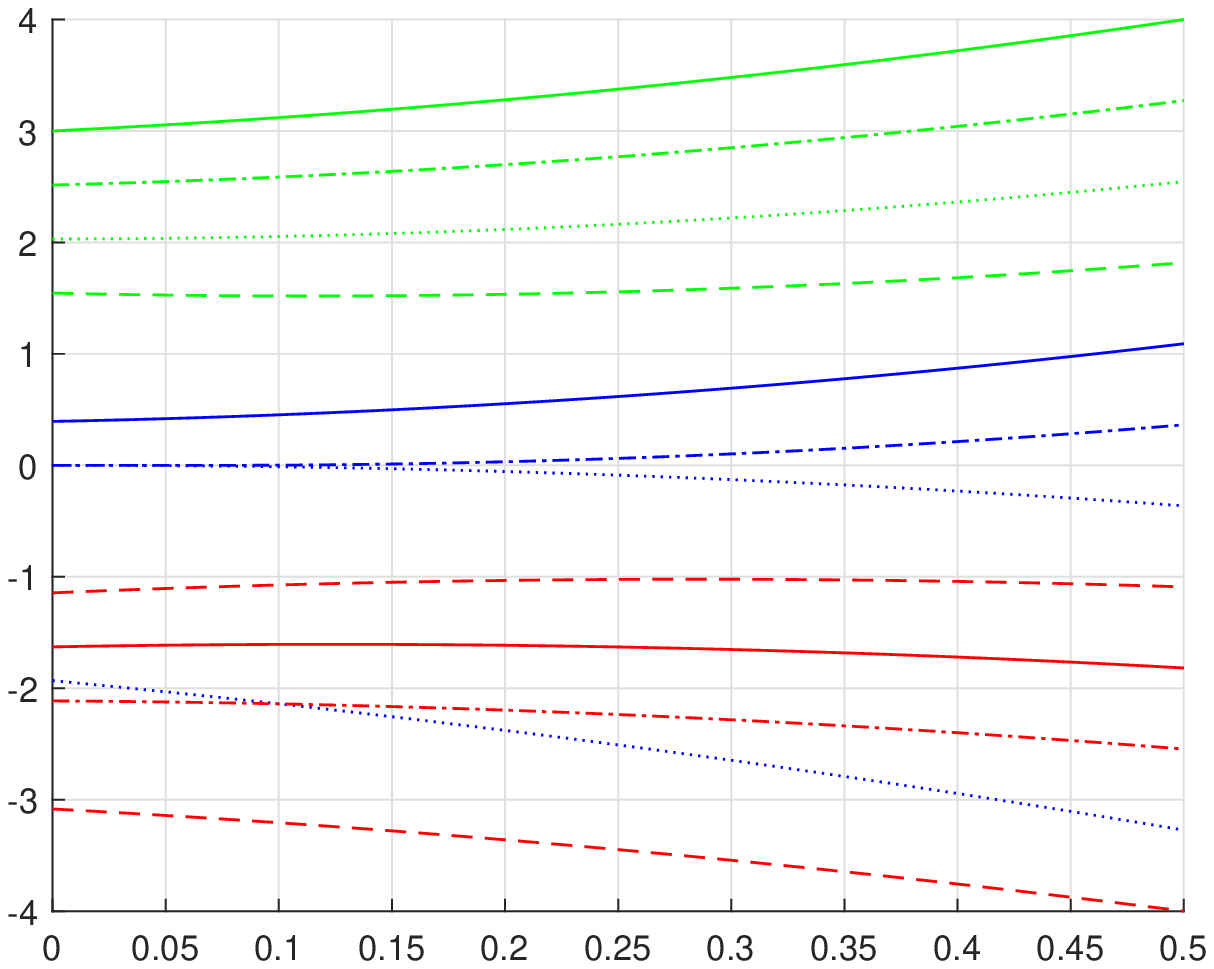}
        \caption{First component, $t=0.5$}
    \end{subfigure}
    
    \begin{subfigure}{0.32\textwidth}
        \centering \includegraphics[width=\textwidth]{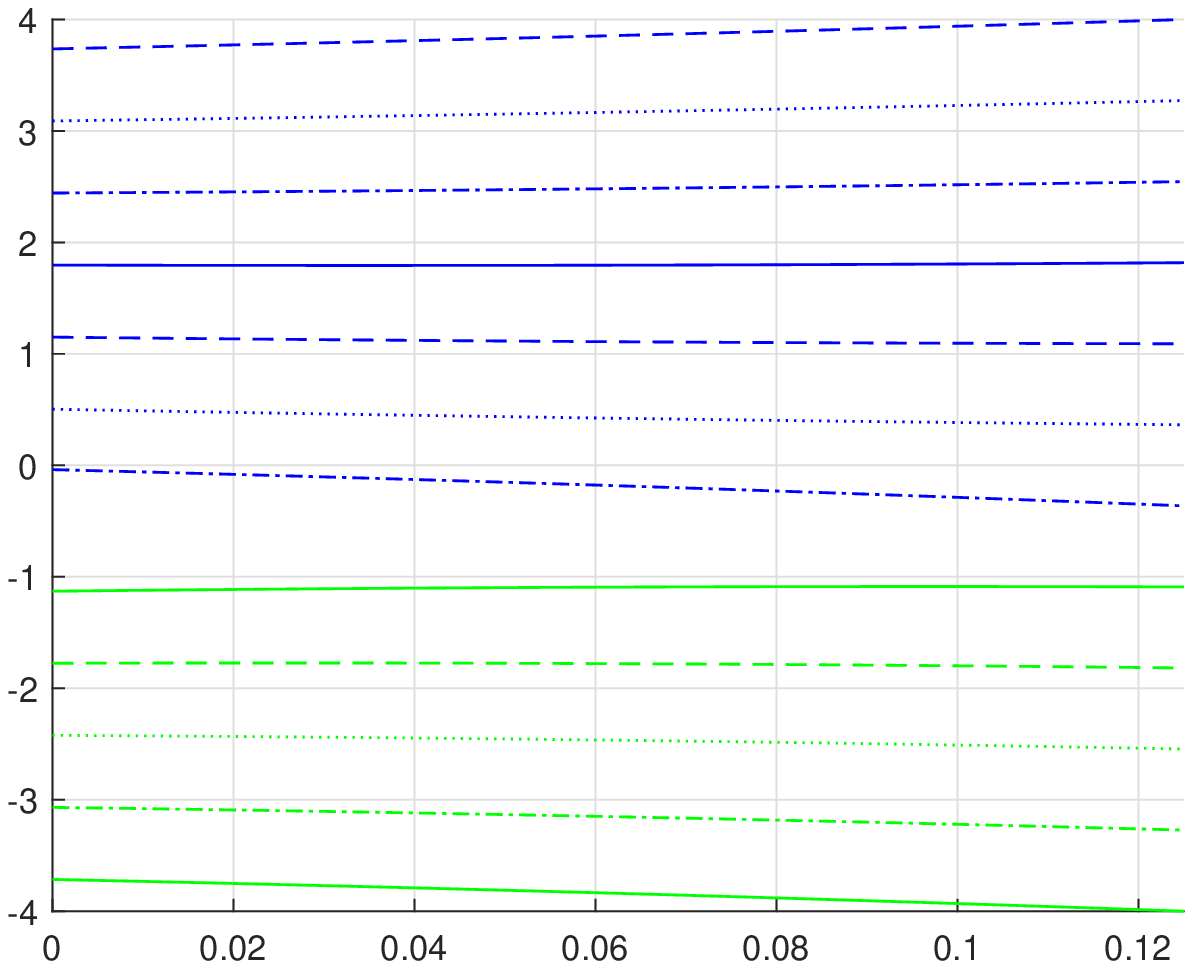}
        \caption{Second component, $t=0.125$}
    \end{subfigure}
    \hfill
    \begin{subfigure}{0.32\textwidth}
        \centering \includegraphics[width=\textwidth]{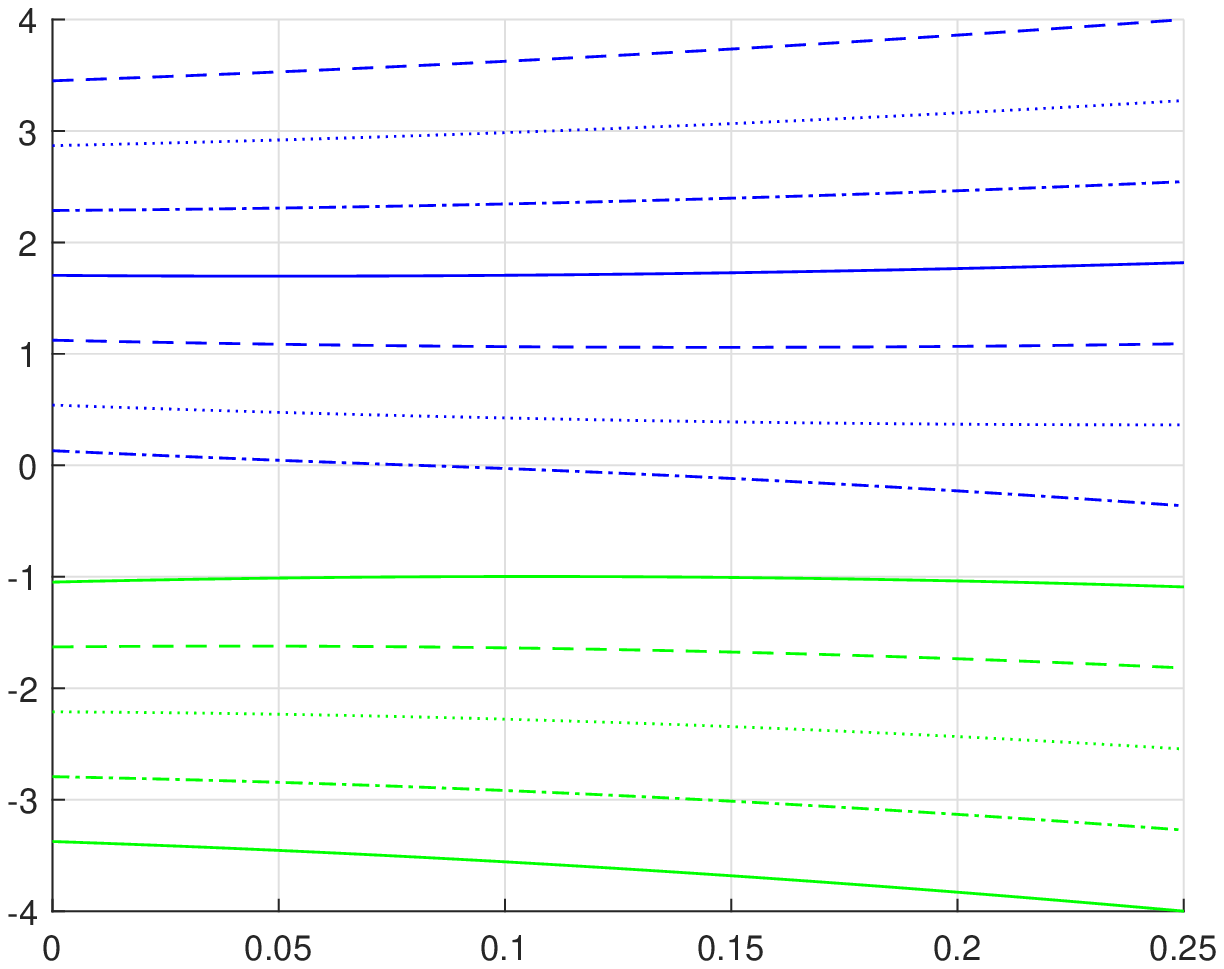}
        \caption{Second component, $t=0.25$}
    \end{subfigure}
    \hfill
    \begin{subfigure}{0.32\textwidth}
        \centering \includegraphics[width=\textwidth]{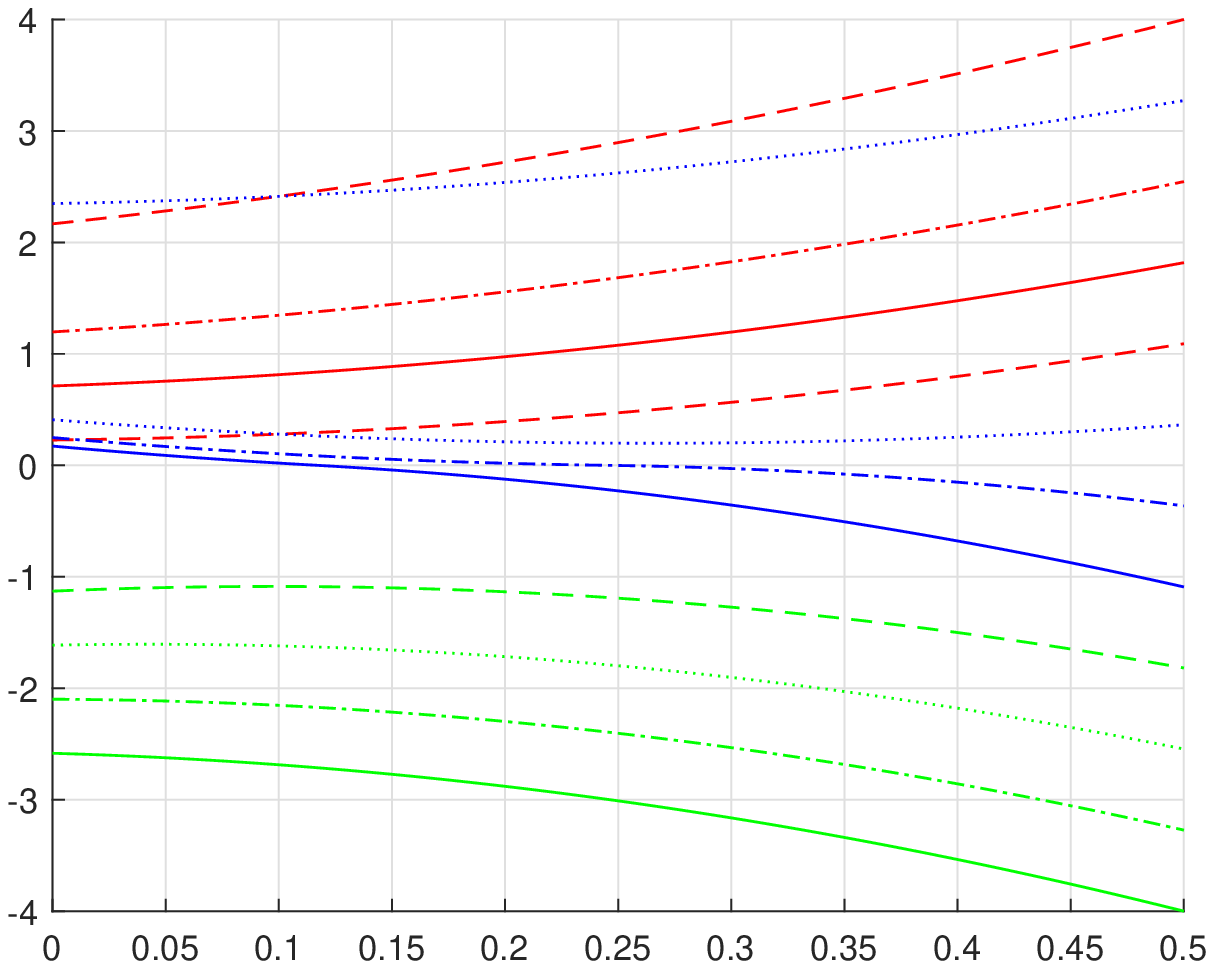}
        \caption{Second component, $t=0.5$}
    \end{subfigure}
    
    \caption{Evaluation of an optimal trajectory $\gmnum(s;(x,-x,0,\dots,0), t)$ of the $10$-dimensional optimal control problem (\ref{eqt: result_optctrl1_hd}) with $\ba = ( 4, 6, 5, \dots, 5)$, $\bb = ( 3, 9, 6, \dots, 6)$, and initial cost $\initcond(\bx) = \min_{j\in\{1,2,3\}} \initcond_j(\bx) = \min_{j\in\{1,2,3\}} \{\frac{1}{2}\|\bx - \by_j\|^2 + \alpha_j\}$, where $\by_1 = (-2,0,\dots,0)$, $\by_2 = (2,-2,-1,0, \dots, 0)$, $\by_3 = (0, 2, 0, \dots, 0)$, $\alpha_1 = -0.5$, $\alpha_2 = 0$, and $\alpha_3 = -1$, versus $s\in[0,t]$ for different terminal positions $(x,-x,0,\dots,0)$ ($x\in[-4,4]$) and different time horizons $t$. The color of the lines denotes which initial cost was used, i.e., $r \in \argmin_{j\in\{1,2,3\}} \Snum_j(\bx,t)$ for $r = 1,2,3$ corresponds to {\textcolor{red}{red}}, {\textcolor{green}{green}}, and {\textcolor{blue}{blue}}, respectively. The different line markers simply differentiate between the different trajectories. Figures (a)-(c) depict the first component of the trajectory versus $s\in[0,t]$ with different time horizons $t$, while (d)-(f) depict the second component of the trajectory versus $s\in[0,t]$ with different time horizons $t$. Plots for time horizons $t = 0.125$, $0.25$, and $0.5$ are depicted in (a)/(d), (b)/(e), and (c)/(f), respectively.}
    \label{fig: HJ1_HD_quad_trajectory}
\end{figure}

In Figure \ref{fig: HJ1_HD_quad_trajectory}, we show several one-dimensional slices of an optimal trajectory $\gmnum(s;(x,-x,0,\dots,0), t)$ of the corresponding optimal control problem with $\ba$ and $\bb$ defined above and initial cost $\initcond$ defined in~\eqref{eqt: hd_initial_data_minofquad}, for different terminal positions $(x,-x,0,\dots,0)$ and different time horizons $t$. We observe that the one-dimensional slices are piecewise quadratic and continuous in $s$, which is consistent with our formulas for $\opttraj$ as
given by~\eqref{eqt: optctrl_defx_1},~\eqref{eqt: optctrl_defx_2},~\eqref{eqt: optctrl_defx_3},~\eqref{eqt: optctrl_defx_4}, and~\eqref{eqt: optctrl_defx_5}.

In Table \ref{tab:timing_minplus}, we present the running time of this example for different dimensions $n$. From Table~\ref{tab:timing_minplus}, we see that it takes less than $5\times 10^{-6}$ seconds on average to compute the solution at one point in a 16-dimensional problem, which demonstrates the efficiency of our proposed algorithm even in high dimensions.

\begin{table}[htbp]
    \centering
    \begin{tabular}{c|c|c|c|c}
        \hline
        \textbf{$\mathbf{n}$} & 4 & 8 & 12 & 16 \\
        \hline
        \textbf{running time (s)} & 9.9695e-07 & 2.2075e-06 & 3.3748e-06 & 4.4818e-06 \\
        \hline
    \end{tabular}
    \hfill 
    
    \caption{Time results in seconds for the average time per call over $102,400$ trials for evaluating the solution of the HJ PDE (\ref{eqt: result_HJ1_hd}) with initial cost $\initcond(\bx) = \min_{j\in\{1,2,3\}} \initcond_j(\bx) = \min_{j\in\{1,2,3\}} \{\frac{1}{2}\|\bx - \by_j\|^2 + \alpha_j\}$, where $\by_1 = (-2,0,\dots,0)$, $\by_2 = (2,-2,-1,0, \dots, 0)$, $\by_3 = (0, 2, 0, \dots, 0)$, $\alpha_1 = -0.5$, $\alpha_2 = 0$, and $\alpha_3 = -1$ for various dimensions $n$.}
    \label{tab:timing_minplus}
\end{table}
\subsection{An FPGA implementation}\label{sec:fpga}

In this section, we describe an FPGA implementation of our building block from Section \ref{subsec:numerical_quad}, i.e., our efficient solver for computing~\eqref{eqt: result_Hopf1_hd}, or equivalently~\eqref{eqt:opt_quad_hd}, exactly. Specifically, we present an FPGA implementation with high throughput. Throughput refers to the amount of data that can be processed in a given amount of time. We achieve a high throughput by designing an implementation with an iteration interval (II) of 1, which means that we can begin processing a new input (in our case, a new $(x_i, t, y_i)\in\R\times[0,\infty)\times\R$ as defined in~\eqref{eqt:opt_quad_hd}) at every FPGA clock cycle. 

FPGAs consist of an array of programmable logic blocks and memory elements connected via reconfigurable interconnects. One of the main constraints in designing an FPGA implementation is the amount of logic resources available on a given FPGA. These resources include general purpose logics such as flip flops (FFs) and lookup tables (LUTs), specialized arithmetic units such as digital signal processing units (DSPs), and memory such as Block Random Access Memory (BRAMs). For a brief overview of FPGAs, we refer the reader to \cite{KastnerFPGA}. 

Our FPGA implementation uses a Xilinx Alveo U280 board with a target design running at 300 MHz. Thus, having an II of 1 means that we can begin processing a new input (i.e., a new $(x_i, t, y_i)$) every 3.3333 nanoseconds. Note that the Alveo U280 board consists of three ``chiplets". Crossing chiplets consumes limited routing resources, which can severely reduce the performance. To avoid this issue, we aim for an FPGA design that uses less than 30\% of any given FPGA resource type to ensure that no chiplet is crossed.

Table \ref{tab:quad_fpga_resources} shows the amount of FPGA resources and the latencies used to implement and run our building block on a Xilinx Alveo U280 board for various dimensions $n$ and $102,400$ points $(\bx, t, \by)\in\Rn\times[0,\infty)\times\Rn$. We observe that since our design streams the points $(\bx, t, \by)$ elementwise (i.e., our FPGA kernel takes the inputs $(x_i, t, y_i)$), the latency scales linearly in the dimension $n$ and the amount of FPGA resources used remains essentially constant in $n$. We also note that we use less than 30\% of the FPGA resources available on the Xilinx Alveo U280 board, which implies that we could either parallelize by implementing multiple copies of our FPGA kernel to maximize usage of the FPGA board or use a smaller (i.e., cheaper) FPGA to implement our building block with similar performance as we report here. 

In Table \ref{tab:quad_fpga_vs_cpu}, we  highlight the performance boost that can be achieved using FPGAs by comparing the performance of a CPU implementation of our building block using C++ to our FPGA implementation. The CPU implementation used here is identical to that used in Section~\ref{subsec:numerical_quad} but using a fixed number of Newton iterations (see Appendix \ref{sec:appendix_numerical_prox} for more details) for better comparability with our FPGA implementation. We observe that our FPGA implementation has a speedup of about 37 to 40 (depending on the dimension $n$) compared to the CPU implementation. Note that using a Xilinx Alveo U280 board, we could parallelize our FPGA implementation for a further speedup of $\times3$ (i.e., a total speedup of about 111 to 122 overall depending on the dimension $n$) by simply replicating our FPGA kernel on each chiplet of the board. This parallelization would not experience any performance degradation due to chiplet crossing, as each copy of the FPGA kernel would be independent, and hence, no kernels/chiplets would need to communicate (as demonstrated in Table \ref{tab:quad_fpga_resources}, each copy of the FPGA kernel would fit fully within a single chiplet since our FPGA kernel uses less than 30\% of the available FPGA resources).

\begin{table}[htbp]
    \centering
    \begin{tabular}{c|c|c|c|c|c}
        \hline
        \textbf{$\mathbf{n}$} & \textbf{Latency (ns)} & \textbf{BRAM} & \textbf{DSPs} & \textbf{FFs} & \textbf{LUTs} \\
        \hline
        4 & 410,604 (1.369e06) & 0 (0\%) & 2,042 (22\%) & 418,105 (16\%) & 154,471 (11\%) \\
        8 & 820,218 (2.734e06) & 0 (0\%) & 2,042 (22\%) & 418,174 (16\%) & 154,442 (11\%) \\
        12 & 1,229,826 (4.099e06) & 0 (0\%) & 2,042 (22\%) & 418,434 (16\%) & 154,465 (11\%) \\
        16 & 1,639,438 (5.464e06) & 0 (0\%) & 2,042 (22\%) & 418,887 (16\%) & 154,611 (11\%) \\
        \hline
    \end{tabular}
    \hfill 
    
    \caption{FPGA resources and latencies in cycles and nanoseconds (ns) for evaluating the solution of the HJ PDE (\ref{eqt: result_HJ1_hd}) with quadratic initial condition~\eqref{eqt: quadratic_IC} at 102,400 points $(\bx,t,\by)\in\Rn\times[0,\infty)\times\Rn$ for various dimensions $n$ using double precision floating points on a Xilinx Alveo U280 board with a frequency of 300 MHz.}
    \label{tab:quad_fpga_resources}
\end{table}

\begin{table}[htbp]
    \centering
    \begin{tabular}{c|c|c|c}
        \hline
        \textbf{$\mathbf{n}$} & \textbf{CPU time (s)} & \textbf{FPGA time (s)} & \textbf{Speedup}  \\
        \hline
        4 & 4.9845e-07 & 1.3369e-08 & 37.2840 \\
        8 & 9.9128e-07 & 2.6699e-08 & 37.1280\\
        12 & 1.6092e-06 & 4.0029e-08 & 40.2009 \\
        16 & 2.1701e-06 & 5.3359e-08 & 40.6698 \\
        \hline
    \end{tabular}
    \hfill 
    
    \caption{Comparison of the average time per call over 102,400 runs for evaluating the solution of the HJ PDE (\ref{eqt: result_HJ1_hd}) with quadratic initial condition~\eqref{eqt: quadratic_IC} for various dimensions $n$ using a CPU implementation on a single Intel Core i5-8250U versus an FPGA implementation on a Xilinx Alveo U280 board with a frequency of 300 MHz.}
    \label{tab:quad_fpga_vs_cpu}
\end{table}

\section{Summary}\label{sec: conclusion}
In this paper, we present analytical solutions to certain optimal control problems with quadratic running costs on the velocity and certain piecewise affine convex running costs on the trajectory. Moreover, we present a Hopf-type representation formula for the corresponding HJ PDE with quadratic kinetic energy and piecewise affine non-positive concave potential function. We also present efficient algorithms for solving these problems with convex initial costs and certain non-convex initial costs. We demonstrate that our algorithms do not suffer from the curse of dimensionality and have promising speedup when implemented on FPGAs in comparison to CPUs. A possible future direction is to combine the proposed algorithms with other building blocks, such as the solver in~\cite{dower2019game} and/or the LQR solver, to handle more general optimal control problems.
\section{Declaration of Competing Interests}

The authors declare that they have no known competing financial interests or personal relationships that could have influenced or appeared to have influenced the work reported in this paper.

\begin{acknowledgements}
This research is supported by  DOE-MMICS SEA-CROGS DE-SC0023191, NSF 1820821, and AFOSR MURI FA9550-20-1-0358. P.C. is supported by the SMART Scholarship, which is funded by USD/R\&E (The Under Secretary of Defense-Research and Engineering), National Defense Education Program (NDEP) / BA-1, Basic Research.
\end{acknowledgements}

\appendix
\newcommand{\yki}{y^k_i}
\newcommand{\pki}{p^k_i}
\newcommand{\xiki}{\xi^k_i}

\section{Some technical lemmas for section~\ref{sec: theory}}

\begin{lem} \label{lem: prop_S}
Let $\valuefn$ be the function defined in~\eqref{eqt: result_S1_1d} and~\eqref{eqt: result_S1_1d_negative}. Let $x\in\R$ and $t,a,b>0$. Then, the function $\R\ni p\mapsto -\valuefn(x,t;p,a,b)\in\R$ is strictly convex, 1-coercive, and continuously differentiable.
\end{lem}
\begin{proof}
We first compute the derivatives of $f_1,f_2, f_3, f_4, f_5$ with respect to $p$ to obtain
\begin{equation}\label{eqt: HJ1_p_deriv}
\begin{split}
    \frac{\partial f_1(x,t;p,a,b)}{\partial p} &= -\frac{at^2}{2} - pt +x,\\
    \frac{\partial f_2(x,t;p,a,b)}{\partial p} &= \frac{bt^2}{2} - pt +x,
    \\
    \frac{\partial f_3(x,t;p,a,b)}{\partial p} &= \frac{a+b}{(a+2b)^2}\left(-(bt-p)^2 + (p-bt)\sqrt{(bt-p)^2 + 2x(a+2b)}\right) + \frac{bx}{a+2b} + \frac{bt^2}{2} - pt,
    \\
    \frac{\partial f_4(x,t;p,a,b)}{\partial p} &= \frac{\partial f_5(x,t;p,a,b)}{\partial p} = - \frac{p^2}{2b}.
\end{split}
\end{equation}
From the above formulas, it is clear that $f_1, \dots, f_5$ are continuously differentiable with respect to $p$ in their domains. Now, we compute their second-order derivatives, which read as follows:
\begin{equation} \label{eqt: lem21_second_derivative}
\begin{split}
    \frac{\partial^2 f_1(x,t;p,a,b)}{\partial p^2} &= \frac{\partial^2 f_2(x,t;p,a,b)}{\partial p^2} = -t <0,
    \\
    \frac{\partial^2 f_3(x,t;p,a,b)}{\partial p^2} &= \frac{2(a+b)}{(a+2b)^2}\left(bt-p +\frac{(bt-p)^2 + x(a+2b)}{\sqrt{(bt-p)^2 + 2x(a+2b)}} \right) - t,
    \\
    \frac{\partial^2 f_4(x,t;p,a,b)}{\partial p^2} &= \frac{\partial^2 f_5(x,t;p,a,b)}{\partial p^2} = - \frac{p}{b} \leq 0.
\end{split}
\end{equation}
Hence, the second-order derivatives of $f_1$ and $f_2$  with respect to $p$ are negative, which implies that $f_1$ and $f_2$ are strongly concave with respect to $p$ in their domains. Note that the domains of $f_4$ and $f_5$ with respect to $p$ are included in $[0,+\infty)$, and hence, their second-order derivatives with respect to $p$ are negative almost everywhere. The strict concavity of $f_4$ and $f_5$ follows. Now, we prove the statement by considering different cases.

First, consider the case where $x\in[0,\frac{at^2}{2}]$. In this case, according to the definitions~\eqref{eqt: result_S1_1d} and~\eqref{eqt: result_S1_1d_negative}, the function $p\mapsto \valuefn(x,t;p,a,b)$ reads
\begin{equation} \label{eqt:lemA1_pf_S_case1}
    \valuefn(x,t;p,a,b) = \begin{dcases}
    f_2(-x,t;-p,b,a) & p< \sqrt{2ax}-at,\\
    f_5(-x,t;-p,b,a) & \sqrt{2ax}-at \leq p< 0,\\
    f_4(x,t;p,a,b) & 0\leq p < b\left(t - \sqrt{\frac{2x}{a}}\right),\\
    f_3(x,t;p,a,b) & p\geq b\left(t - \sqrt{\frac{2x}{a}}\right).
    \end{dcases}
\end{equation}
By straightforward calculation, the function $p\mapsto \valuefn(x,t;p,a,b)$ is continuously differentiable. Now, we prove the 1-coercivity of $-\valuefn$ by computing the limit of the derivatives as $p$ approaches $+\infty$ or $-\infty$. As $p$ approaches $-\infty$, we have
\begin{equation} \label{eqt: lem21_pf_lim1}
\begin{split}
    \lim_{p\to -\infty} \frac{\partial \valuefn(x,t; p,a,b)}{\partial p} &= -\lim_{q\to +\infty} \frac{\partial f_2(-x,t; q,b,a)}{\partial q} = -\lim_{q\to +\infty} \left\{\frac{at^2}{2} - qt - x\right\} = +\infty,
\end{split}
\end{equation}
where the change of variable $q = -p$ is performed to obtain the first equality. Similarly, when $p$ approaches $+\infty$, we obtain
\begin{equation}\label{eqt: lem21_pf_lim2}
\begin{split}
    &\lim_{p\to +\infty} \frac{\partial \valuefn(x,t; p,a,b)}{\partial p} =
    \lim_{p\to +\infty} \frac{\partial f_3(x,t; p,a,b)}{\partial p}\\
    =\ &\lim_{p\to +\infty} \frac{a+b}{(a+2b)^2}\left(-(bt-p)^2 + (p-bt)\sqrt{(bt-p)^2 + 2x(a+2b)}\right) + \frac{bx}{a+2b} + \frac{bt^2}{2} - pt\\
    = \ &\lim_{p\to +\infty} \frac{a+b}{(a+2b)^2}(p-bt)\frac{2x(a+2b)}{\sqrt{(bt-p)^2 + 2x(a+2b)} - bt + p} + \frac{bx}{a+2b} + \frac{bt^2}{2}- pt\\
    = \ &\lim_{p\to +\infty} \frac{2x(a+b)}{a+2b}\frac{1}{\sqrt{1 + \frac{2x(a+2b)}{(bt-p)^2}} +1} + \frac{bx}{a+2b} + \frac{bt^2}{2}- pt = -\infty,
\end{split}
\end{equation}
where the last equality holds since we assume $t>0$. By~\eqref{eqt: lem21_pf_lim1} and~\eqref{eqt: lem21_pf_lim2}, we conclude that $p\mapsto -\valuefn(x,t;p,a,b)$ is 1-coercive when $x\in[0,\frac{at^2}{2}]$. 
It remains to prove the strict concavity of $\valuefn$ with respect to $p$ in this case. Since $\valuefn$ is defined piecewise and continuously differentiable, it suffices to prove that each piece is strictly concave. Recall that we proved the strict concavity of $f_2,f_4,f_5$ using~\eqref{eqt: lem21_second_derivative}, and hence, the first three lines in~\eqref{eqt:lemA1_pf_S_case1} are strictly concave functions.
It remains to consider $f_3$. Some computation shows that $\frac{\partial^2 f_3(x,t;p,a,b)}{\partial p^2}\leq 0$ holds if and only if there holds
\begin{equation}\label{eqt: lem21_pf_f3_ineqt}
    (a+b)^2x^2 - t\left((a+b)p+\frac{a^2t}{4}\right)((p-bt)^2 + 2x(a+2b))\leq 0.
\end{equation}
The left-hand side in~\eqref{eqt: lem21_pf_f3_ineqt} is a second-order polynomial with respect to $x$ with positive leading coefficient. As a result, to check that the inequality in~\eqref{eqt: lem21_pf_f3_ineqt} holds for all $x\in[0,\frac{at^2}{2}]$, it suffices to prove that the inequality holds at $x=0$ and $x=\frac{at^2}{2}$. Denote the left-hand side of~\eqref{eqt: lem21_pf_f3_ineqt} by $P(x)$. Then, we have that
\begin{equation} \label{eqt: lem21_pf_P_endpt}
    \begin{split}
    P(0) &= -t\left((a+b)p+\frac{a^2t}{4}\right)(p-bt)^2 \leq 0,\\
    P\left(\frac{at^2}{2}\right) &= -tp\left((a+b)(bt-p)^2 + \frac{a^2pt}{4} + \left(\left(a+\frac{b}{2}\right)(a^2+2ab)+ab^2\right)t^2\right)\leq 0,
    \end{split}
\end{equation}
where $P(0)=0$ holds if and only if $p=bt$ and where $P(\frac{at^2}{2})=0$ holds if and only if $p=0$. In other words, $\frac{\partial^2 f_3(x,t;p,a,b)}{\partial p^2}\leq 0$ holds in its corresponding domain $p>bt-b\sqrt{\frac{2x}{a}}$ except at finitely many points, which implies the strict concavity of $\valuefn$ with respect to $p$ in this domain. Therefore, we conclude that the function $p\mapsto -\valuefn(x,t;p,a,b)$ is strictly convex, 1-coercive, and continuously differentiable for $x\in [0,\frac{at^2}{2}]$.

Next, we consider the case where $x> \frac{at^2}{2}$. In this case, the function $p\mapsto \valuefn(x,t;p,a,b)$ reads:
\begin{equation}\label{eqt:lemA1_pf_S_case2}
    \valuefn(x,t;p,a,b) = \begin{dcases}
     f_2(-x,t;-p,b,a) & p < 0,\\
     f_1(x,t;p,a,b) & 0\leq p\leq \frac{x}{t} - \frac{at}{2},\\
     f_3(x,t;p,a,b) & p> \frac{x}{t} - \frac{at}{2}.
    \end{dcases}
\end{equation}
By straightforward calculation, $p\mapsto \valuefn(x,t;p,a,b)$ is continuously differentiable. Note that \eqref{eqt: lem21_pf_lim1} and~\eqref{eqt: lem21_pf_lim2} still hold in this case, and hence, the function $p\mapsto -\valuefn(x,t;p,a,b)$ is 1-coercive. As in the first case, the strict concavity of the first two lines in~\eqref{eqt:lemA1_pf_S_case2} follows from~\eqref{eqt: lem21_second_derivative}. Hence, to prove the strict concavity of $p\mapsto \valuefn(x,t;p,a,b)$, it suffices to show that $\frac{\partial^2 f_3(x,t;p,a,b)}{\partial p^2}\leq 0$ when $x> \frac{at^2}{2}$ and $p>\frac{x}{t} - \frac{at}{2}$, i.e., it suffices to check~\eqref{eqt: lem21_pf_f3_ineqt} for all $x\in (\frac{at^2}{2}, tp + \frac{at^2}{2})$. Again, denote the left-hand side of~\eqref{eqt: lem21_pf_f3_ineqt} by $P(x)$. Since $P(x)$ is a second-order polynomial with respect to $x$ with positive leading coefficient, it suffices to show that $P(x)\leq 0$ at $x=\frac{at^2}{2}$ and $x=tp + \frac{at^2}{2}$. According to~\eqref{eqt: lem21_pf_P_endpt}, $P(\frac{at^2}{2})\leq 0$. After some calculations, we get that
\begin{equation*}
    P\left(tp + \frac{at^2}{2}\right) = -pt^3\left(b(a+b)^2 + \frac{a^2}{2}(a+b)\right) - p^2t^2 \left((a+b)^2 + \frac{a^2}{4}\right) - (a+b)p^3t \leq 0,
\end{equation*}
where the inequality holds since we have $a,b,p, t>0$. Therefore, for $p>\frac{x}{t} - \frac{at}{2}$, the function $p\mapsto f_3(x,t;p,a,b)$ is strictly concave.
As a result, we have shown that the function $p\mapsto -\valuefn(x,t;p,a,b)$ is strictly convex, 1-coercive, and continuously differentiable for $x\in (\frac{at^2}{2}, +\infty)$.

Now, we consider the case where $x<0$.
By~\eqref{eqt: result_S1_1d_negative}, for all $x,p\in\R$, $t\geq 0$, 
we have $\valuefn(x,t;p,a,b) = \valuefn(-x, t; -p,b,a)$. Therefore, the conclusion for $x<0$ also holds since we have already proved that the function $-p\mapsto \valuefn(-x, t; -p,b,a)$ is strictly convex, 1-coercive, and continuously differentiable.
\end{proof}

\begin{lem} \label{prop: S_linearJ}
Let $\initcond\colon \R\to\R$ be defined by $\initcond(x)=px$ for some $p \in\R$. Let $\potentfn\colon \R\to\R$ be the function defined by~\eqref{eqt: result_defV} for some constants $a,b>0$. Let $\valuefn$ be the function defined in~\eqref{eqt: result_S1_1d} and~\eqref{eqt: result_S1_1d_negative}.
Then, the following statements hold:
\begin{enumerate}
    \item[(a)] The function $(x,t)\mapsto \valuefn(x,t;p,a,b)$ is a continuously differentiable solution to the corresponding HJ PDE~\eqref{eqt: result_HJ1_1d}.
    \item[(b)] For all $t\geq 0$, the function $x\mapsto \valuefn(x,t;p,a,b)$ is convex.
\end{enumerate}
\end{lem}

\begin{proof}
(a) We first prove the statement for non-negative $p$. Assume $p\geq 0$. By straightforward calculation, the derivatives of $f_3$ with respect to $x$ and $t$ read:
\begin{equation}\label{eqt:propA1_pf_df_dxt_1}
    \begin{split}
        \frac{\partial f_3(x,t;p, a,b)}{\partial t} &= \frac{b(a+b)}{(a+2b)^2} \left((bt-p)^2 + (bt-p)\sqrt{(bt-p)^2 + 2x(a+2b)}\right) \\
        &\quad - \frac{b^2x}{a+2b} - \frac{b^2t^2}{2} + bpt - \frac{p^2}{2}, \\
    \frac{\partial f_3(x,t;p, a,b)}{\partial x} &= \frac{a+b}{a+2b} \sqrt{(bt-p)^2 + 2x(a+2b)} - \frac{b(bt-p)}{a+2b}.
    \end{split}
\end{equation}
The derivatives of the functions $f_1,f_2,f_4,f_5$ read:
\begin{equation}\label{eqt:propA1_pf_df_dxt_2}
\begin{aligned}
    \frac{\partial f_1(x,t;p, a,b)}{\partial t} &= -\frac{a^2t^2}{2} - apt + ax - \frac{p^2}{2}, \quad\quad\quad\quad & \frac{\partial f_1(x,t;p, a,b)}{\partial x} &= at + p,\\
    \frac{\partial f_2(x,t;p, a,b)}{\partial t} &= -\frac{b^2t^2}{2} + bpt - bx - \frac{p^2}{2}, & \frac{\partial f_2(x,t;p, a,b)}{\partial x} &= -bt + p,\\
    \frac{\partial f_4(x,t;p, a,b)}{\partial t} &= 0, &
    \frac{\partial f_4(x,t;p, a,b)}{\partial x} &= \sqrt{2ax},\\
    \frac{\partial f_5(x,t;p, a,b)}{\partial t} &= 0, &
    \frac{\partial f_5(x,t;p, a,b)}{\partial x} &= -\sqrt{-2bx}.
\end{aligned}
\end{equation}
It is straightforward to check that these five functions all satisfy the differential equation in~\eqref{eqt: result_HJ1_1d}. 
By straightforward calculation, the function $(x,t)\mapsto \valuefn(x,t;p,a,b)$ is continuously differentiable and satisfies the initial condition in~\eqref{eqt: result_HJ1_1d}. Therefore, $\valuefn$ is a continuously differentiable solution to the HJ PDE~\eqref{eqt: result_HJ1_1d}.

Next, we consider negative $p$. Assume $p<0$. For all $\alpha, \beta >0$, denote by $\potentfn_{\alpha,\beta}(x)$ the function $\potentfn$ defined in~\eqref{eqt: result_defV} with constants $a=\alpha$ and $b=\beta$. Note that there holds $\potentfn_{\alpha,\beta}(x) = \potentfn_{\beta, \alpha}(-x)$ for all $x\in\R$.
By~\eqref{eqt: result_S1_1d_negative}, we have 
\begin{equation*}
\begin{split}
    &\frac{\partial \valuefn(x,t;p,a,b)}{\partial t} + \frac{1}{2}(\nabla_x \valuefn(x,t;p,a,b))^2 + \potentfn_{a,b}(x)\\
    =\ & \frac{\partial \valuefn(-x,t;-p,b,a)}{\partial t} + \frac{1}{2}(-\nabla_x \valuefn(-x,t;-p,b,a))^2 + \potentfn_{b,a}(-x)\\
    =\ & \frac{\partial \valuefn(y,t;-p,b,a)}{\partial t} + \frac{1}{2}(\nabla_x \valuefn(y,t;-p,b,a))^2 + \potentfn_{b,a}(y)\\
    = \ & 0,
\end{split}
\end{equation*}
where the second equality follows from the change of variable $y=-x$ and the last equality holds since the function $(y,t)\mapsto \valuefn(y,t;-p,b,a)$ is a solution to the HJ PDE~\eqref{eqt: result_HJ1_1d} with potential energy $\potentfn_{b,a}$, according to the proof above for positive $p$.
We check the initial condition as follows:
\begin{equation*}
    \valuefn(x,0; p,a,b) = \valuefn(-x,0; -p,b,a) = (-p)(-x)= px. 
\end{equation*}
Therefore, the function $(x,t)\mapsto \valuefn(x,t;p,a,b)$ defined in~\eqref{eqt: result_S1_1d_negative} solves~\eqref{eqt: result_HJ1_1d}. It is continuously differentiable since it equals $\valuefn(-x,t;-p,b,a)$, which is continuously differentiable with respect to $(-x,t)$.

(b) When $t=0$, we have $\valuefn(x,t;p,a,b) = px$, which is convex with respect to $x$. It remains to consider the case when $t>0$. Recall that we have proved above that the function $\valuefn$ is continuously differentiable. To show the convexity of $\valuefn$ with respect to $x$, it suffices to show the convexity of each $f_i$ in its domain. Let $p\geq 0$. By straightforward calculation, we obtain
\begin{equation*}
\begin{aligned}
\frac{\partial^2 f_1(x,t;p, a,b)}{\partial x^2} &= \frac{\partial^2 f_2(x,t;p, a,b)}{\partial x^2} = 0,&
\frac{\partial^2 f_4(x,t;p, a,b)}{\partial x^2} &= \frac{a}{\sqrt{2ax}}>0,\\
\frac{\partial^2 f_3(x,t;p, a,b)}{\partial x^2} &= \frac{a+b}{\sqrt{(bt-p)^2 + 2x(a+2b)}} >0,&
\frac{\partial^2 f_5(x,t;p, a,b)}{\partial x^2} &= \frac{b}{\sqrt{-2bx}}>0.
\end{aligned}
\end{equation*}
Therefore, the function $x\mapsto \valuefn(x,t;p,a,b)$ is convex for all $p\geq 0$. The convexity of $x\mapsto \valuefn(x,t;p,a,b)$ for $p<0$ follows from~\eqref{eqt: result_S1_1d_negative}.
\end{proof}

\begin{lem}\label{lem: appendix_S_bd}
Let $a,b>0$ be any positive scalars and $\valuefn$ be the function defined in~\eqref{eqt: result_S1_1d} and~\eqref{eqt: result_S1_1d_negative}. Then, for all $p\in\R$, there exists a constant $C_p\in\R$, such that there holds
\begin{equation} \label{eqt: lemA3_eq}
px + C_p\leq \valuefn(x,t;p,a,b)\leq px + t(a+b)|x| \quad \forall x\in\R, t\geq 0.
\end{equation}
\end{lem}
\begin{proof}
In this proof, whenever there is no ambiguity, we write $\valuefn(x,t;p)$ instead of $\valuefn(x,t;p,a,b)$.

We begin by proving the first, leftmost inequality in~\eqref{eqt: lemA3_eq}. Let $p\geq 0$ be an arbitrary non-negative number. We first prove the statement for this $p$. For this, it suffices to prove that
\begin{equation*}
\inf_{\substack{(x,t)\in\R\times[0,+\infty)\\s.t.\,(x,t,p)\in\Omega_i}} \{\valuefn(x,t;p) - px\} > -\infty, \quad\quad \forall i\in\{1,2,3,4,5\}.
\end{equation*}
For $(x,t,p)\in\Omega_1$, we have $x\geq pt + \frac{at^2}{2}$ and 
\begin{equation*}
\begin{split}
\valuefn(x,t;p) - px &= -\frac{a^2t^3}{6} - \frac{apt^2}{2}  + atx - \frac{p^2 t}{2}
\geq -\frac{a^2t^3}{6} - \frac{apt^2}{2}  + at\left(pt + \frac{at^2}{2}\right) - \frac{p^2 t}{2}\\
&= \frac{a^2t^3}{3} + \frac{apt^2}{2} - \frac{p^2t}{2}.
\end{split}
\end{equation*}
Note that the function $t\mapsto \frac{a^2t^3}{3} + \frac{apt^2}{2} - \frac{p^2t}{2}$ is a third-order polynomial with positive leading coefficient $\frac{a^2}{3}$. Thus, this function is continuous and coercive for $t\in[0,+\infty)$, and hence, it has a finite lower bound.
In other words, we have
\begin{equation*}
\inf_{\substack{(x,t)\in\R\times[0,+\infty)\\s.t.\,(x,t,p)\in\Omega_1}} \{\valuefn(x,t;p) - px\}
\geq \inf_{t\geq 0} \left\{\frac{a^2t^3}{3} + \frac{apt^2}{2} - \frac{p^2t}{2}\right\} > -\infty.
\end{equation*}
For $(x,t,p)\in\Omega_2$ and $t< \frac{p}{b}$, we have $x< 0$, and hence,
\begin{equation*}
\begin{split}
\valuefn(x,t;p) - px &= -\frac{b^2t^3}{6} + \frac{bpt^2}{2} - \frac{p^2 t}{2} - btx
\geq -\frac{b^2t^3}{6} + \frac{bpt^2}{2} - \frac{p^2 t}{2} \\
&\geq -\frac{p^3}{6b} +0 - \frac{p^3}{2b} = - \frac{2p^3}{3b}.
\end{split}
\end{equation*}
For $(x,t,p)\in\Omega_2$ and $t\geq \frac{p}{b}$, we have $x< -\frac{b}{2}(t-\frac{p}{b})^2$, and hence,
\begin{equation*} 
\begin{split}
\valuefn(x,t;p) - px &= -\frac{b^2t^3}{6} + \frac{bpt^2}{2} - \frac{p^2 t}{2} - btx
\geq -\frac{b^2t^3}{6} + \frac{bpt^2}{2} - \frac{p^2 t}{2}  + \frac{b^2t}{2}\left(t-\frac{p}{b}\right)^2\\
&= \frac{b^2t^3}{3} - \frac{bpt^2}{2}.
\end{split}
\end{equation*}
The function $t\mapsto \frac{b^2t^3}{3} - \frac{bpt^2}{2}$ is a third-order polynomial with positive leading coefficient, which implies that this function is bounded from below for $t\in[0,+\infty)$. Therefore, we obtain
\begin{equation*}
\inf_{\substack{(x,t)\in\R\times[0,+\infty)\\s.t.\,(x,t,p)\in\Omega_2}} \{\valuefn(x,t;p) - px\}
\geq \min\left\{- \frac{2p^3}{3b}, \inf_{t\geq 0} \left\{\frac{b^2t^3}{3} - \frac{bpt^2}{2}\right\}\right\} > -\infty.
\end{equation*}
Consider $(x,t,p)\in\Omega_3$ and let $\Delta:= (bt-p)^2 + 2x(a+2b)$.
If $0\leq t< \frac{p}{b}$, we have $0\leq x< pt + \frac{at^2}{2}$. By Lemma~\ref{prop: S_linearJ}, $\valuefn(x,t;p) - px$ is continuous with respect to $(x,t)$, and hence, it is bounded from below in the compact domain $\{(x,t)\in\R^2\colon 0\leq t\leq \frac{p}{b},\,0\leq x\leq pt + \frac{at^2}{2}\}$. 
If $t\geq \frac{p}{b}$, we have $\frac{a}{2}(t-\frac{p}{b})^2\leq x< pt+\frac{at^2}{2}$, and hence, $\Delta \geq (bt-p)^2 + a(a+2b)(t-\frac{p}{b})^2 = (1+\frac{a}{b})^2(bt-p)^2$. Therefore, we obtain
\begin{equation} \label{eqt: lemA3_leftest_3}
\begin{split}
&\valuefn(x,t;p) - px \\
=& 
\frac{a+b}{3(a+2b)^2}\left((bt-p)^3 + \Delta^{3/2}\right)
- \frac{b}{a+2b}(bt-p)x - \frac{b^2t^3}{6} + \frac{bpt^2}{2} - \frac{p^2t}{2} -px\\
\geq& \frac{(a+b)(bt-p)^3}{3(a+2b)^2} \left(1+\left(1+\frac{a}{b}\right)^3\right)
- \frac{(b^2t + (a+b)p)x}{a+2b}- \frac{b^2t^3}{6} + \frac{bpt^2}{2}- \frac{p^2t}{2}\\
\geq& \frac{(a+b)(bt-p)^3}{3(a+2b)^2} \left(1+\left(1+\frac{a}{b}\right)^3\right)
- \frac{b^2t + (a+b)p}{a+2b}\left(pt+\frac{at^2}{2}\right) - \frac{b^2t^3}{6} + \frac{bpt^2}{2}- \frac{p^2t}{2}\\
=:& g(t),
\end{split}
\end{equation}
where the function $g$ defined in the last line is a third-order polynomial with respect to $t$ whose leading coefficient reads:
\begin{equation*}
\begin{split}
    &\frac{(a+b)b^3}{3(a+2b)^2}\left(1+\left(1+\frac{a}{b}\right)^3\right) - \frac{b^2a}{2(a+2b)} - \frac{b^2}{6}\\
    =\,& \frac{(a+b)b^3+(a+b)^4}{3(a+2b)^2} - \frac{3ab^2 + b^2(a+2b)}{6(a+2b)}\\
    =\,& \frac{(a+b)(a^2+ab+b^2) - b^2(2a+b)}{3(a+2b)} = \frac{a^2}{3}>0.
\end{split}
\end{equation*}
Therefore, the function $g$ is bounded from below for $t\in[0,+\infty)$, which implies
\begin{equation*}
\begin{split}
\inf_{\substack{(x,t)\in\R\times[0,+\infty)\\s.t.\,(x,t,p)\in\Omega_3}} \{\valuefn(x,t;p) - px\}
\geq & \min\left\{\inf_{0\leq t\leq \frac{p}{b}, 0\leq x\leq pt + \frac{at^2}{2}} \{\valuefn(x,t;p) - px\}, \inf_{t\geq 0} g(t)\right\}\\
>& -\infty.
\end{split}
\end{equation*}
If $(x,t,p)\in\Omega_4$, we have $x\geq 0$. By~\eqref{eqt:propA1_pf_df_dxt_2}, we obtain
\begin{equation*}
    \frac{\partial}{\partial x}(\valuefn(x,t;p) - px) =\sqrt{2ax} - p\in \begin{dcases}
    [0,+\infty) & x\geq \frac{p^2}{2a},\\
    (-\infty,0] &0\leq x< \frac{p^2}{2a}.
    \end{dcases}
\end{equation*}
Hence, the minimal value of $\valuefn(x,t;p) - px$ is attained at $x= \frac{p^2}{2a}$, and we have
\begin{equation*}
\begin{split}
&\inf_{\substack{(x,t)\in\R\times[0,+\infty)\\s.t.\,(x,t,p)\in\Omega_4}} \{\valuefn(x,t;p) - px\}
\geq \valuefn\left(\frac{p^2}{2a},t;p\right) - \frac{p^3}{2a} = \frac{p^3}{3a} - \frac{p^3}{6b} - \frac{p^3}{2a} >-\infty.
\end{split}
\end{equation*}
If $(x,t,p)\in\Omega_5$, we have $x\leq 0$. By~\eqref{eqt:propA1_pf_df_dxt_2}, we obtain
\begin{equation*}
    \frac{\partial}{\partial x}(\valuefn(x,t;p) - px) = -\sqrt{-2bx} - p \leq 0.
\end{equation*}
Then, the minimal value is attained at $x= 0$, and we have
\begin{equation*}
\begin{split}
&\inf_{\substack{(x,t)\in\R\times[0,+\infty)\\s.t.\,(x,t,p)\in\Omega_5}} \{\valuefn(x,t;p) - px\}
\geq \valuefn\left(0,t;p\right) = - \frac{p^3}{6b}  >-\infty.
\end{split}
\end{equation*}
Therefore, we have $\inf_{x\in\R,t\geq 0}\{\valuefn(x,t;p) - px\} >-\infty$ for all $p\geq 0$, and hence, the first inequality in~\eqref{eqt: lemA3_eq} holds for all $p\geq 0$. If $p<0$, we get
\begin{equation*}
\begin{split}
    \inf_{x\in\R,t\geq 0}\{\valuefn(x,t;p,a,b) - px\}
    &= \inf_{x\in\R,t\geq 0}\{\valuefn(-x,t;-p,b,a) - px\} \\
    &= \inf_{y\in\R,t\geq 0}\{\valuefn(y,t;q,b,a) - qy\} >-\infty,
\end{split}
\end{equation*}
where the first equality holds by~\eqref{eqt: result_S1_1d_negative},
the second equality holds by the change of variables $q=-p$ and $y=-x$, 
and the last inequality holds since we have already proved the first inequality in~\eqref{eqt: lemA3_eq} for the positive $q$ case. Therefore, the first inequality in~\eqref{eqt: lemA3_eq} also holds for all $p<0$.

Now, we prove the second inequality in~\eqref{eqt: lemA3_eq}. We first consider the case where $p\geq 0$. Since the function $x\mapsto px + t(a+b)|x|$ is linear in $[0,+\infty)$ and $(-\infty, 0]$ and the function $x\mapsto \valuefn(x,t;p)$ is convex by Lemma~\ref{prop: S_linearJ}(b), then to prove that $\valuefn(x,t;p)\leq px + t(a+b)|x|$ holds for all $x\in\R$, it suffices to prove it for $x=0$, $x\geq C$, and $x\leq -C$ for some large scalar $C>0$. In other words, it suffices to consider the cases where $x=0$, $(x,t,p)\in\Omega_1$, and $(x,t,p)\in\Omega_2$. Note that for all $p\geq 0$ and $t\leq \frac{p}{b}$, the value $\valuefn(x,t;p)$ at $x=0$ equals the function $f_2(0,t;p,a,b)$, according to the continuity of $\valuefn$. Therefore, we only need to consider the following three cases:
\begin{enumerate}
    \item If $x=0$ and $t>\frac{p}{b}$,
we have
\begin{equation*}
    \valuefn(x,t;p) = \valuefn(0,t;p) =f_4(0,t;p,a,b)= -\frac{p^3}{6b} \leq 0 = px + t(a+b)|x|.
\end{equation*}
\item If $(x,t,p)\in\Omega_1$, we have $x\geq 0$ and 
\begin{equation*}
    \valuefn(x,t;p) =-\frac{a^2t^3}{6} - \frac{apt^2}{2}  + atx - \frac{p^2 t}{2} + px \leq atx + px \leq px + t(a+b)|x|.
\end{equation*}
\item If $(x,t,p)\in\Omega_2$, we have $x\leq 0$ and 
\begin{equation*}
\begin{split}
    \valuefn(x,t;p) &= -\frac{b^2t^3}{6} + \frac{bpt^2}{2} - \frac{p^2 t}{2} - btx + px = px + bt|x| - \frac{t}{6}\left(bt - \frac{3p}{2}\right)^2 - \frac{p^2t}{8} \\
    &\leq px + bt|x| \leq px + t(a+b)|x|.
\end{split}
\end{equation*}
\end{enumerate}
Therefore, the second inequality in~\eqref{eqt: lemA3_eq} holds for all $t,p\geq 0$, and $x\in\R$.
If $p<0$, letting $q=-p$ and $y=-x$, we have
\begin{equation*}
    \valuefn(x,t;p,a,b) = \valuefn(y,t;q,b,a)\leq qy + t(a+b)|y| = px + t(a+b)|x|,
\end{equation*}
where the first equality holds by~\eqref{eqt: result_S1_1d_negative} and the first inequality holds since we have already proved it above for the positive $q$ case.
Therefore, the second inequality in~\eqref{eqt: lemA3_eq} also holds for all $p<0$.
\end{proof}

\begin{lem}\label{lem: appendix_dSdx_bd}
Let $a,b>0$ be positive scalars and $\valuefn$ be the function defined in~\eqref{eqt: result_S1_1d} and~\eqref{eqt: result_S1_1d_negative}. Then, there holds
\begin{equation}\label{eqt: lemA4_bd_dSdx}
    \left|\frac{\partial \valuefn(x,t; p,a,b)}{\partial x}\right| \leq C(R_x, R_t, R_p) < +\infty,
\end{equation}
for all $x\in [-R_x,R_x]$, $t\in[0,R_t]$, and $p\in [-R_p,R_p]$.
\end{lem}
\begin{proof}
We first consider the case where $p\geq 0$.
Let $|x|\leq R_x$, $t\in[0,R_t]$, and $p\in[0,R_p]$. 
We obtain the derivatives of $\valuefn$ with respect to $x$ in~\eqref{eqt:propA1_pf_df_dxt_1} and~\eqref{eqt:propA1_pf_df_dxt_2}.
If $(x,t,p)\in\Omega_1$, we have
\begin{equation*}
\left|\frac{\partial \valuefn(x,t; p,a,b)}{\partial x}\right|
= |at+p|\leq aR_t + R_p =:C_1.
\end{equation*}
If $(x,t,p)\in\Omega_2$, we have
\begin{equation*}
\left|\frac{\partial \valuefn(x,t; p,a,b)}{\partial x}\right|
= |-bt+p|\leq bR_t + R_p =: C_2.
\end{equation*}
If $(x,t,p)\in\Omega_3$, we have
\begin{equation*}
\begin{split}
\left|\frac{\partial \valuefn(x,t; p,a,b)}{\partial x}\right|
&= \left|\frac{(a+b)\sqrt{(bt-p)^2 + 2x(a+2b)}}{a+2b}
- \frac{b}{a+2b}(bt-p)\right|\\
&\leq \frac{(a+b)\sqrt{2b^2R_t^2 + 2R_p^2 + 2(a+2b)R_x} + b^2R_t + bR_p}{a+2b} =: C_3.
\end{split}
\end{equation*}
If $(x,t,p)\in\Omega_4$, we have
\begin{equation*}
    \left|\frac{\partial \valuefn(x,t;p, a,b)}{\partial x}\right| = \sqrt{2ax}\leq \sqrt{2aR_x} =:C_4.
\end{equation*}
If $(x,t,p)\in\Omega_5$, we have
\begin{equation*}
    \left|\frac{\partial \valuefn(x,t;p, a,b)}{\partial x}\right| = \left|-\sqrt{-2bx}\right|\leq \sqrt{2bR_x} =:C_5.
\end{equation*}
Therefore, the bound in~\eqref{eqt: lemA4_bd_dSdx} holds at $(x,t,p)$ for the constant $C(R_x,R_t,R_p)=C_{a,b}$ defined by
\begin{equation*}
    C_{a,b}:= \max\left\{C_1,C_2,C_3,C_4,C_5 \right\}.
\end{equation*}
Now, we consider the case where $p<0$. Let $|x|\leq R_x$, $t\in[0,R_t]$, and $p\in[-R_p,0)$. Let $q=-p$ and $y=-x$. Hence, we have $q\in (0,R_p]$ and $y\in [-R_x,R_x]$. By~\eqref{eqt: result_S1_1d_negative}, we have
\begin{equation*}
    \left|\frac{\partial \valuefn(x,t;p, a,b)}{\partial x}\right| = \left|-\frac{\partial \valuefn(y,t;q, b,a)}{\partial y}\right|\leq C_{b,a},
\end{equation*}
where the inequality was proved in the beginning of this proof since the parameter $q$ is positive.
Therefore, the inequality~\eqref{eqt: lemA4_bd_dSdx} holds for all $x\in [-R_x,R_x]$, $t\in[0,R_t]$, and $p\in [-R_p,R_p]$ with the constant $C(R_x,R_t,R_p):=\max\{C_{a,b}, C_{b,a}\}$.
\end{proof}

\begin{lem}\label{lem: appendix_optval_equalS}
Let $a,b>0$, $x,p\in\R$, and $t>0$. Let $\valuefn$ be the function defined in~\eqref{eqt: result_S1_1d} and~\eqref{eqt: result_S1_1d_negative} and $\potentfn$ be the function defined in~\eqref{eqt: result_defV} with parameters $a,b$.
Let $[0,t]\ni s\mapsto \opttraj(s;x,t,p ,a,b)\in\R$ be the trajectory defined in~\eqref{eqt: optctrl_defx_1},~\eqref{eqt: optctrl_defx_2},~\eqref{eqt: optctrl_defx_3},~\eqref{eqt: optctrl_defx_4},~\eqref{eqt: optctrl_defx_5}, and~\eqref{eqt: optctrl_defx_neg} for different cases. Then, we have
\begin{equation}\label{eqt: lemA1_optval1_equalS}
    \int_0^t \left(\frac{1}{2} \left(\frac{d}{ds}\opttraj(s;x,t,p ,a,b)\right)^2 - \potentfn(\opttraj(s;x,t,p ,a,b))\right) ds + p \opttraj(0;x,t,p ,a,b) = \valuefn(x,t; p ,a,b).
\end{equation}
\end{lem}
\begin{proof}
If $(x,t,p )\in\Omega_1$, the left-hand side of~\eqref{eqt: lemA1_optval1_equalS} equals
\begin{equation*}
\begin{split}
&\int_0^t \left(\frac{1}{2} \left(p  + as\right)^2 + a\left(x-p (t-s)-\frac{a}{2}(t^2-s^2)\right)\right) ds
+ p  \left(x - p t-\frac{at^2}{2}\right)\\
=\, & \int_0^t \left(a^2s^2 + 2ap s + \frac{p ^2}{2} + ax-ap t - \frac{a^2t^2}{2} \right)ds + p x - p ^2t - \frac{ap t^2}{2}\\
=\, & \frac{a^2t^3}{3} + ap t^2 + \frac{p^2t}{2} + atx - apt^2 - \frac{a^2t^3}{2} + px - p^2t - \frac{apt^2}{2}\\
=\, & -\frac{a^2t^3}{6} - \frac{apt^2}{2}  + atx - \frac{p^2 t}{2} + px\\
=\, & \valuefn(x,t;p,a,b).
\end{split}
\end{equation*}
If $(x,t,p)\in\Omega_2$, the left-hand side of~\eqref{eqt: lemA1_optval1_equalS} equals
\begin{equation*}
\begin{split}
    &\int_0^t \left(\frac{1}{2} \left(p-bs\right)^2 - b\left(x - p(t-s) + \frac{b}{2}(t^2-s^2)\right)\right) ds + p\left(x-pt + \frac{bt^2}{2}\right)\\
    =\, & \int_0^t\left(b^2s^2 -2bps + \frac{p^2}{2} - bx +bpt - \frac{b^2t^2}{2} \right)ds + px - p^2t + \frac{bpt^2}{2}\\
    =\, & \frac{b^2t^3}{3} - bpt^2 +\frac{p^2t}{2}- bxt + bpt^2 - \frac{b^2t^3}{2} + px - p^2t + \frac{bpt^2}{2}\\
    =\, & -\frac{b^2t^3}{6} + \frac{bpt^2}{2} - \frac{p^2 t}{2} + px - btx\\
    =\, & \valuefn(x,t;p,a,b).
\end{split}
\end{equation*}
If $(x,t,p)\in\Omega_3$, let $\Delta\in\R$ be defined by $\Delta:= (bt-p)^2 + 2(2b+a)x$ and  $\tau$ be the scalar defined in~\eqref{eqt:def_opttraj_case3_tau}. Then, the left-hand side of~\eqref{eqt: lemA1_optval1_equalS} equals
\begin{equation*}
\begin{split}
    &\int_0^\tau \left(\frac{1}{2} \left(p-bs\right)^2 - b\left(-p(\tau - s) + \frac{b}{2}(\tau^2 - s^2)
    \right)\right) ds + p\left(-p\tau + \frac{b\tau^2}{2}\right) \\
    &\quad\quad  +\int_\tau^t \left(\frac{1}{2} \left(p-b\tau + a(s - \tau)\right)^2 + a\left((p - b\tau)(s-\tau) + \frac{a}{2}(s-\tau)^2
    \right) \right)ds \\
    =\, & \int_0^\tau\left( b^2s^2 -2bps + \frac{p^2}{2}+bp\tau - \frac{b^2\tau^2}{2} \right)ds - p^2\tau + \frac{bp\tau^2}{2}\\ &\quad\quad +\int_0^{t-\tau} \left(\frac{1}{2} \left(p-b\tau + as\right)^2 + a\left((p - b\tau)s + \frac{a}{2}s^2
    \right) \right)ds \\
    =\, & \frac{b^2\tau^3}{3} - bp\tau^2 + \frac{p^2\tau}{2} + bp\tau^2 - \frac{b^2\tau^3}{2} - p^2\tau +\frac{bp\tau^2}{2} \\
    &\quad\quad +\int_0^{t-\tau} \left(a^2s^2 + 2a(p-b\tau)s + \frac{(p-b\tau)^2}{2}\right)ds\\
    =\, & -\frac{b^2\tau^3}{6} + \frac{bp\tau^2}{2} - \frac{p^2\tau}{2} + \frac{a^2(t-\tau)^3}{3} + a(p-b\tau)(t-\tau)^2 + \frac{(p-b\tau)^2(t-\tau)}{2}.
\end{split}
\end{equation*}
Let $r := t-\tau = \frac{bt-p + \sqrt{\Delta}}{2b+a}$ and $A:= bt-p$. After some calculations, the left-hand side of~\eqref{eqt: lemA1_optval1_equalS} equals
\begin{equation*}
\begin{split}
&\frac{(a+b)(a+2b)r^3}{3} + (p-bt)\frac{(2a+3b)r^2}{2} + (bt - p)^2r - \frac{b^2t^3}{6} + \frac{bpt^2}{2} - \frac{p^2t}{2}\\
=\,& \frac{a+b}{3(a+2b)^2}\left(4A^3 + 3A^2\sqrt{\Delta} + 6(a+2b)xA+\sqrt{\Delta}^3\right)\\
&\quad \quad \quad \quad - \frac{2a+3b}{(a+2b)^2}\left(A^3 +A^2\sqrt{\Delta} +(a+2b)xA\right) \\
&\quad \quad \quad \quad + \frac{A^3+A^2\sqrt{\Delta}}{a+2b}- \frac{b^2t^3}{6} + \frac{bpt^2}{2} - \frac{p^2t}{2}
\\
=\,& \frac{a+b}{3(a+2b)^2}\left(A^3 + \sqrt{\Delta}^3\right) - \frac{bxA}{a+2b}- \frac{b^2t^3}{6} + \frac{bpt^2}{2} - \frac{p^2t}{2}\\
    =\, & \valuefn(x,t;p,a,b).
\end{split}
\end{equation*}
If $(x,t,p)\in\Omega_4$, the left-hand side of~\eqref{eqt: lemA1_optval1_equalS} equals
\begin{equation*}
\begin{split}
&\int_{t-\sqrt{\frac{2x}{a}}}^t \left(\frac{1}{2} \left(a\left(s-t-\sqrt{\frac{2x}{a}}\right)\right)^2 + a\left(\frac{a}{2}\left(s-t+\sqrt{\frac{2x}{a}}\right)^2\right)\right)ds \\
& \quad\quad\quad \quad +\int_0^{\frac{p}{b}} \left(\frac{1}{2} \left(p-bs\right)^2 + b\left(  \frac{1}{2b}(p - bs)^2\right)\right)ds + p\left(-\frac{p^2}{2b}\right)\\
=\,&  \int_{t-\sqrt{\frac{2x}{a}}}^t  a^2\left(s-t+\sqrt{\frac{2x}{a}}\right)^2ds + \int_0^{\frac{p}{b}} \left(p-bs\right)^2 ds - \frac{p^3}{2b}\\
=\,& \int_0^{\sqrt{\frac{2x}{a}}}  a^2s^2ds + \int_0^{\frac{p}{b}} b^2s^2 ds - \frac{p^3}{2b}\\
=\,& \frac{a^2}{3}\left(\frac{2x}{a}\right)^{3/2} + \frac{b^2}{3}\left(\frac{p}{b}\right)^3 - \frac{p^3}{2b}\\
=\, & \valuefn(x,t;p,a,b).
\end{split}
\end{equation*}
If $(x,t,p)\in\Omega_5$, the left-hand side of~\eqref{eqt: lemA1_optval1_equalS} equals
\begin{equation*}
\begin{split}
& \int_{t-\sqrt{\frac{2|x|}{b}}}^t \left(\frac{1}{2} \left(b\left(s-t-\sqrt{\frac{2|x|}{b}}\right)\right)^2 + b\left(\frac{b}{2}\left(s-t+\sqrt{\frac{2|x|}{b}}\right)^2\right)\right)ds\\
&\quad\quad\quad \quad + \int_0^{\frac{p}{b}} \left(\frac{1}{2} \left(p-bs\right)^2 + b\left(  \frac{1}{2b}(p - bs)^2\right)\right)ds + p\left(-\frac{p^2}{2b}\right)\\
=\,& \int_{t-\sqrt{\frac{2|x|}{b}}}^t  b^2\left(s-t+\sqrt{\frac{2|x|}{b}}\right)^2ds + \int_0^{\frac{p}{b}} \left(p-bs\right)^2 ds - \frac{p^3}{2b}\\
=\,& \int_0^{\sqrt{\frac{2|x|}{b}}}  b^2s^2ds + \int_0^{\frac{p}{b}} b^2s^2 ds - \frac{p^3}{2b}\\
=\,& \frac{b^2}{3}\left(\frac{2|x|}{b}\right)^{3/2} + \frac{b^2}{3}\left(\frac{p}{b}\right)^3 - \frac{p^3}{2b}\\
=\, & \valuefn(x,t;p,a,b).
\end{split}
\end{equation*}
Therefore,~\eqref{eqt: lemA1_optval1_equalS} holds for any $x\in\R$, $p,t\geq 0$.

Now, we consider the case where $p< 0$. By definition, we have $\opttraj(s;x,t,p,a,b) = -\opttraj(s;-x,t,-p,b,a)$ and $\valuefn(x,t;p,a,b)= \valuefn(-x,t;-p,b,a)$. Let $\potentfn_{b,a}$ denote the following function:
\begin{equation}\label{eqt:lem_gamma_def_Vba}
    \potentfn_{b,a}(y) := \potentfn(-y) = \begin{dcases}
    -by & y\geq 0,\\
    ay & y< 0.
    \end{dcases}
\end{equation}
Then, the left-hand side of~\eqref{eqt: lemA1_optval1_equalS} equals
\begin{equation*}
\begin{split}
    &\int_0^t\left( \frac{1}{2} \left(-\frac{d}{ds}\opttraj(s;-x,t,-p,b,a)\right)^2 - \potentfn(-\opttraj(s;-x,t,-p,b,a)) \right)ds - p\opttraj(0;-x,t,-p,b,a)\\
    =\,&\int_0^t \left(\frac{1}{2} \left(\frac{d}{ds}\opttraj(s;y,t,q,b,a)\right)^2 - \potentfn_{b,a}(\opttraj(s;y,t,q,b,a)) \right)ds +q\opttraj(0;y,t,q,b,a)\\
    =\,& \valuefn(y,t;q,b,a) =\valuefn(-x,t;-p,b,a) = \valuefn(x,t;p,a,b),
\end{split}
\end{equation*}
where the first equality holds by~\eqref{eqt:lem_gamma_def_Vba} and the change of variables $y=-x$ and $q=-p$ and the second equality holds since we have already proved above that~\eqref{eqt: lemA1_optval1_equalS} holds in the positive $q$ case.
Therefore,~\eqref{eqt: lemA1_optval1_equalS} holds for all $x,p\in\R$ and $t\geq 0$.
\end{proof}

\begin{lem} \label{lem: appendix_gradpS_gamma0}
Let $a,b>0$, $x,p\in\R$, and $t>0$. Let $\valuefn$ be the function defined in~\eqref{eqt: result_S1_1d} and~\eqref{eqt: result_S1_1d_negative} and $\potentfn$ be the function defined in~\eqref{eqt: result_defV} with parameters $a,b$.
Let $[0,t]\ni s\mapsto \opttraj(s;x,t,p,a,b)\in\R$ be the trajectory defined in~\eqref{eqt: optctrl_defx_1},~\eqref{eqt: optctrl_defx_2},~\eqref{eqt: optctrl_defx_3},~\eqref{eqt: optctrl_defx_4},~\eqref{eqt: optctrl_defx_5}, and~\eqref{eqt: optctrl_defx_neg} for different cases. Then, we have
\begin{equation}\label{eqt: lemA2}
    \frac{\partial \valuefn}{\partial p}(x,t;p,a,b) = \opttraj(0;x,t,p,a,b).
\end{equation}
\end{lem}
\begin{proof}
We prove~\eqref{eqt: lemA2} using~\eqref{eqt: HJ1_p_deriv} and straightforward calculation.
If $(x,t,p)\in \Omega_1$, we have
\begin{equation*}
\frac{\partial \valuefn}{\partial p}(x,t;p,a,b) = -\frac{at^2}{2} - pt + x = \opttraj(0;x,t,p,a,b).
\end{equation*}
If $(x,t,p)\in \Omega_2$, we have
\begin{equation*}
\frac{\partial \valuefn}{\partial p}(x,t;p,a,b) = \frac{bt^2}{2} - pt + x = \opttraj(0;x,t,p,a,b).
\end{equation*}
If $(x,t,p)\in \Omega_3$, let $\Delta:=(bt-p)^2 + 2x(2b+a)$. Then, we have
\begin{equation}\label{eqt: lemA2_case3_eq1}
\begin{split}
\frac{\partial \valuefn}{\partial p}(x,t;p,a,b) &= \frac{a+b}{(a+2b)^2}\left(-(bt-p)^2 + (p-bt)\sqrt{\Delta}\right) 
+ \frac{bx}{a+2b} + \frac{bt^2}{2} - pt\\
&= -\frac{(a+b)(bt-p)^2}{(a+2b)^2} + \frac{bx}{a+2b} + \frac{bt^2}{2} - pt + \frac{(a+b)(p-bt)\sqrt{\Delta}}{(a+2b)^2}\\
&= -\frac{(a+b)p^2 }{(a+2b)^2}
-\frac{(a^2+2ab+2b^2)pt }{(a+2b)^2} + \frac{(a^2+2ab+2b^2)bt^2}{2(a+2b)^2}\\
&\quad\quad\quad\quad + \frac{bx}{a+2b} + \frac{(a+b)(p-bt)\sqrt{\Delta}}{(a+2b)^2}.
\end{split}
\end{equation}
Let $c\in\R$ be defined by $ c:= (a+b)t+p$ and $\tau$ be the scalar defined in~\eqref{eqt:def_opttraj_case3_tau}, which equals $\frac{c - \sqrt{\Delta}}{2b+a}$. Then, by definition, we have
\begin{equation} \label{eqt: lemA2_case3_eq2}
\begin{split}
\opttraj(0;x,t,p,a,b) &= -p\tau + \frac{b\tau^2}{2}
= \frac{p\sqrt{\Delta} - pc}{2b+a} + \frac{bc^2 + b\Delta - 2bc\sqrt{\Delta}}{2(2b+a)^2}\\
&= \frac{bc^2+b\Delta - 2pc(2b+a)}{2(2b+a)^2} + \frac{2(2b+a)p-2bc}{2(2b+a)^2}\sqrt{\Delta}.
\end{split}
\end{equation}
By some calculations, we have
\begin{equation}\label{eqt: lemA2_case3_eq3}
\begin{split}
    &bc^2+b\Delta - 2pc(2b+a)\\
    =\,& b\left((a+b)^2t^2 + p^2 + 2(a+b)tp + (bt-p)^2 +2x(2b+a)\right) \\
    &\quad\quad\quad\quad - 2(2b+a)(a+b)tp - 2(2b+a)p^2\\
    =\,& b(bt-p)^2 -(3b+2a)p^2 - 2(a+b)^2tp + b(a+b)^2t^2 + 2bx(2b+a)\\
    =\,& -2(a+b)p^2 -2((a+b)^2+b^2)tp + ((a+b)^2+b^2)bt^2 + 2(2b+a)bx,
\end{split}
\end{equation}
and 
\begin{equation}\label{eqt: lemA2_case3_eq4}
    2(2b+a)p - 2bc = 2(2b+a)p - 2b(a+b)t - 2bp = 2(a+b)p - 2(a+b)bt.
\end{equation}
Combining~\eqref{eqt: lemA2_case3_eq1}, \eqref{eqt: lemA2_case3_eq2}, \eqref{eqt: lemA2_case3_eq3}, and~\eqref{eqt: lemA2_case3_eq4}, we obtain~\eqref{eqt: lemA2} for $(x,t,p)\in\Omega_3$.

If $(x,t,p)\in \Omega_4\cup \Omega_5$, we have
\begin{equation*}
\frac{\partial \valuefn}{\partial p}(x,t;p,a,b) = -\frac{p^2}{2b}= \opttraj(0;x,t,p,a,b).
\end{equation*}
Thus, we have proved~\eqref{eqt: lemA2} for any $x\in\R$ and $t,p\geq 0$.

If $p <0$, let $q=-p$. Then, we have 
\begin{equation*}
\frac{\partial \valuefn}{\partial p}(x,t;p,a,b)
= -\frac{\partial \valuefn}{\partial q}(-x,t;q,b,a)
= -\opttraj(0;-x,t,q,b,a) = \opttraj(0;x,t,p,a,b),
\end{equation*}
where the first equality holds by~\eqref{eqt: result_S1_1d_negative}, the second equality holds since we have already proved~\eqref{eqt: lemA2} in the positive $q$ case, and the third equality holds by~\eqref{eqt: optctrl_defx_neg}.
Hence, \eqref{eqt: lemA2} also holds for any $x\in\R$, $t\geq 0$, and $p<0$.
\end{proof}

\begin{lem} \label{lem: appendix_lowerbd_ell}
Let $\initcond\colon \Rn\to\R$ be a convex function. Let $\{a_i,b_i\}_{i=1}^n$ be positive constants and $\potentfn_i\colon \R\to\R$ be the function defined by~\eqref{eqt: result_defV} with constants $a = a_i$ and $b = b_i$ for each $i\in\{1,\dots, n\}$. Define the function $\ell\colon \Rn\times [0,+\infty) \times \Rn\to \R\cup\{+\infty\}$ by
\begin{equation} \label{eqt:lemA4_def_ell}
    \ell(\bx,t,\bp):= -\sum_{i=1}^n \valuefn(x_i, t; p_i, a_i,b_i) + \initcond^*(\bp),
\end{equation}
for all $\bx=(x_1,\cdots,x_n)\in\Rn$, $\bp=(p_1,\dots, p_n)\in\Rn$, and $t\geq 0$, where each function $\valuefn(x_i, t; p_i, a_i,b_i)$ on the right-hand side is the function defined in~\eqref{eqt: result_S1_1d} and~\eqref{eqt: result_S1_1d_negative}.
Let $M > 0$ and $\alpha \in\R$ be arbitrary scalars. Then, there exists $M_{\bp}>0$, such that
$
    \ell(\bx,t,\bp)\geq \alpha
$
holds for all $t\in [0,M]$ and for all $\bx, \bp\in\Rn$ satisfying $\|\bx\|\leq M$ and $\|\bp\|\geq M_{\bp}$.
\end{lem}
\begin{proof}
Since $\initcond$ is a finite-valued convex function, its Legendre-Fenchel transform $\initcond^*$ is 1-coercive, and hence, the function $\bp\mapsto \initcond^*(\bp) - M\|\bp\|$ is also 1-coercive. As a result, there exists $M_{\bp}>0$, such that 
\begin{equation}\label{eqt:lemA9_pf_1coercive}
\initcond^*(\bp) - M\|\bp\| \geq M^2\sum_{i=1}^n(a_i+b_i) + \alpha,
\end{equation}
for all $\bp\in\Rn$ satisfying $\|\bp\|\geq M_{\bp}$. By straightforward calculation, for all $t\in[0,M]$ and for all $\bp,\bx\in\Rn$ satisfying $\|\bp\|\geq M_{\bp}$ and $\|\bx\|\leq M$, we have
\begin{equation*}
\begin{split}
    \ell(\bx,t,\bp) &= -\sum_{i=1}^n \valuefn(x_i, t; p_i, a_i,b_i) + \initcond^*(\bp)\geq -\sum_{i=1}^n(p_ix_i + t(a_i+b_i)|x_i|) + \initcond^*(\bp)\\
    &\geq -\|\bp\|\|\bx\| - t\|\bx\|\sum_{i=1}^n(a_i+b_i) + \initcond^*(\bp)\\
    &\geq -M\|\bp\| - M^2\sum_{i=1}^n(a_i+b_i) + \initcond^*(\bp)\\
    &\geq \alpha,
\end{split}
\end{equation*}
where the first inequality holds by Lemma~\ref{lem: appendix_S_bd} and the last inequality follows from~\eqref{eqt:lemA9_pf_1coercive}.
\end{proof}

\begin{lem} \label{lem: continuity_SJ_pde1_hd}
Let $\initcond\colon \Rn\to\R$ be a convex function. Let $\{a_i,b_i\}_{i=1}^n$ be positive constants and $\potentfn_i\colon \R\to\R$ be the function defined by~\eqref{eqt: result_defV} with constants $a = a_i$ and $b = b_i$ for each $i\in\{1,\dots, n\}$. Let $\valuefn\colon \R^n\times [0,+\infty)\to\R$ be the function defined in~\eqref{eqt: result_Hopf1_hd}. Then, the function $\valuefn$ is continuous in $\Rn\times [0,+\infty)$.
\end{lem}
\begin{proof}
First, we show the existence of the maximizer in~\eqref{eqt: result_Hopf1_hd}.
Let $\ell\colon \Rn\times [0,+\infty)\times \Rn\to\R\cup\{+\infty\}$ be the function defined in~\eqref{eqt:lemA4_def_ell}.
The function $\bp\mapsto -\ell(\bx,t,\bp)$ is the objective function in the maximization problem~\eqref{eqt: result_Hopf1_hd} at $(\bx,t)$.
Since $\initcond$ is finite-valued and convex, its Legendre-Fenchel transform $\initcond^*$ is convex, lower semi-continuous, and 1-coercive. By Lemma~\ref{lem: prop_S}, if $t>0$, the function $p_i\mapsto -\valuefn(x_i, t; p_i,a_i,b_i)$ is convex for each $i\in\{1,\dots,n\}$. If $t=0$, by straightforward calculation, we have
$
    -\valuefn(x_i, t; p_i,a_i,b_i) = -p_ix_i,
$
which is a convex function with respect to $p_i$. Therefore, for all $\bx\in\Rn$ and $t\geq 0$, the function $\ell$ is convex, lower semi-continuous, and 1-coercive with respect to $\bp$. As a result, the maximizer in~\eqref{eqt: result_Hopf1_hd} exists for all $\bx\in\Rn$ and $t\geq 0$ (see~\cite[Definition~IV.3.2.6]{Hiriart1993Convex}), and hence, the function $\valuefn$ is finite-valued in $\Rn\times [0,+\infty)$.

Let $\bx\in\Rn$ and $t\geq 0$. Now, we show the continuity of the function $\valuefn$ at $(\bx,t)$ by showing it is lower semi-continuous and upper semi-continuous. We begin by proving lower semi-continuity. Let $\bp^*$ be a maximizer in~\eqref{eqt: result_Hopf1_hd} at $(\bx,t)$. 
By Lemma~\ref{prop: S_linearJ}, each function $(x_i,t)\mapsto \valuefn(x_i,t;p_i,a_i,b_i)$ is continuous, and hence, the function $(\bx,t)\mapsto\ell(\bx,t,\bp)$ is continuous for all $\bp$ in the domain of $\initcond^*$.
For any sequence $\{(\bx^k,t^k)\}\subseteq \Rn\times [0,+\infty)$ converging to $(\bx,t)$, we have
\begin{equation*}
    \liminf_{k\to\infty} \valuefn(\bx^k,t^k) \geq \liminf_{k\to\infty} -\ell(\bx^k,t^k,\bp^*) = -\ell(\bx,t,\bp^*) = \valuefn(\bx,t),
\end{equation*}
where the inequality holds by definition of $\valuefn$ in~\eqref{eqt: result_Hopf1_hd}, the first equality holds since the function $(\bx,t)\mapsto\ell(\bx,t,\bp^*)$ is continuous, and the last equality holds since $\bp^*$ is a maximizer in~\eqref{eqt: result_Hopf1_hd} at $(\bx,t)$. Therefore, the function $\valuefn$ is lower semi-continuous at $(\bx,t)$. 

Now, we prove that $\valuefn$ is upper semi-continuous at $(\bx,t)$. Let $\delta>0$ be an arbitrary positive number. Our goal is to find a neighborhood $\mathcal{N}_{\bx,t}$ of $(\bx,t)$ in $\Rn\times [0,+\infty)$, such that there holds
\begin{equation}\label{eqt:lemA2_pf_upper_ineqt1}
\valuefn(\by,s) \leq \valuefn(\bx,t)+\delta    \quad \forall (\by,s)\in\mathcal{N}_{\bx,t}.
\end{equation}
Note that by definition of $\valuefn$ in~\eqref{eqt: result_Hopf1_hd}, the inequality in~\eqref{eqt:lemA2_pf_upper_ineqt1} holds if and only if there holds
\begin{equation}\label{eqt:lemA2_pf_upper_ineqt2}
    -\ell(\by,s,\bp)\leq \valuefn(\bx,t)+\delta \quad \forall \bp\in\Rn.
\end{equation}
Therefore, it suffices to find a neighborhood $\mathcal{N}_{\bx,t}$ of $(\bx,t)$ in $\Rn\times [0,+\infty)$, such that~\eqref{eqt:lemA2_pf_upper_ineqt2} holds for all $(\by,s)\in\mathcal{N}_{\bx,t}$ and $\bp\in\Rn$.
Let $M>0$ be a scalar. 
By Lemma~\ref{lem: appendix_lowerbd_ell} with $\alpha = -\valuefn(\bx,t)-\delta$, there exists $M_{\bp} >0$, such that for all $\by,\bp\in\Rn$ and $s\geq 0$ satisfying $\|\by\|\leq M$, $s\leq M$, and $\|\bp\|\geq M_{\bp}$, we have 
\begin{equation}\label{eqt: lem22_eq1}
\begin{split}
    \ell(\by,s,\bp) \geq \alpha = -\valuefn(\bx,t) - \delta.
\end{split}
\end{equation}
Now, we consider the case when $\|\bp\| < M_{\bp}$ holds. By Lemmas~\ref{lem: prop_S} and~\ref{prop: S_linearJ}, the function $(y,s,p)\mapsto \valuefn(y,s;p,a_i,b_i)$ is continuous in $\R\times [0,+\infty)\times \R$ for each $i\in\{1,\dots,n\}$. Hence, there exists a neighborhood $\tilde{\mathcal{N}}_{\bx,t}\subset \Rn\times [0,+\infty)$ of $(\bx,t)$, such that $|\valuefn(y_i,s;p_i,a_i,b_i) - \valuefn(x_i,t;p_i,a_i,b_i)|\leq \frac{\delta}{n}$ holds for all $(\by,s)\in \tilde{\mathcal{N}}_{\bx,t}$, $\|\bp\|\leq M_{\bp}$ and $i \in\{1,\dots,n\}$. Recall that we use $x_i$, $y_i$, and $p_i$ to denote the $i$-th component of the vectors $\bx$, $\by$, and $\bp$, respectively. Therefore, for all $(\by,s)\in \tilde{\mathcal{N}}_{\bx,t}$ and $\|\bp\|\leq M_{\bp}$, we have
\begin{equation}\label{eqt: lem22_eq2}
\begin{split}
    \ell(\by,s,\bp) &= -\sum_{i=1}^n \valuefn(y_i, s; p_i, a_i,b_i) + \initcond^*(\bp)\\
    &\geq
    -\sum_{i=1}^n \left(\valuefn(x_i, t; p_i, a_i,b_i) + \frac{\delta}{n}\right) + \initcond^*(\bp)\\
    & = -\sum_{i=1}^n \valuefn(x_i, t; p_i, a_i,b_i) + \initcond^*(\bp) - \delta \\
    & = \ell(\bx,t,\bp)-\delta \geq -\valuefn(\bx,t) - \delta.
\end{split}
\end{equation}
Combining~\eqref{eqt: lem22_eq1} and~\eqref{eqt: lem22_eq2}, we conclude that~\eqref{eqt:lemA2_pf_upper_ineqt2} holds for all $\bp\in\Rn$ and $(\by,s)\in \tilde{\mathcal{N}}_{\bx,t}\cap \left(B_M(\R^n)\times [0,M]\right)$.
Therefore, the function $\valuefn$ is upper semi-continuous at $(\bx,t)$.

Since $(\bx,t)$ is an arbitrary point in $\Rn\times [0,+\infty)$, we conclude that
the function $\valuefn$ is continuous in $\Rn\times [0,+\infty)$.
\end{proof}

\begin{lem} \label{lem: prop_barp_hopf_hd}
Let $\initcond\colon \Rn\to\R$ be a convex function. Let $\{a_i,b_i\}_{i=1}^n$ be positive constants and $\potentfn_i\colon \R\to\R$ be the function defined by~\eqref{eqt: result_defV} with constants $a=a_i$ and $b=b_i$ for each $i\in\{1,\dots, n\}$. 
Let $\bp^*(\bx,t)$ be the set of maximizers in~\eqref{eqt: result_Hopf1_hd}. When the maximizer is unique, we abuse notation and also use $\bp^*(\bx,t)$ to denote the unique maximizer (as opposed to the singleton containing the maximizer), whenever there is no ambiguity.
Then, the following statements hold:
\begin{itemize}
    \item[(a)] For any $\bx\in\Rn$ and $t>0$, the maximizer in~\eqref{eqt: result_Hopf1_hd} exists and is unique.
    \item[(b)] For any $\bx\in\Rn$ and $t=0$, the maximizer in~\eqref{eqt: result_Hopf1_hd} exists, and we have $\bp^*(\bx,0) = \partial \initcond(\bx)$. 
    \item[(c)] For any $\bx\in\Rn$, $t\geq 0$, and any neighborhood $\mathcal{N}_{\bp}$ of $\bp^*(\bx,t)$, there exists a neighborhood $\mathcal{N}_{\bx,t}$ of $(\bx,t)$ in $\Rn\times[0,+\infty)$, such that $\bp^*(\by,s)\subseteq \mathcal{N}_{\bp}$ holds for all $(\by,s)\in\mathcal{N}_{\bx,t}$.
    \item[(d)] 
    The function $\bp^*\colon \Rn\times (0,+\infty) \to \Rn$ is continuous.
    \item[(e)] For any $R>0$, the set $\left\{\bp\colon \bx\in B_R(\Rn),\ t\in [0,R],\ \bp\in \bp^*(\bx,t)\right\}$ is bounded, where $B_R(\Rn)$ denotes the closed ball in $\Rn$ centered at $\mathbf{0}$ with radius $R$. 
\end{itemize}
\end{lem}
\begin{proof}
In the proof of Lemma~\ref{lem: continuity_SJ_pde1_hd}, we have proved the existence of the maximizer in~\eqref{eqt: result_Hopf1_hd} for all $\bx\in\Rn$ and $t\geq 0$. Hence, the set $\bp^*(\bx,t)$ is non-empty for all $\bx\in\Rn$ and $t\geq 0$. Let $\ell$ be the function defined in~\eqref{eqt:lemA4_def_ell}. Recall that the function $\ell$ is convex and lower semi-continuous with respect to $\bp$.
Now, we prove the statements as follows:
\begin{itemize}
    \item[(a)] Let $\bx\in\Rn$ and $t>0$. By Lemma~\ref{lem: prop_S}, the function $p_i\mapsto -\valuefn(x_i,t;p_i,a_i,b_i)$ is strictly convex for each $i\in\{1,\dots,n\}$. Therefore, the function $\bp\mapsto \ell(\bx,t,\bp)$ is strictly convex, and hence, the maximizer in~\eqref{eqt: result_Hopf1_hd} is unique.
    We have already proved the existence of the maximizer in the proof of Lemma~\ref{lem: continuity_SJ_pde1_hd}.
    As a result, the set $\bp^*(\bx,t)$ is a singleton.
    
    \item[(b)] Let $\bx\in\Rn$ and $t=0$. We have $\valuefn(x_i,0;p_i,a_i,b_i) = p_ix_i$ for each $i\in\{1,\dots,n\}$ and any $p_i\in\R$. Hence, the maximization problem in~\eqref{eqt: result_Hopf1_hd} becomes
    \begin{equation} \label{eqt: lem23_SJ_initialcond}
        \valuefn(\bx,0) = \sup_{\bp\in\Rn}\left\{\sum_{i=1}^n x_ip_i - \initcond^*(\bp)\right\} = \sup_{\bp\in\Rn}\left\{\langle \bx,\bp\rangle - \initcond^*(\bp)\right\} = \initcond(\bx).
    \end{equation}
    The set of maximizers equals $\partial \initcond(\bx)$, which is a non-empty, convex, and compact set.
    
    \item[(c)]
    We prove this statement by contradiction. Assume it does not hold. Then, there exist $\bx\in\Rn$, $t\geq 0$, an open neighborhood $\mathcal{N}_{\bp}$ of $\bp^*(\bx,t)$, and a sequence $\{(\by^k, s^k)\}$ in $\Rn\times [0,+\infty)$ converging to $(\bx,t)$, such that $\bp^k\not\in \mathcal{N}_{\bp}$ holds for all $k\in\N$, where $\bp^k$ is a maximizer in~\eqref{eqt: result_Hopf1_hd} at $(\by^k,s^k)$. Let $\delta >0$ be any positive scalar. Since the function $\valuefn$ is continuous according to Lemma~\ref{lem: continuity_SJ_pde1_hd}, we assume 
    \begin{equation}\label{eqt:lemA3_cd_pf_Scont}
        |\valuefn(\by^k,s^k) - \valuefn(\bx,t)|< \delta\quad \forall k\in\N
    \end{equation}
    by taking the tail of the sequence $\{(\by^k,s^k)\}$. Note that the sequence $\{(\by^k,s^k)\}$ is bounded. Then, according to Lemma~\ref{lem: appendix_lowerbd_ell}, there exists $M_{\bp} >0$, such that there holds
    \begin{equation} \label{eqt:lemA3_cd_pf_bd_ell}
        \ell(\by^k,s^k,\bp) \geq -\valuefn(\bx,t) + \delta,
    \end{equation}
    for all $\bp\in\Rn$ satisfying $\|\bp\|\geq M_{\bp}$. Since $\bp^k$ is a maximizer in~\eqref{eqt: result_Hopf1_hd} at $(\by^k,s^k)$, we have
    \begin{equation}\label{eqt:lemA3_cd_pf_bd_ell2}
        \ell(\by^k,s^k,\bp^k) = -\valuefn(\by^k,s^k) < -\valuefn(\bx,t) + \delta \quad \forall k\in\N,
    \end{equation}
    where the inequality follows from~\eqref{eqt:lemA3_cd_pf_Scont}.
    Note that~\eqref{eqt:lemA3_cd_pf_bd_ell2} contradicts with~\eqref{eqt:lemA3_cd_pf_bd_ell}, and hence, we get $\|\bp^k\|< M_{\bp}$ for all $k\in\N$. Therefore, the sequence $\{\bp^k\}$ is bounded. By taking a convergent subsequence, we assume $\{\bp^k\}$ converges to a point in $\Rn$ denoted by $\bp^0$. After some calculation, we obtain
    \begin{equation*}
        \valuefn(\bx,t) = \lim_{k\to\infty} \valuefn(\by^k, s^k) = -\lim_{k\to\infty}\ell(\by^k, s^k, \bp^k)
        \leq -\ell(\bx,t,\bp^0),
    \end{equation*}
    where the first equality holds since the function $\valuefn$ is continuous by Lemma~\ref{lem: continuity_SJ_pde1_hd}, the second equality holds since $\bp^k$ is a maximizer in~\eqref{eqt: result_Hopf1_hd}, and the last inequality holds since $\ell$ is lower semi-continuous and $\{(\by^k, s^k, \bp^k)\}$ converges to $(\bx,t,\bp^0)$. Then, by definition of $\valuefn(\bx,t)$, we conclude that $\bp^0$ is a maximizer in~\eqref{eqt: result_Hopf1_hd} at $(\bx, t)$. Hence, we have $\bp^0\in \bp^*(\bx,t) \subseteq \mathcal{N}_{\bp}$. However, since $\bp^0$ is the limit of $\{\bp^k\}$ and each $\bp^k$ is not in the open set $\mathcal{N}_{\bp}$, the limit point $\bp^0$ is also not in the set $\mathcal{N}_{\bp}$, which gives a contradiction.
    
    \item[(d)]
    By (a), the maximizer in~\eqref{eqt: result_Hopf1_hd} at any $(\bx,t)\in\Rn\times (0,+\infty)$ is unique, and we denote the unique maximizer by $\bp^*(\bx,t)$. According to (c), for any $\bx\in\Rn$, $t> 0$, and any neighborhood $\mathcal{N}_{\bp}$ of $\bp^*(\bx,t)$, there exists a neighborhood $\mathcal{N}_{\bx,t}$ of $(\bx,t)$ in $\Rn\times (0,+\infty)$, such that $\bp^*(\by,s)\in \mathcal{N}_{\bp}$ holds for all $(\by,s)\in\mathcal{N}_{\bx,t}$. This proves the continuity of the function $\bp^*\colon\Rn\times (0,+\infty)\to\Rn$.
    
    \item[(e)] We first prove the boundedness of $\bp^*(\bx,t)$ for all $\bx\in\Rn$ and $t\geq 0$. If $t>0$, the set $\bp^*(\bx,t)$ is a singleton by (a), and hence, it is a bounded set. If $t=0$, the set $\bp^*(\bx,t)$ equals $\partial \initcond(\bx)$, which is a bounded set since $\initcond$ is a finite-valued convex function. 
    Now, we prove (e) by contradiction. Assume~(e) does not hold. Then, there exist sequences $\{\bx^k\}\subset B_R(\Rn)$, $\{t^k\}\subset[0,R]$, and $\{\bp^k\}\subset \Rn$, such that $\bp^k$ is in $\bp^*(\bx^k,t^k)$ for all $k\in\N$ and $\|\bp^k\|$ increases to infinity as $k$ goes to infinity. Since $\{(\bx^k,t^k)\}$ is a bounded sequence, by replacing it with a convergent subsequence, we can assume that $\{(\bx^k,t^k)\}$ converges to a point denoted by $(\bx^*,t^*)\in\Rn\times [0,+\infty)$.
    Let $\mathcal{N}_{\bp}$ be a bounded neighborhood of $\bp^*(\bx^*,t^*)$, which exists since we have proved the boundedness of the set $\bp^*(\bx^*,t^*)$. According to (c), there exists a neighborhood $\mathcal{N}_{\bx^*,t^*}$ of $(\bx^*,t^*)$ such that $\bp^*(\by,s)\subseteq \mathcal{N}_{\bp}$ holds for all $(\by,s)\in\mathcal{N}_{\bx^*,t^*}$. Since the sequence $\{(\bx^k,t^k)\}$ converges to $(\bx^*,t^*)$, by replacing it with a tail sequence, we can assume that $\{(\bx^k,t^k)\}$ is in the neighborhood $\mathcal{N}_{\bx^*,t^*}$.
    Hence, the set $\bp^*(\bx^k,t^k)$ is included in the bounded set $\mathcal{N}_{\bp}$ for all $k\in\N$, which contradicts the assumption that $\bp^k$ is in $\bp^*(\bx^k,t^k)$ for all $k\in\N$ and $\|\bp^k\|$ increases to infinity as $k$ goes to infinity. 
    Therefore, statement~(e) holds.
\end{itemize}
\end{proof}

\begin{lem}\label{lem: lem24_grad_SJ}
Let $\initcond\colon \Rn\to\R$ be a convex function. Let $\{a_i,b_i\}_{i=1}^n$ be positive constants and $\potentfn_i\colon \R\to\R$ be the function defined by~\eqref{eqt: result_defV} with constants $a = a_i$ and $b = b_i$ for each $i\in\{1,\dots, n\}$. Let $\valuefn\colon \R^n\times [0,+\infty)\to\R$ be the function defined in~\eqref{eqt: result_Hopf1_hd}. Then, the function $\valuefn$ is continuously differentiable in $\Rn\times (0,+\infty)$, and its gradient at $(\bx,t)\in\Rn\times (0,+\infty)$ equals
\begin{equation} \label{eqt: lem24_grad_SJ}
\begin{split}
    \nabla \valuefn(\bx,t) = \Bigg(\frac{\partial \valuefn(x_1,t;p^*_1(\bx,t),a_1,b_1)}{\partial x},\dots, \frac{\partial \valuefn(x_n,t;p^*_n(\bx,t), a_n,b_n)}{\partial x},\quad \quad \\
    \sum_{i=1}^n\frac{\partial \valuefn(x_i,t;p^*_i(\bx,t), a_i,b_i)}{\partial t}\Bigg),
\end{split}
\end{equation}
where $\bp^*(\bx,t) = (p^*_1(\bx,t), \dots, p^*_n(\bx,t))$ denotes the unique maximizer in~\eqref{eqt: result_Hopf1_hd} at $(\bx,t)$ and the functions $\frac{\partial \valuefn(x_i,t;p^*_i(\bx,t),a_i,b_i)}{\partial x}$ and $\frac{\partial \valuefn(x_i,t;p^*_i(\bx,t),a_i,b_i)}{\partial t}$ on the right-hand side denote the derivatives of the function defined in~\eqref{eqt: result_S1_1d} and~\eqref{eqt: result_S1_1d_negative} with respect to $x$ and $t$, respectively.
\end{lem}
\begin{proof}
In this proof, whenever there is no ambiguity, we write $\bp^* = (p^*_1,\dots, p^*_n)$ instead of $\bp^*(\bx,t) = (p^*_1(\bx,t), \dots, p^*_n(\bx,t))$, and we write $\valuefn_i(x,t;p)$ instead of $\valuefn(x,t;p,a_i,b_i)$, for simplicity.
Let $\bx\in\Rn$ and $t>0$. 
Let $\bp^* = (p^*_1,\dots, p^*_n)$ be the maximizer in~\eqref{eqt: result_Hopf1_hd} at $(\bx,t)$. By Lemma~\ref{lem: prop_barp_hopf_hd}(a), the maximizer $\bp^*$ exists and is unique.

We first compute the directional derivative of $\valuefn$ at $(\bx,t)$. Consider the spatial direction $\by\in\Rn$ and time direction $s\in\R$. Let $(\by^k, s^k)\in\Rn\times (0,+\infty)$ be the perturbed vector around $(\bx,t)$ along the direction $(\by,s)$. To be specific, define $\by^k := \bx + \alpha^k \by$ and $s^k := t + \alpha^k s$ where $\{\alpha^k\}$ is a sequence of positive numbers converging to zero. Let $\bp^k\in\Rn$ be the unique maximizer in~\eqref{eqt: result_Hopf1_hd} at $(\by^k,s^k)$. 
Denote the $i$-th component of $\by^k$ and $\bp^k$ by $\yki$ and $\pki$, respectively.
By definition of $\valuefn$, we have
\begin{equation*}
\begin{split}
    \valuefn(\by^k,s^k) &= \sum_{i=1}^n \valuefn_i(\yki, s^k; \pki) - \initcond^*(\bp^k)\geq \sum_{i=1}^n \valuefn_i(\yki, s^k; p^*_i) - \initcond^*(\bp^*),\\
    \valuefn(\bx,t) &= \sum_{i=1}^n \valuefn_i(x_{i}, t; p^*_{i}) - \initcond^*(\bp^*)\geq \sum_{i=1}^n \valuefn_i(x_{i}, t; \pki) - \initcond^*(\bp^k),
\end{split}
\end{equation*}
for all $k\in\N$. On the one hand, we have
\begin{equation*}
\begin{split}
    &\liminf_{k\to\infty} \frac{\valuefn(\by^k,s^k) - \valuefn(\bx,t)}{\alpha^k} \\
    \geq& \liminf_{k\to\infty} \frac{\sum_{i=1}^n \valuefn_i(\yki, s^k; p^*_i) - \initcond^*(\bp^*) - \left(\sum_{i=1}^n \valuefn_i(x_{i}, t; p^*_{i}) - \initcond^*(\bp^*)\right)}{\alpha^k}\\
    =& \sum_{i=1}^n \liminf_{k\to\infty} \frac{\valuefn_i(\yki, s^k; p^*_i) - \valuefn_i(x_{i}, t; p^*_{i})}{\alpha^k}\\
    =& \sum_{i=1}^n \langle \nabla \valuefn_i(x_i, t; p^*_i), (y_i, s)\rangle,
\end{split}
\end{equation*}
where $\nabla \valuefn_i(x_i, t; p^*_i)$ denotes the gradient of the function $(x,t)\mapsto \valuefn(x, t; p^*_i,a_i,b_i)$ at $(x_i,t)$ and $y_i$ denotes the $i$-th component of $\by$. On the other hand, we have
\begin{equation*}
    \begin{split}
    &\limsup_{k\to\infty} \frac{\valuefn(\by^k,s^k) - \valuefn(\bx,t)}{\alpha^k} \\
    \leq& \limsup_{k\to\infty} \frac{\sum_{i=1}^n \valuefn_i(\yki, s^k; \pki) - \initcond^*(\bp^k) - \left(\sum_{i=1}^n \valuefn_i(x_{i}, t; \pki) - \initcond^*(\bp^k)\right)}{\alpha^k}\\
    =& \sum_{i=1}^n \limsup_{k\to\infty} \frac{\valuefn_i(\yki, s^k; \pki) - \valuefn_i(x_{i}, t; \pki)}{\alpha^k}\\
    =& \sum_{i=1}^n \limsup_{k\to\infty} \langle \nabla \valuefn_i(x_i + \xiki y_i, t + \xiki s; \pki), (y_i, s)\rangle,\\
    =& \sum_{i=1}^n \langle \nabla \valuefn_i(x_i, t; p^*_i), (y_i, s)\rangle,
\end{split}
\end{equation*}
where the second equality holds for some constants $\xiki\in [0, \alpha^k]$ for each $i\in\{1,\dots, n\}$ and $k\in\N$ by Taylor's theorem and the last equality holds since we have $\lim_{k\to\infty} \xiki = 0$ for each $i\in\{1,\dots,n\}$, $\lim_{k\to\infty} \bp^k = \bp^*$ by Lemma~\ref{lem: prop_barp_hopf_hd}(d), and each function $\valuefn_i$ is continuously differentiable with respect to $(x,t,p)$ by Lemmas~\ref{lem: prop_S} and~\ref{prop: S_linearJ}. Therefore, we conclude that 
\begin{equation*}
\begin{split}
&\lim_{k\to\infty} \frac{\valuefn(\by^k,s^k) - \valuefn(\bx,t)}{\alpha^k}\\
=& \sum_{i=1}^n \langle \nabla \valuefn_i(x_i, t; p^*_i), (y_i, s)\rangle\\
=& \left\langle \left(\frac{\partial \valuefn_1}{\partial x}(x_1,t;p^*_1), \dots, \frac{\partial \valuefn_n}{\partial x}(x_n,t;p^*_n), \sum_{i=1}^n\frac{\partial \valuefn_i}{\partial t}(x_i,t;p^*_i)\right), (\by, s)\right\rangle.
\end{split}
\end{equation*}
This equality holds for any direction $(\by, s)\in\Rn\times (0,+\infty)$. As a result, the function $\valuefn$ is differentiable at $(\bx,t)$, and the gradient satisfies~\eqref{eqt: lem24_grad_SJ}.
Note that $(\bx,t)$ is an arbitrary point in $\Rn\times (0,+\infty)$, and hence, the function $\valuefn$ is differentiable in $\Rn\times (0,+\infty)$ with gradient equal to~\eqref{eqt: lem24_grad_SJ}.

It remains to prove the continuity of the gradient of $\valuefn$. 
By Lemmas~\ref{lem: prop_S} and~\ref{prop: S_linearJ}, each function $(x,t,p)\mapsto \valuefn(x,t;p,a_i,b_i)$ on the right-hand side of~\eqref{eqt: lem24_grad_SJ} is continuously differentiable in $\R\times (0,+\infty)\times \R$. Moreover, the function $\bp^*$ is also continuous by Lemma~\ref{lem: prop_barp_hopf_hd}(d). Therefore, the right-hand side of~\eqref{eqt: lem24_grad_SJ} is continuous with respect to $(\bx,t)\in\Rn\times (0,+\infty)$. 
As a result, the function $\valuefn$ is continuously differentiable with respect to $(\bx,t)$ in $\Rn\times (0,+\infty)$, and the gradient satisfies~\eqref{eqt: lem24_grad_SJ}.
\end{proof}

\begin{lem}\label{lem:A10_uniqueness}
Assume $\initcond\colon\Rn\to\R$ is a continuous function satisfying~\eqref{eqt: HJ1_hd_condJ} and is bounded from below by an affine function. Let $\potentfn\colon\Rn\to(-\infty,0]$ be a Lipschitz continuous function and $M$ be a symmetric positive definite matrix with $n$ rows and $n$ columns. Let $\valuefn\colon\Rn\times [0,+\infty)\to \R\cup\{-\infty\}$ be the value function defined by~\eqref{eqt: result_optctrl1_hd_general}.
Then, the following statements holds:
\begin{itemize}
    \item[(a)] The function $\valuefn$ is finite-valued for all $\bx\in\Rn$ and $t\geq 0$, and it is a viscosity solution to the HJ PDE~\eqref{eqt: HJhd_1_general}
    in the solution set $\mathcal{G}$ defined in~\eqref{eqt: HJ1_hd_solset_G}.

    \item[(b)] Assume that for all $\balp\in\Rn$ such that $\bx\mapsto \initcond(\bx)-\langle \balp, \bx\rangle $ is bounded from below, 
    the function $(\bx,t)\mapsto \valuefn(\bx,t) - \langle \balp, \bx\rangle $ is also bounded from below. Then, the function $\valuefn$ is the unique viscosity solution to the HJ PDE~\eqref{eqt: result_optctrl1_hd_general} in the solution set $\mathcal{G}$.
\end{itemize}
\end{lem}

\begin{proof}
(a) Since $\initcond$ is bounded from below by an affine function, there exists a vector $\balp\in\Rn$ and a scalar $\beta\in\R$ satisfying $\initcond(\bx)\geq \langle \balp, \bx\rangle + \beta$ for all $\bx\in\Rn$. Define $\tilde{\initcond}\colon\Rn\to\R$ by
\begin{equation}\label{eqt:lemA10a_pf_def_tildeJ}
    \tilde{\initcond}(\bx):= \initcond(\bx)-\langle\balp, \bx\rangle \quad \forall \bx\in\Rn.
\end{equation}
Then, the function $\tilde{\initcond}$ is bounded from below. 
Now, consider another optimal control problem, which reads:
\begin{equation} \label{eqt:lemA10a_pf_newoptctrl}
\begin{split}
 \tilde{\valuefn}(\bx,t) &= \inf \left\{\int_0^t \ell(\bx(s),\bu(s)) ds + \tilde{\initcond}(\bx(0)) \colon \dot{\bx}(s) = f(\bx(s),\bu(s)) \,\,\forall s\in(0,t), \,\, \bx(t) = \bx\right\},
\end{split}
\end{equation}
where the Lagrangian function $\ell\colon \Rn\times A\to\R$ and the source term $f\colon \Rn\times A\to\R$ are defined by
\begin{equation}\label{eqt:lemA10a_pf_def_ell_f}
    \ell(\bx,\bu) := \frac{1}{2}\|\bu\|_{M^{-1}}^2 - \potentfn(\bx) - \frac{1}{2}\|\balp\|_M^2,\quad\quad f(\bx,\bu) := \bu+M\balp,\quad \forall\,\bx,\bu\in\Rn.
\end{equation}
Here and in the rest of this proof, the set $A$ denotes the domain of the control variable $\bu$, which equals $\Rn$ in our case.
By straightforward calculation, the cost in~\eqref{eqt:lemA10a_pf_newoptctrl} equals
\begin{equation*}
\begin{split}
    &\int_0^t \left(\frac{\|\dot{\bx}(s)-M\balp\|_{M^{-1}}^2}{2} - \potentfn(\bx(s))- \frac{\|\balp\|_M^2}{2}\right) ds + \tilde{\initcond}(\bx(0))\\
    =\ &\int_0^t \left(\frac{\|\dot{\bx}(s)\|_{M^{-1}}^2}{2} - \potentfn(\bx(s))- \langle\balp, \dot{\bx}(s)\rangle\right) ds + \initcond(\bx(0)) - \langle \balp,\bx(0)\rangle\\
    =\ &\int_0^t \left(\frac{\|\dot{\bx}(s)\|_{M^{-1}}^2}{2} - \potentfn(\bx(s))\right) ds - \langle \balp, \bx(t)-\bx(0)\rangle + \initcond(\bx(0)) - \langle \balp,\bx(0)\rangle\\
    =\ &\int_0^t \left(\frac{\|\dot{\bx}(s)\|_{M^{-1}}^2}{2} - \potentfn(\bx(s))\right) ds + \initcond(\bx(0)) - \langle \balp,\bx(t)\rangle.
\end{split}
\end{equation*}
As a result, there holds
\begin{equation*}
\begin{split}
 \tilde{\valuefn}(\bx,t) &= \inf \Big\{\int_0^t \left(\frac{\|\bu(s)\|_{M^{-1}}^2}{2} - \potentfn(\bx) - \frac{\|\balp\|_M^2}{2}\right) ds + \tilde{\initcond}(\bx(0)) \\
 &\quad\quad\quad\quad\colon \dot{\bx}(s) = \bu(s)+M\balp \,\,\forall s\in(0,t), \,\, \bx(t) = \bx\Big\}\\
 &= \inf \left\{\int_0^t \left(\frac{\|\dot{\bx}(s) - M\balp\|_{M^{-1}}^2}{2} - \potentfn(\bx) - \frac{\|\balp\|_M^2}{2}\right) ds + \tilde{\initcond}(\bx(0)) \colon \bx(t) = \bx\right\}\\
 &= \inf \left\{\int_0^t \left(\frac{\|\dot{\bx}(s)\|_{M^{-1}}^2}{2} - \potentfn(\bx(s))\right) ds + \initcond(\bx(0)) - \langle \balp,\bx(t)\rangle \colon \bx(t) = \bx\right\}\\
 &= \inf \left\{\int_0^t \left(\frac{\|\dot{\bx}(s)\|_{M^{-1}}^2}{2} - \potentfn(\bx(s))\right) ds + \initcond(\bx(0)) \colon \bx(t) = \bx\right\} - \langle \balp,\bx\rangle\\
 &= \valuefn(\bx,t) - \langle \balp,\bx\rangle.
\end{split}
\end{equation*}
Therefore, we have
\begin{equation}\label{eqt:lemA10a_pf_eqt_tildeS_S}
    \tilde{\valuefn}(\bx,t) = \valuefn(\bx,t) - \langle \balp,\bx\rangle
    \quad \forall\bx\in\Rn, t\geq 0.
\end{equation}
Now, by applying~\cite[Theorem 3.2]{Bardi1997Bellman}, we will prove that the function $\tilde{\valuefn}$ defined in~\eqref{eqt:lemA10a_pf_newoptctrl} is the unique viscosity solution to the HJ PDE defined by
\begin{equation} \label{eqt: lemA10_newHJ1}
\begin{dcases} 
\frac{\partial \tilde{\valuefn}}{\partial t}(\bx,t)+\frac{1}{2}\|\nabla_{\bx}\tilde{\valuefn}(\bx,t) + \balp\|_M^2 + \potentfn(\bx) = 0 & \bx\in\Rn, t\in(0,+\infty),\\
\tilde{\valuefn}(\bx,0)=\tilde{\initcond}(\bx) & \bx\in\Rn,
\end{dcases}
\end{equation}
in the solution set $\tilde{\mathcal{G}}$ defined by
\begin{equation}\label{eqt:prop23_pf_def_tildeG}
\begin{split}
    \tilde{\mathcal{G}} := \{W \in C(\Rn\times [0, +\infty))\colon W\text{ is bounded from below in }\Rn\times [0,T] \,\forall T>0\\
    \text{ and }\|W\|_R < +\infty \ \forall R > 0\},
\end{split}
\end{equation}
where $\|\cdot\|_R$ is defined in~\eqref{eqt:def_norm_R}.
Here, we need to check the assumptions of~\cite[Theorem 3.2]{Bardi1997Bellman}, which include:
\begin{itemize}
    \item[$(H_0)$] The domain $A$ of the control variable $\bu$ is a closed subset of a normed space. 
    \item[$(H_1)$] $f\colon\Rn \times A \to \Rn$ is continuous, and there exists a constant $L>0$, such that
    \begin{equation*}
        \langle f(\bx, \bu) - f(\by, \bu), \bx - \by\rangle \leq L\|\bx - \by\|^2 \quad \forall \bx, \by \in \Rn, \bu\in A.
    \end{equation*}

    \item[$(H_2)$] $\ell\colon \Rn \times A \to \R$ is continuous and bounded from below.
    \item[$(H_3)$] $\tilde{\initcond}\colon\Rn \to \R$ is continuous and bounded from below.
    \item[$(H_4)$] $\lambda\geq 0$. Note that $\lambda$ is the discounting factor in~\cite{Bardi1997Bellman}. In our case, there is no discounting factor, and hence, $\lambda$ is always zero.
    \item[$(H_5)$] There exists $\sigma\geq 1$, such that for any compact $K \subseteq \Rn$, there exists a constant
    $f_K > 0$ satisfying
    \begin{equation*}
        \|f(\bx, \bu)\|\leq f_K(1 + \|\bu\|^\sigma) \quad \forall (\bx, \bu) \in K \times A.
    \end{equation*}

    \item[$(H_6)$] There exist $\ell_0 > 0$, $C_0 \geq 0$, $\delta_1>\sigma$ (where $\sigma$ is the constant in (H5)), such that
    \begin{equation*}
        \ell(\bx,\bu) \geq \ell_0\|\bu\|^{\delta_1} - C_0 \quad \forall (\bx, \bu) \in \Rn \times A.
    \end{equation*}

    \item[(2.1)] There exist $\delta_2\geq 0$ and $\bar{\ell}>0$, such that there holds
    \begin{equation}\label{eqt:lemA10_pf_eqt2p1}
    |\ell(\bx,\bu) - \ell(\by,\bu)| \leq \bar{\ell}\|\bx-\by\|(1 + \|\bu\|^{\delta_1} + \|\bx\|^{\delta_2} + \|\by\|^{\delta_2}) \quad \forall \bx,\by\in\Rn,  \bu\in A.
    \end{equation}
    
    \item[(2.13)] There exists $\delta_2\geq 0$ and $\bar{g}>0$, such that there holds
    \begin{equation}\label{eqt:lemA10_pf_eqt2p13}
        |\tilde{\initcond}(\bx) - \tilde{\initcond}(\by)| \leq \bar{g} \|\bx-\by\|(1 + \|\bx\|^{\delta_2} + \|\by\|^{\delta_2}) \quad \forall \bx,\by\in\Rn,  \bu\in A.
    \end{equation}
\end{itemize}
In our case, the set $A$ equals $\Rn$, the function $f$ does not depend on the state variable $\bx$, and there is no discounting factor in the optimal control problem (i.e., the parameter $\lambda$ in~\cite{Bardi1997Bellman} is zero). Hence, assumptions $(H_0)$, $(H_1)$, and $(H_4)$ hold. By definition of $f$ in~\eqref{eqt:lemA10a_pf_def_ell_f}, assumption $(H_5)$ holds with $\sigma=1$ and $f_K = \max\{1,\|M\balp\|\}$.
Since the potential energy $\potentfn$ is continuous and non-positive, the function $\ell$ is continuous and bounded from below by $- \frac{1}{2}\|\balp\|_M^2$, which implies $(H_2)$.
Moreover, we obtain a lower bound for $\ell$ by
\begin{equation*}
\ell(\bx,\bu) = \frac{1}{2}\|\bu\|_{M^{-1}}^2 - \potentfn(\bx) - \frac{1}{2}\|\balp\|_M^2
\geq \frac{1}{2}\|\bu\|_{M^{-1}}^2  - \frac{1}{2}\|\balp\|_M^2 \geq \ell_0 \|\bu\|^2 - \frac{1}{2}\|\balp\|_M^2,
\end{equation*}
for any $\bx,\bu\in\Rn$, 
where $\ell_0$ is the smallest eigenvalue of the matrix $M^{-1}$, which is positive since $M^{-1}$ is symmetric positive definite. 
Hence, assumption $(H_6)$ holds for $\delta_1=2$ and $C_0 = \frac{1}{2}\|\balp\|_M^2$. Note that we have $\sigma=1$ in $(H_5)$ and $\delta_1 = 2$ in $(H_6)$, and thus, the inequality $\sigma < \delta_1$ in $(H_6)$ is satisfied.
After some calculation, we have
\begin{equation*}
\begin{split}
    |\ell(\bx,\bu)-\ell(\by,\bu)| &= |\potentfn(\bx)-\potentfn(\by)|\leq \bar{\ell}\|\bx-\by\| \quad \forall \bx,\by,\bu\in\Rn,
\end{split}
\end{equation*}
where $\bar{\ell}$ is the Lipschitz constant of $\potentfn$, and hence,~\eqref{eqt:lemA10_pf_eqt2p1} holds. 
Now, it remains to check the assumptions for the initial condition $\tilde{\initcond}$. Recall that $\tilde{\initcond}$ is continuous and bounded from below, and thus, $(H_3)$ is satisfied.
By assumption, $\initcond$ satisfies~\eqref{eqt: HJ1_hd_condJ}. Then, by~\eqref{eqt:lemA10a_pf_def_tildeJ}, we obtain
\begin{equation*}
    |\tilde{\initcond}(\bx)-\tilde{\initcond}(\by)|\leq |\initcond(\bx) - \initcond(\by)| + \|\balp\|\|\bx-\by\|\leq C\|\bx-\by\|(1 + \|\bx\|^\delta + \|\by\|^\delta)+ \|\balp\|\|\bx-\by\|,
\end{equation*}
for any $\bx,\by\in\R^n$.
Hence,~\eqref{eqt:lemA10_pf_eqt2p13} holds with $\bar{g} = C+ \|\balp\|$ and $\delta_2=\delta$.
Therefore, all the assumptions in~\cite[Theorem 3.2]{Bardi1997Bellman} are satisfied. Then, applying~\cite[Theorem 3.2]{Bardi1997Bellman} (with the time reversal technique applied to the optimal control problem) to any time horizon $T>0$, the function $\tilde{\valuefn}$ defined in~\eqref{eqt:lemA10a_pf_newoptctrl} is the unique viscosity solution in the solution set $\tilde{\mathcal{G}}$ in~\eqref{eqt:prop23_pf_def_tildeG} to the HJ PDE whose initial condition is $\tilde{\initcond}$, and the Hamiltonian equals
\begin{equation*}
\begin{split}
    H(\bx,\bp) &= \sup_{\bu\in\Rn} \{\langle f(\bx,\bu), \bp\rangle - \ell(\bx, \bu)\}\\
    &= \sup_{\bu\in\Rn} \left\{\langle \bu+M\balp, \bp\rangle - \frac{1}{2}\|\bu\|_{M^{-1}}^2 + \potentfn(\bx) + \frac{1}{2}\|\balp\|_M^2\right\}\\
    &= \frac{1}{2}\|\bp + \balp\|_M^2  + \potentfn(\bx).
\end{split}
\end{equation*}
Therefore, the function $\tilde{\valuefn}$ is the unique viscosity solution to the HJ PDE~\eqref{eqt: lemA10_newHJ1} in the solution set $\tilde{\mathcal{G}}$.
By straightforward calculation using~\eqref{eqt:lemA10a_pf_eqt_tildeS_S}, we obtain
\begin{equation}\label{eqt:lemA10a_pf_subdiff_superdiff}
    D^- \valuefn(\bx, t) = D^- \tilde{\valuefn}(\bx, t) + \{(\balp,0)\}, \quad D^+ \valuefn(\bx, t) = D^+ \tilde{\valuefn}(\bx, t) + \{(\balp,0)\}, \quad \forall \bx\in\Rn, t> 0,
\end{equation}
where $D^- W$ and $D^+ W$ respectively denote the subdifferential and superdifferential of a continuous function $W$.
By applying~\eqref{eqt:lemA10a_pf_subdiff_superdiff},~\eqref{eqt:lemA10a_pf_def_tildeJ}, and~\eqref{eqt:lemA10a_pf_eqt_tildeS_S} to~\eqref{eqt: lemA10_newHJ1}, we conclude that the function $\valuefn$ is a viscosity solution to the HJ PDE~\eqref{eqt: HJhd_1_general}. 

To prove statement (a), it remains to prove that the function $\valuefn$ is in the solution set $\mathcal{G}$. 
Since the function $\tilde{\valuefn}$ is in the solution set $\tilde{\mathcal{G}}$, the function $\tilde{\valuefn}$ is continuous in $\Rn\times [0,+\infty)$, and hence, $\valuefn$ is also in continuous in $\Rn\times [0,+\infty)$. For all $T>0$, the function $\tilde{\valuefn}$ is bounded from below in $\Rn\times [0,T]$, and thus, the function $(\bx,t)\mapsto \valuefn(\bx,t)-\langle \balp, \bx\rangle$ is bounded from below in $\Rn\times [0,T]$. Moreover, by straightforward calculation, for all $R>0$, there holds
\begin{equation*}
\begin{split}
    \|\valuefn\|_R &= \sup\left\{|\valuefn(\bx, t)| + \|\bp\|\colon (\bx, t)\in B_R(\R^n)\times [0,R], \bp \in D_{\bx}^-\valuefn(\bx, t)\right\}\\
    &= \sup\left\{\left|\tilde{\valuefn}(\bx, t) + \langle \balp, \bx\rangle \right| + \|\bq + \balp\|\colon (\bx, t)\in B_R(\R^n)\times [0,R], \bq \in D_{\bx}^-\tilde{\valuefn}(\bx, t)\right\}\\
    &\leq \sup\left\{\left|\tilde{\valuefn}(\bx, t)\right| + \|\bq\|\colon (\bx, t)\in B_R(\R^n)\times [0,R], \bq \in D_{\bx}^-\tilde{\valuefn}(\bx, t)\right\} + \|\balp\|R + \|\balp\|\\
    &= \|\tilde{\valuefn}\|_R + \|\balp\|R + \|\balp\| < +\infty,
\end{split}
\end{equation*}
where the second equality follows from~\eqref{eqt:lemA10a_pf_eqt_tildeS_S} and~\eqref{eqt:lemA10a_pf_subdiff_superdiff}.
Therefore, the function $\valuefn$ is in the solution set $\mathcal{G}$, and the statement is proved.

(b)
Let $W$ be a viscosity solution to~\eqref{eqt: HJhd_1_general} in the solution set $\mathcal{G}$. By definition of $\mathcal{G}$, there exists $\balp\in\Rn$, such that for all $T>0$, the function $(\bx,t)\mapsto W(\bx,t)-\langle \balp,\bx\rangle$ is bounded from below in $\Rn\times [0,T]$. 
By picking $t=0$, we get that the function $\bx\mapsto \initcond(\bx)-\langle \balp,\bx\rangle =W(\bx,0)-\langle \balp,\bx\rangle$ is bounded from below in $\Rn$. Hence, by assumption, 
the function $(\bx,t)\mapsto \valuefn(\bx,t) - \langle \balp, \bx\rangle $ is also bounded from below in $\Rn\times [0,+\infty)$. 
Define the functions $\tilde{\valuefn}, \tilde{W}\colon \Rn\times [0,+\infty)\to\R$ by
\begin{equation}\label{eqt:prop23_pf_def_tildeS_tildeW}
    \tilde{\valuefn}(\bx,t) = \valuefn(\bx,t)-\langle\balp, \bx\rangle,\quad\quad
    \tilde{W}(\bx,t) = W(\bx,t)-\langle\balp, \bx\rangle,\quad \forall \bx\in\Rn,\, t\geq 0.
\end{equation}
By definition,~\eqref{eqt:lemA10a_pf_subdiff_superdiff} holds for both functions $\valuefn$ and $W$.
Note that $\valuefn$ and $W$ are both viscosity solutions to the HJ PDE~\eqref{eqt: result_HJ1_1d} in the solution set $\mathcal{G}$ (recall that $\valuefn$ was proved to be a viscosity solution in (a)).
By straightforward calculation using~\eqref{eqt:lemA10a_pf_subdiff_superdiff}, the functions $\tilde{\valuefn}$ and $\tilde{W}$ are two viscosity solutions to the HJ PDE~\eqref{eqt: lemA10_newHJ1} in the solution set $\tilde{\mathcal{G}}$. However, by~\cite[Theorem 3.2]{Bardi1997Bellman}, the viscosity solution to~\eqref{eqt: lemA10_newHJ1} is unique in $\tilde{\mathcal{G}}$, which implies $\tilde{\valuefn} = \tilde{W}$, and hence, $\valuefn= W$ holds. Therefore, $\valuefn$ is the unique viscosity solution to the HJ PDE~\eqref{eqt: HJhd_1_general} in the solution set $\mathcal{G}$.
\end{proof}

\section{Some numerical computations}\label{sec: ADMMcomp_HJ1}

\subsection{Proximal point of $p\mapsto -\frac{\valuefn(x,t;p,a,b)}{\lambda}$} \label{sec:appendix_numerical_prox}
In this section, we provide a numerical method for solving the following problem:
\begin{equation}\label{eqt: ADMM_HJ1_di_update_simplified}
    p^* = \argmin_{p\in\mathbb{R}} \left\{-\valuefn(x,t; p,a,b) + \frac{\lambda}{2}(p - c)^2\right\},
\end{equation}
where $x,c\in\R$ and $t,a,b,\lambda>0$ are some parameters and $\valuefn$ is the function defined in~\eqref{eqt: result_S1_1d} and~\eqref{eqt: result_S1_1d_negative}.

By Lemma~\ref{lem: prop_S}, the objective function in \eqref{eqt: ADMM_HJ1_di_update_simplified} is strictly convex, 1-coercive, and continuously differentiable. Therefore, the minimizer exists and is unique. Moreover, $p^*$ is the minimizer if and only if the objective function has zero first-order derivative at $p^*$, i.e., if and only if $p^*$ satisfies
\begin{equation}\label{eqt:appendixB_optcond}
     0 = \left.\frac{\partial (-\valuefn(x,t;p,a,b) + \frac{\lambda}{2} (p - c)^2)}{\partial p}\right|_{p=p^*} = \left.\left(-\frac{\partial \valuefn(x,t;p,a,b)}{\partial p} + \lambda(p - c)\right)\right|_{p=p^*}.
\end{equation}
Recall that the function $V$ is defined piecewise. 
Thus, we compute the first-order derivative using~\eqref{eqt: HJ1_p_deriv} and solve~\eqref{eqt:appendixB_optcond} on each piece. 

First, we consider the case where $p\geq 0$. If $(x,t,p) \in \Omega_1$, then we have
\begin{equation*}
    -\frac{\partial \valuefn(x,t;p,a,b)}{\partial p} + \lambda(p - c) = \frac{at^2}{2} + pt -x  + \lambda(p-c),
\end{equation*}
which has the root
\begin{equation*}
    p_1(x,t,a,b,c) = \frac{-\frac{at^2}{2} +x + \lambda c}{t+\lambda}. 
\end{equation*}
If $(x,t,p) \in \Omega_2$, then we have
\begin{equation*}
    -\frac{\partial \valuefn(x,t;p,a,b)}{\partial p} + \lambda(p - c) = -\frac{bt^2}{2} + pt -x  + \lambda(p-c),
\end{equation*}
which has the root
\begin{equation*}
    p_2(x,t,a,b,c) = \frac{\frac{bt^2}{2} +x+ \lambda c}{t+\lambda}. 
\end{equation*}
If $(x,t,p) \in \Omega_3$, then we have
\begin{equation} \label{eqt:appendix_B_der_3}
\begin{split}
    &-\frac{\partial \valuefn(x,t;p,a,b)}{\partial p} + \lambda(p - c) \\
    =\,& -\frac{a+b}{(a+2b)^2}\left(-(bt-p)^2 + (p-bt)\sqrt{(bt-p)^2 + 2x(a+2b)}\right) \\
    &\quad\quad\quad\quad - \frac{bx}{a+2b} - \frac{bt^2}{2} + pt + \lambda(p-c).
\end{split}
\end{equation}
We apply Newton's method to compute the root $p_3$ of \eqref{eqt:appendix_B_der_3}. In Newton's method, the value is iteratively updated as
\begin{equation*}
\begin{split}
    p_3^{k+1} &= p_3^k - \frac{-\frac{\partial \valuefn(x,t;p,a,b)}{\partial p} + \lambda(p - c)}{-\frac{\partial^2 \valuefn(x,t;p,a,b)}{\partial p^2} + \lambda}
    = p_3^k - \frac{A_1}{A_2},
\end{split}
\end{equation*}
where the numerator $A_1$ equals
\begin{equation*}
\begin{split}
    A_1 = -\frac{a+b}{(a+2b)^2}\left(-(bt-p)^2 + (p-bt)\sqrt{(bt-p)^2 + 2x(a+2b)}\right) \\
    - \frac{bx}{a+2b} - \frac{bt^2}{2} + pt + \lambda(p-c),
\end{split}
\end{equation*}
and the denominator $A_2$ equals 
\begin{equation*}
    A_2 = -\frac{2(a+b)}{(a+2b)^2}\left(bt-p +\frac{(bt-p)^2 + x(a+2b)}{\sqrt{(bt-p)^2 + 2x(a+2b)}} \right) + t + \lambda.
\end{equation*}
To make Newton's method more robust, we also enforce a lower bound $\underline{p}$ defined by
\begin{equation*}
\underline{p} = \max\left\{\frac{x}{t} - \frac{at}{2}, bt-b\sqrt{\frac{2|x|}{a}}, 0\right\}    
\end{equation*}
and an upper bound $\bar{p}$ defined by
\begin{equation*}
\bar{p} = \max\{bt+|c|+1, \underline{p}\}.    
\end{equation*}
The lower bound $\underline{p}$ is given by the definition of $\Omega_3$. The upper bound $\bar{p}$ is set to be $\max\{bt+|c|+1, \underline{p}\}$ since the corresponding function value in~\eqref{eqt:appendix_B_der_3} is positive and the function in~\eqref{eqt:appendix_B_der_3} is increasing with respect to $p$ for $p\geq \underline{p}$, which implies that no root in $\Omega_3$ can be larger than $bt+|c|+1$. If the function value corresponding to $\underline{p}$ in~\eqref{eqt:appendix_B_der_3} is also positive, then there is no root in $\Omega_3$, and we simply set $p_3$ to be $\underline{p}$. Here, the choice of initialization for Newton's method is not crucial. For convenience, in Sections \ref{subsec:numerical_quad} and \ref{sec: ADMM_nonconvex}, we initialize Newton's method with $p_3^0 = \underline{p}+1$, and in Section \ref{sec:numerical_convex}, we initialize Newton's method with the value of $p_3$ found in the previous iteration of ADMM in Algorithm \ref{alg:admm_ver1}.

If $(x,t,p) \in \Omega_4\cup \Omega_5$, we have 
\begin{equation*}
     -\frac{\partial \valuefn(x,t;p,a,b)}{\partial p} + \lambda(p - c) = \frac{p^2}{2b} + \lambda(p-c),
\end{equation*}
whose largest root is given by
\begin{equation*}
    p_4(x,t,a,b,c) = -b\lambda + \sqrt{b^2\lambda^2 + 2b\lambda c},
\end{equation*}
which is non-negative if and only if $c\geq 0$. In other words, if $c$ is negative, $p_4$ is not a candidate for the minimizer.

Now, we consider the case of $p<0$. 
By~\eqref{eqt: result_S1_1d_negative}, we have
\begin{equation} \label{eqt:optcond_negp}
    0 = -\frac{\partial \valuefn(x,t;p,a,b)}{\partial p} + \lambda(p - c)
    = \frac{\partial \valuefn(-x,t;q,b,a) }{\partial q} + \lambda(-q - c),
\end{equation}
where we apply the change of variable $q = -p$. From~\eqref{eqt:optcond_negp}, we observe that the vector $q$ is the positive solution of~\eqref{eqt:appendixB_optcond} where the parameters $(x,t,a,b,c)$ are replaced by $(-x,t,b,a,-c)$. Hence, the roots are in the set $\{p_i'\colon i=1,2,3,4\}$, where each $p'_i$ is defined by
\begin{equation*}
    p'_i(x,t,a,b,c) = -p_i(-x,t,b,a,-c).
\end{equation*}
Therefore, we get several candidates for the minimizer $p^*$. We denote the set of candidates by $C$, which is defined by
\begin{equation*}
    C = \left\{p_i(x,t,a,b,c),\, p'_i(x,t,a,b,c)\colon i=1,2,3,4\right\}.
\end{equation*}

Note that we can simplify the set $C$. 
By straightforward calculation, we obtain
\begin{equation*}
    p'_1(x,t,a,b,c) = p_2(x,t,a,b,c),
    \quad 
    p'_2(x,t,a,b,c) = p_1(x,t,a,b,c).
\end{equation*}
Moreover, we have
\begin{equation*}
    p'_4(x,t,a,b,c) = a\lambda - \sqrt{a^2\lambda^2 -2a\lambda c},
\end{equation*}
which is negative if and only if $c<0$. In other words, if $c$ is non-negative, $p'_4$ is not a candidate for the minimizer. Define $\tilde{p}_4(x,t,a,b,c)$ by 
\begin{equation*}
    \tilde{p}_4(x,t,a,b,c)
    = \begin{dcases}
    p_4(x,t,a,b,c) & c\geq 0,\\
    p'_4(x,t,a,b,c) & c< 0.
    \end{dcases}
\end{equation*}
Then, $\tilde{p}_4(x,t,a,b,c)$ can replace $p_4(x,t,a,b,c)$ and $p'_4(x,t,a,b,c)$ as candidate minimizers.
Similarly, if $x$ is negative, then $(x,t,p)$ cannot be in the set $\Omega_3$, and hence $p_3$ is not a possible candidate. Meanwhile, if $x$ is non-negative, then $(-x,t,-p)$ cannot be in the set $\Omega_3(b,a)$ (which denotes the set $\Omega_3$ with parameters $b,a$, instead of $a,b$), and hence, $p'_3$ is not a possible candidate. Thus, we define $\tilde{p}_3(x,t,a,b,c)$ by 
\begin{equation*}
    \tilde{p}_3(x,t,a,b,c)
    = \begin{dcases}
    p_3(x,t,a,b,c) & x\geq 0,\\
    p'_3(x,t,a,b,c) & x< 0,
    \end{dcases}
\end{equation*}
and we use $\tilde{p}_3$ to replace $p_3$ and $p'_3$ as candidates.
As a result, we simplify the set $C$ to
\begin{equation*}
    C = \{p_1(x,t,a,b,c),\, p_2(x,t,a,b,c),\, \tilde{p}_3(x,t,a,b,c),\, \tilde{p}_4(x,t,a,b,c)\}.
\end{equation*}
Finally, the minimizer $p^*$ is selected among the four candidates as follows:
\begin{equation*}
    p^* = \argmin_{p\in C} \left\{-\valuefn(x,t; p,a,b) + \frac{\lambda}{2}(p - c)^2\right\}.
\end{equation*}

\subsection{An equivalent expression for $\valuefn$}
Let $\valuefn$ be the function defined in~\eqref{eqt: result_S1_1d} and~\eqref{eqt: result_S1_1d_negative}. In this section, we present an equivalent expression for $\valuefn$, which is used in our numerical implementation. 
When $p\geq 0$, the function $\valuefn$ can be written as follows: 
\begin{equation*}
\begin{split}
    \valuefn(x,t; \pzero,a,b) &= \begin{dcases}
    f_1 & (x,t,\pzero)\in\Omega_1,\\
    f_2 & (x,t,\pzero)\in\Omega_2,\\
    \max\{f_3, f_4\} & (x,t,\pzero)\in\Omega_3 \cup \Omega_4,\\
    f_5 & (x,t,\pzero)\in\Omega_5,
    \end{dcases}\\
    &= \begin{dcases}
    f_1 & x\geq pt + \frac{at^2}{2},\\
    f_2 & \left(t<\frac{p}{b}, x<0\right) \text{ or } \left(x < -\frac{b}{2}\left(t-\frac{p}{b}\right)^2, t\geq\frac{p}{b}\right),\\
    \max\{f_3, f_4\} & 0\leq x < pt + \frac{at^2}{2},\\
    f_5 & t-\frac{p}{b} \geq \sqrt{\frac{2|x|}{b}}, x<0, 
    \end{dcases}
\end{split}
\end{equation*}
where $f_1,\dots, f_5$ are the functions defined in~\eqref{eqt: case1_deff}, or equivalently in \eqref{eqt: equiv_form_fi}, below. Here, we use the notation $f_i$ instead of $f_i(x,t; \pzero,a,b)$ for simplicity. 
By straightforward calculation, when $\pzero <0$, the formula reads:
\begin{equation*}
    \valuefn(x,t; \pzero,a,b) = \begin{dcases}
    f_2 & x\leq pt - \frac{bt^2}{2},\\
    f_1 & \left(t<-\frac{p}{a}, x>0\right) \text{ or } \left( x> \frac{a}{2}\left(t+\frac{p}{a}\right)^2, t\geq -\frac{p}{a}\right),\\
    \max\{f'_3, f'_4\} & pt - \frac{bt^2}{2} < x \leq 0,\\
    f'_5 & t+\frac{p}{a}\geq \sqrt{\frac{2|x|}{a}}, x > 0, 
    \end{dcases}
\end{equation*}
where $f'_i$ denotes $f_i(-x,t; -p,b,a)$. The above equivalent expressions for $\valuefn$ are advantageous since they slightly reduce the amount of conditional branching required by our implementation as well as simplify the conditions that need to be checked. Thus, both of these differences help promote the performance of our implementation.
Furthermore, the equivalent formulas for $f_i$ and $f'_i$ in our implementation are given as follows:
\begin{equation}\label{eqt: equiv_form_fi}
\begin{split}
    f_1 &= \frac{p^3}{6a} - \frac{(at+p)^3}{6a} + x(at+p),\\
    f_2 &= -\frac{p^3}{6b} - \frac{(bt-p)^3}{6b} - x(bt-p),\\
    f_3 &= \frac{a+b}{3(a+2b)^2}\left((bt-p)^3 + \sqrt{\Delta}^3\right) - \frac{bx(bt-p)}{a+2b} -\frac{p^3}{6b} - \frac{(bt-p)^3}{6b},\\
    f'_3 &= \frac{a+b}{3(2a+b)^2}\left((at+p)^3 + \sqrt{\Delta'}^3\right) + \frac{ax(at+p)}{2a+b} + \frac{p^3}{6a} - \frac{(at+p)^3}{6a},\\
    f_4 &= \frac{\sqrt{8a|x|^3}}{3} - \frac{p^3}{6b},\\
    f'_4 &= \frac{\sqrt{8b|x|^3}}{3} + \frac{p^3}{6a},\\
    f_5 &= \frac{\sqrt{8b|x|^3}}{3} - \frac{p^3}{6b},\\
    f'_5 &= \frac{\sqrt{8a|x|^3}}{3} + \frac{p^3}{6a},
\end{split}
\end{equation}
where we define $\Delta = (bt-p)^2 + 2x(a+2b)$ and $\Delta' = (at+p)^2 - 2x(2a+b)$. These equivalent expressions are formulated in order to reduce extraneous arithmetic operations in the implementation as well as to avoid any potential complications with undefined square roots.

\subsection{Proof of Proposition~\ref{prop:convergence_ADMM_convexJ}} \label{subsec:appendix_pf_prop31}
Let $\bv^N$, $\bd^N$, and $\bp^N$ be the corresponding vectors in the algorithm at the $N$-th iteration. Define the function $F\colon \Rn\to\R$ by
\begin{equation*}
    F(\bp) = \sum_{i=1}^n \valuefn(x_i, t; p_i,a_i,b_i) \quad \forall \bp=(p_1,\dots,p_n)\in\Rn.
\end{equation*}
Then, the objective function in~\eqref{eqt: result_Hopf1_hd} equals $F-\initcond^*$.
Let $\bp^*=(p^*_1,\dots,p^*_n)$ be the optimizer of the optimization problem in~\eqref{eqt: result_Hopf1_hd}, which exists and is unique by Lemma~\ref{lem: prop_barp_hopf_hd}(a).
By~\cite[Section~3.2]{Boyd2011Distributed}, we have
\begin{equation*}
    \lim_{N\to\infty} \left(F(\bd^N) - \initcond^*(\bv^N)\right) = F(\bp^*) - \initcond^*(\bp^*) = \Sanaly(\bx,t).
\end{equation*}
According to~\cite[Theorem~2.2]{Deng2016global} whose assumptions are proved using~\cite[Remark~2.2]{Deng2016global}, both $\bv^N$ and $\bd^N$ converge to the point $\bp^*$ as $N$ approaches infinity, and hence, $\bp^N$ in Algorithm~\ref{alg:admm_ver1} also converges to $\bp^*$. 
By Lemma~\ref{lem: prop_S}, the function $F$ is continuously differentiable, and hence, it is also Lipschitz in any compact domain. Let $L$ be its Lipschitz constant on a compact domain containing $\{\bd^N\}$ and $\{\bv^N\}$. Then, we have
\begin{equation*}
\begin{split}
\lim_{N\to\infty} |\Snum(\bx,t) - \Sanaly(\bx,t)|
&\leq \lim_{N\to\infty} \left|F(\bd^N) - \initcond^*(\bv^N) -\Sanaly(\bx,t)\right| + \lim_{N\to\infty}\left|F(\bd^N) - F(\bv^N)\right|
\\
&\leq \lim_{N\to\infty} L\left\|\bd^N - \bv^N\right\| = 0.
\end{split}
\end{equation*}

Now, it remains to prove the second formula in~\eqref{eqt:prop31_conv}. By definition of $\opttraj$ in~\eqref{eqt: optctrl_defx_1},~\eqref{eqt: optctrl_defx_2},~\eqref{eqt: optctrl_defx_3}, \eqref{eqt: optctrl_defx_4},~\eqref{eqt: optctrl_defx_5}, and~\eqref{eqt: optctrl_defx_neg}, the function $p\mapsto \opttraj(s;x,t,p,a,b)$ is continuous. Since $\bp^N$ converges to $\bp^*$, $\gmnum$ converges to $\gmanaly$ pointwise for any $s\in[0,t]$. 
Moreover, by straightforward calculation, the function $s\mapsto \opttraj(s;x,t,p,a,b)$ is continuously differentiable with respect to $s$, and its derivative is bounded by $|p| + (a+b)t$. Also, the function $s\mapsto \opttraj(s;x,t,p,a,b)$ is bounded by $|x| + |p|t + (a+b)t^2$. Thus, the second formula in~\eqref{eqt:prop31_conv} holds by the Arzela-Ascoli theorem.
\qed 

\bibliographystyle{spmpsci}      
\bibliography{biblist_new}

\end{document}